\documentclass{amsart}
\usepackage[margin=1.3in]{geometry} %for adjusting the margins
\usepackage{graphicx} % Required for inserting images
\usepackage{amsthm}
\usepackage{amsfonts}
\usepackage{amsmath}
\usepackage{enumerate}
\usepackage{comment}
\usepackage{mathtools}
\usepackage{amssymb}
\usepackage{listings}
\usepackage{amsmath,amssymb}
\usepackage{array}
\usepackage{float}
\usepackage{tikz}
\usetikzlibrary{tikzmark}
\usepackage{hyperref}
\usepackage{cleveref}

\allowdisplaybreaks %allows align to break over pages

\usepackage[backend=biber]{biblatex}
\newtheorem{thm}{Theorem}[section]
\newtheorem{cor}{Corollary}[section]
\newtheorem{prop}{Proposition}[section]
\newtheorem{conj}{Conjecture}[section]

\newtheorem*{thm*}{Theorem}
\newtheorem*{conj*}{Conjecture}
\newtheorem*{prop*}{Proposition}

\theoremstyle{definition}
\newtheorem{rem}{Remark}[section]

\crefname{prop}{proposition}{propositions}
\Crefname{prop}{Proposition}{Propositions}
\crefname{thm}{theorem}{theorems}
\Crefname{thm}{Theorem}{Theorems}

\title{Indecomposable Uniform Families on the Adams 4-Line}
\author{Jordan Benson, Runji Li}
\date{November 2026}

\begin{document}

\begin{abstract}
We show that the $f_j$, $p_j$, $D_3(j)$, and $p_j'$ families on the 4-line of the classical Adams spectral sequence all conform to the Uniform Doomsday Conjecture. We do so by first using Bruner's formula to establish uniform families of differentials on these classes. We then study the Burklund-Xu and Cartan-Eilenberg spectral sequences to ensure that the differential targets are nonzero. The $f_j$ family requires further analysis of possible Adams $d_2$ differentials since the $f_j$ support uniform $d_3$ differentials. As a result of our work, only two indecomposable $\text{Sq}^0$ families remain to be settled in Adams filtration 4.
\end{abstract}

\maketitle

\tableofcontents

\section{Introduction}

A fundamental tool in the study of the 2-primary stable homotopy groups of spheres is the Adams spectral sequence
\[
\text{Ext}_\mathcal{A}(\mathbb{F}_2,\mathbb{F}_2) \Rightarrow \pi_\ast\mathbb{S}_2,
\]
which converges to the 2-completed stable stems. The homological filtration on the $E_2$ page, arising from any projective resolution used to compute these Ext groups, corresponds to the Adams filtration on the stable stems.

Various finiteness properties have been conjectured about the set of classes in any fixed Adams filtration. For example, in 1971, Milgram (\cite{Milgram}) stated the Doomsday Conjecture of Joel Cohen:

\begin{conj}[Doomsday]
In any fixed Adams filtration, only finitely many elements on the Adams $E_2$ page for the 2-completed sphere $\mathbb{S}_2$ survive to homotopy.
\end{conj}

In 1977, Mahowald (\cite{MahowaldEtaJ}) disproved this conjecture by constructing an infinite family $\{\eta_j\}$ of nontrivial homotopy classes detected by $\{h_1h_j\}$ in Adams filtration 2. In response, Minami (\cite{Minami}) stated the New Doomsday Conjecture, which gives a more precise formulation of Cohen's original idea. Rather than addressing entire Adams filtrations, the New Doomsday Conjecture concerns specific families of elements:

\begin{conj}[New Doomsday]
For any family $\{(\text{Sq}^0)^n x\}_{n\geq0}$ on the Adams $E_2$ page for $\mathbb{S}_2$, only finitely many members survive to homotopy.
\end{conj}

Adams (\cite{Adams}) implicitly verified this conjecture in filtration 1 with his work on the Hopf invariant one problem:

\begin{thm}[\cite{Adams}]
For $j\ge 4$, $d_2(h_j) = h_0 h_{j-1}^2 \neq 0$. The elements $h_0,h_1,h_2,h_3$ are permanent cycles, which detect the classes $2,\eta,\nu,\sigma$ in $\pi_*(\mathbb{S}_2)$ of Hopf invariant 1.
\end{thm}

Further empirical evidence for the New Doomsday Conjecture came from the work of Hill-Hopkins-Ravenel (\cite{HHR}) on the Kervaire invariant one problem in the early 2010s. In their work, they showed that the elements $h_j^2$ all support nontrivial differentials for $j\geq 7$, but the explicit differentials are yet unknown. For $0\le j\le 6$, the $h_j^2$ are known to be permanent cycles, for example by \cite{BMT} for $j=4$, \cite{BJM62,Xu62} for $j=5$, and \cite{LWX} for $j=6$. When $j=4$, an explicit manifold representative is constructed by Jones in \cite{Jones}, but for larger $j$, there are no known manifolds corresponding to $h_j^2$. Since the only $\text{Sq}^0$ families in Adams filtrations 1 and 2 are of the form $h_j$ or $h_j^2$, this confirms that the New Doomsday Conjecture holds in these filtrations.

Beyond these cases, the plausibility of the New Doomsday Conjecture is supported by Bruner's formula describing the interaction between Adams differentials and algebraic Steenrod operations:

\begin{thm}[\cite{Bruner}]
Let $x\in E_r^{s,n}$ be an element on the Adams $E_r$ page for the sphere. Then
\[
d_\ast \text{Sq}^i x = \text{Sq}^{i+r-1}d_rx \dotplus
\begin{cases}
0 & v > k+1 \text{ or } 2r-2 < v < k \\
axd_rx & v = k+1 \\
a\text{Sq}^{i+v}x & v = k \text{ or } (v < k \text{ and } v\leq10)
\end{cases}
\]
%where $k=s-i$, $v=v(n+k)$, and $a\in E_\infty^{\phi,v-1}$ detects a generator of $\text{im }J$ in $\pi_{v-1}S^0$. Here, $v(\ell) = 8a+2^b$ if the exponent of 2 in the prime factorization of $\ell+1$ is $4a+b$ with $0\leq b\leq 3$. The notation $a\dotplus b$ should be interpreted as picking out the lowest-filtration part of $a+b$.
\end{thm}

The notation in the statement of this formula will be explained in \Cref{thm:bruner-formula} below. If one considers classes in low filtration, Bruner's formula suggests that elements belonging to a common $\text{Sq}^0$ family support similar differentials. Based on this observation, the result of Adams \cite{Adams}, and the result of Hill-Hopkins-Ravenel \cite{HHR}, Xu and his collaborators formulated the following refinement of the New Doomsday Conjecture, called the Uniform Doomsday Conjecture (see \cite{ICM}, \cite{BX}):

\begin{conj}[Uniform Doomsday Conjecture]
Let $\{a_j\}_{j\geq0}$ be a nontrivial $\text{Sq}^0$ family on the Adams $E_2$ page. Then there exists another $\text{Sq}^0$ family $\{b_j\}_{j\geq0}$, an $r\geq2$, and a class $c$ such that
\[
d_r(a_j) = c\cdot b_j \neq 0
\]
for $j>>0$.
\end{conj}

In other words, not only should $\text{Sq}^0$ families eventually die before reaching the Adams $E_\infty$ page, as in the New Doomsday Conjecture, but they should do so in a uniform way. As further evidence of this Uniform Doomsday Conjecture, Burklund-Xu (\cite{BX}) deduce the following uniform family of differentials:
\begin{thm}[\cite{BX}]
There are differentials
\[
d_4(h_j^3) = h_0^3g_{j-2} \neq 0
\]
for $j\geq6$.
\end{thm}

In this work, we settle the fates of four $\text{Sq}^0$ families on the 4-line of the classical 2-primary Adams spectral sequence and show that they all conform to the Uniform Doomsday Conjecture of Burklund-Xu:
\begin{thm}[Main Theorem]\label{thm:main}
The following differentials are nonzero for all $j\geq0$:
\[
\begin{split}
    d_3(f_{j+2}) &= h_1h_{j+4}Q_3(j) \\ %tight bound
    d_2(p_{j+4}) &= h_0h_{j+3}p_{j+3}' \\%tight bound
    d_2(D_3(j+2)) &= h_0K_{j+1} \\ %tight bound
    d_2(p_{j+4}') &= h_0T_{j+3} %maybe p_2' or p_3' works
\end{split}
\]
\end{thm}

\begin{rem}
By the finite-stem calculations in \cite{charts} and \cite{Weinan}, the first three families of differentials are maximal, in the sense that no earlier members of those families support the stated differential. For the final family, it is not known whether $p_2'$ or $p_3'$ supports the stated differential.
\end{rem}

By \cite{Chen}, the indecomposable elements of filtration 4 on the Adams $E_2$ page can be grouped into seven $\text{Sq}^0$ families:
\[
\begin{array}{c c c c c}
d_j &&e_j &&f_j \\ \\
g_{j+1} &&p_j &&D_3(j) \\ \\
&&p_j' &&
\end{array}
\]
An immediate result of our work is that more than half the indecomposables here will have been shown to conform to the Uniform Doomsday Conjecture of Burklund-Xu. Of the remaining families, the $e_j$ family has been settled by \cite{Li-Li}. Then all that remain are the following two families:
\[
\begin{array}{c c c}
d_j &&g_{j+1}
\end{array}
\]

The fate of the $d_j$ family is settled in joint work in progress with Yueshi Hou, Yuxuan Li, and Zhouli Xu. At present, there is no indication that the remaining $g_{j+1}$ family can be approached in a similar way.

Because Bruner's formula provides many predictions of uniform families of Adams differentials, the main obstacle to establishing further families has been the computation of higher filtrations of the Adams $E_2$ page. Computing these higher filtrations becomes prohibitively difficult very quickly, with the computation of the Adams 5-line in \cite{Chen} requiring hundreds of families of differentials in the Lambda algebra. We avoid this difficulty by using the Burklund-Xu and Cartan-Eilenberg spectral sequences to target specific bidegrees in Adams filtrations 6 and 7, where practically none of the multiplicative or even additive structure of the Adams $E_2$ page is known.

We end the introduction by giving a possible interpretation of our results. The arguments in \Cref{sec:differentials} shed light on large-scale patterns in the Adams spectral sequence. For example, if one applies Bruner's formula to $\text{Sq}^0$ families, one finds that such families \textit{should} support uniform differentials of a very particular form, barring algebraic obstructions. In particular, we expect that $\text{Sq}^0$ families eventually support multiplication by a fixed class:

\begin{conj}
Let $\{a_j\}_{j\geq0}$ be a nontrivial $\text{Sq}^0$ family on the Adams $E_2$ page and $b$ be a fixed nonzero class on the Adams $E_2$ page. Then
\[
b\cdot a_j \neq 0
\]
for $j>>0$.
\end{conj}

This property is known to hold through the Adams 5-line (see Theorem 1.2 of \cite{Chen}), where every relation fits into a $\text{Sq}^0$ family. For example, we have $h_0h_j = 0$ if and only if $j=1$. We believe that this pattern holds in higher filtrations and gives rise to uniform families
\[
d_2(a_{j+1}) = h_0\text{Sq}^1(a_j) \neq 0 \text{, }j>>0
\]
of $d_2$ differentials, so long as $\text{Sq}^1(a_j)$ itself is nonzero for all $j$. If one studies closely the case of the $f_j$ family, one might formulate a similar conjecture for $d_3$ differentials. However, this conjecture would necessarily be much more complicated than for $d_2$ differentials since shorter differentials could intervene.

%In all, empirical evidence in low filtrations suggests that the Uniform Doomsday Conjecture, which predicts that only finitely many members of any $\text{Sq}^0$ family can survive to homotopy, should be mirrored by the existence of only finitely many $\text{Sq}^0$ families of \textit{relations} in any filtration of the Adams spectral sequence.

This paper is structured as follows. After reviewing the Burklund-Xu and Cartan-Eilenberg spectral sequences in \Cref{sec:recollections}, we use Bruner's formula in \Cref{sec:differentials} to derive formulas for the four uniform families of differentials in \Cref{thm:main}. The remaining four sections handle each family individually to show that the targets of the differentials are nonzero on the Adams $E_2$ or $E_3$ page. In the appendix we provide the commands needed to use the software of Bruner and Chatham-Chua-Beauvais-Feisthauer.

\section{Acknowledgments}
The authors would like to thank their advisor, Zhouli Xu, for suggesting this project and for careful feedback on early drafts of this paper.

\section{Recollections on the Burklund-Xu Method}\label{sec:recollections}
In this section we briefly review a method used in \cite{WangXu61} and \cite{BX} (see also \cite{Benson}, \cite{Li-Li} for detailed expositions) to determine whether certain classes survive to the classical Adams spectral sequence beyond the 5-line, where our knowledge of entire Adams filtrations ends. First recall that the classical Adams spectral sequence is of the form
\[
\text{Ext}_\mathcal{A}(\mathbb{F}_2,\mathbb{F}_2) \Longrightarrow (\pi_\ast S^0)_2,
\]
where the subscript denotes 2-completion and $\mathcal{A} = \mathbb{F}_2[\xi_1,\xi_2,\dots]$ denotes the dual mod-2 Steenrod algebra. In order to compute the $E_2$ page, one can study the Cartan-Eilenberg spectral sequence associated to the Hopf algebra extension
\[
\mathcal{P}\to\mathcal{A}\to\mathcal{Q},
\]
where $\mathcal{P} = \mathbb{F}_2[\xi_1^2,\xi_2^2,\dots]$ denotes the Hopf subalgebra of $\mathcal{A}$ generated by the squares of the $\xi_i$ and $\mathcal{Q}$ denotes the exterior algebra on the $\xi_i$. This Cartan-Eilenberg spectral sequence is of the form
\[
\text{Ext}_\mathcal{P}(\mathbb{F}_2,\text{Ext}_\mathcal{Q}(\mathbb{F}_2,\mathbb{F}_2)) \Longrightarrow \text{Ext}_\mathcal{A}(\mathbb{F}_2,\mathbb{F}_2),
\]
and since $\mathcal{Q}$ is an exterior algebra, the $E_2$ page can be identified with
\[
\text{Ext}_\mathcal{P}(\mathbb{F}_2,\mathbb{F}_2[q_0,q_1,\dots]),
\]
where $q_i$ corresponds to $\xi_{i+1}\in\mathcal{Q}$. This $E_2$ page is simpler to compute than the full Adams $E_2$ page, in part because there is a splitting
\[
\mathbb{F}_2[q_0,q_1,\dots] \cong \mathbb{F}_2\oplus\bigoplus\limits_i\mathbb{F}_2\langle q_i\rangle\oplus\bigoplus\limits_{i,j}\mathbb{F}_2\langle q_iq_j\rangle\oplus\dots
\]
as $\mathcal{P}$-comodules. This allows one to compute the Cartan-Eilenberg $E_2$ page at each fixed ``length," the number of $q$ factors. Because of this splitting, one may further filter by the Burklund-Xu degree
\[
\text{deg} : \left(a\in \text{Ext}_\mathcal{P}(\mathbb{F}_2,\mathbb{F}_2\langle q_{i_1}q_{i_2}\dots q_{i_k}\rangle)\right) \mapsto i_1+i_2+\dots+i_k
\]
to obtain the following spectral sequence:
\[
E_1 = \mathbb{F}_2[q_0,q_1,\dots]\otimes\text{Ext}_\mathcal{P}(\mathbb{F}_2,\mathbb{F}_2) \Longrightarrow \text{Ext}_\mathcal{P}(\mathbb{F}_2,\mathbb{F}_2[q_0,q_1,\dots])
\]
We will call this spectral sequence the \textit{Burklund-Xu spectral sequence}. The differentials in this spectral sequence are of the form
\[
d_r: E_r^{s,k,t,\text{deg}} \to E_r^{s+1,k,t,\text{deg}-r},
\]
where the individual degrees are defined as follows. If one has a class
\[
q_{i_1}q_{i_2}\dots q_{i_k}\cdot a \in E_1,
\]
then $s$ is the homological degree of $a\in\text{Ext}^s_\mathcal{P}(\mathbb{F}_2,\mathbb{F}_2)$. The $k$ degree is the Cartan-Eilenberg filtration, the number of $q$ factors. The $t$ degree is given by
\[
t(q_{i_1}q_{i_2}\dots q_{i_k}\cdot a) = 2^{i_1+1}-1 + 2^{i_2+1}-1 + \dots + 2^{i_k+1}-1 + t(a),
\]
where $t(a)$ is the internal degree of $a\in\text{Ext}_\mathcal{P}(\mathbb{F}_2,\mathbb{F}_2)$. Note that we are considering the internal degree of $a$ as an element in the cohomology of $\mathcal{P}$, rather than $\mathcal{A}$, so this internal degree is twice the expected degree. For example, if $a = h_0\in\text{Ext}_\mathcal{P}(\mathbb{F}_2,\mathbb{F}_2)$, then $t(a) = 2^{0+1} = 2$. This is because $h_0$ is identified with $h_1$ under the doubling-up isomorphism
\[
\text{Ext}_\mathcal{P}^{s,t}(\mathbb{F}_2,\mathbb{F}_2) \cong \text{Ext}_\mathcal{A}^{s,2t}(\mathbb{F}_2,\mathbb{F}_2).
\]
Finally, the Burklund-Xu degree is, as described above in setting up the spectral sequence, given by
\[
\text{deg}(q_{i_1}q_{i_2}\dots q_{i_k}\cdot a) = i_1+i_2+\dots+i_k.
\]

Some low-page differentials are computed in \cite{BX}, and we now summarize the most important differentials for our work:
\begin{enumerate}
\item $d_1(q_{i+1}\cdot a) = q_i\cdot h_i a$
\item $d_2(q_{i+2}\cdot a) = q_i\cdot \langle h_i,h_{i+1},a\rangle$
\item $d_2(q_{i+1}^2\cdot a) = q_i^2\cdot h_{i+1}a$
\end{enumerate}

In order to apply these differentials to classes in Cartan-Eilenberg filtration at least 2, we may use the following result of \cite{Benson}:
\begin{prop}\label{prop:leibniz}
Burklund-Xu differentials follow the Leibniz rule. In particular, if $q_{i_1}\dots q_{i_{k'}}\cdot a$ and $q_{i_{k'+1}}\dots q_{i_k}\cdot a$ are both $E_{r-1}$-cycles, then
\[
d_r(q_{i_1}\dots q_{i_{k'}}q_{i_{k'+1}}\dots q_{i_k}\cdot a) = \left(d_r(q_{i_1}\dots q_{i_{k'}})\right)(q_{i_{k'+1}}\dots q_{i_k})\cdot a + (q_{i_1}\dots q_{i_{k'}})\left(d_r(q_{i_{k'+1}}\dots q_{i_k})\right)\cdot a
\]
\end{prop}

As an example of how one might go about computing such differentials, consider the element
\[
q_1q_2^2\cdot h_2^2
\]
on the Burklund-Xu $E_1$ page. We may apply the Leibniz rule to compute the $d_1$ differential:
\[
\begin{split}
d_1(q_1q_2^2\cdot h_2^2) &= \left(d_1(q_1)\right)(q_2^2)\cdot h_2^2 + (q_1)\left(d_1(q_2)\right)(q_2)\cdot h_2^2 + (q_1q_2)\left(d_1(q_2)\right)\cdot h_2^2 \\
&= q_0q_2^2\cdot h_0h_2^2 + q_1^2q_2\cdot h_1h_2^2 + q_1^2q_2\cdot h_1h_2^2 \\
&= 0
\end{split}
\]

Now we note that, in computing this $d_1$ differential, we found that both $q_1\cdot h_2^2$ and $q_2^2\cdot h_2^2$ were individually $d_1$ cycles. As a result, we may now compute $d_2$ using the Leibniz rule:
\[
\begin{split}
d_2(q_1q_2^2\cdot h_2^2) &= \left(d_2(q_1)\right)(q_2^2)\cdot h_2^2 + (q_1)\left(d_2(q_2)\right)\cdot h_2^2 \\
&= q_{-1}q_2^2\cdot\langle h_{-1},h_0,h_2^2\rangle + q_1^3\cdot h_2^3 \\
&= q_1^3\cdot h_2^3
\end{split}
\]
Here, we are using the fact that 0 is the lowest possible index of $q$ and $h$, forcing the first term to vanish.

\begin{rem}
While this is a valid expression for this $d_2$ differential, the attentive reader may notice that the target is trivial on the Burklund-Xu $E_2$ page:
\[
d_1(q_1^2q_2\cdot h_1h_3) = q_1^3\cdot h_1^2h_3 = q_1^3\cdot h_2^3
\]
As a result, the class $q_1q_2^2\cdot h_2^2$ survives to $E_3$, where it may or may not be trivial.
\end{rem}

In order to show that a class in the Cartan-Eilenberg spectral sequence is not killed by differentials, we can rule out all possible sources by showing that they support or receive Burklund-Xu differentials. Since Cartan-Eilenberg differentials preserve internal degree, we can enumerate all possible sources by searching for classes with a fixed internal degree.

\section{Obtaining the Differentials}\label{sec:differentials}

In this section, we apply Bruner's formula to the $\text{Sq}^0$ families of interest. We obtain expressions for certain Adams differentials on classes $p_j$, $D_3(j)$, $p_j'$, and $f_j$. This will establish the differentials claimed in \Cref{thm:main}, and the remaining sections will be devoted to showing that the differentials are in fact nontrivial. Recall that Bruner's formula is given as follows:

\begin{thm}[\cite{Bruner}]\label{thm:bruner-formula}
Let $x\in E_r^{s,t-s}$ be an element on the Adams $E_r$ page for the sphere. Then
\[
d_\ast \text{Sq}^i x = \text{Sq}^{i+r-1}d_rx \dotplus
\begin{cases}
0 & v > k+1 \text{ or } 2r-2 < v < k \\
axd_rx & v = k+1 \\
a\text{Sq}^{i+v}x & v = k \text{ or } (v < k \text{ and } v\leq10)
\end{cases}
\]
where $k=s-i$, $v=v(t-i)$, and $a\in E_\infty^{\phi,v-1}$ detects a generator of $\text{im }J$ in $\pi_{v-1}S^0$. Here, $v(\ell) = 8a+2^b$ if the exponent of 2 in the prime factorization of $\ell+1$ is $4a+b$ with $0\leq b\leq 3$. The notation $a\dotplus b$ should be interpreted as picking out the lowest-filtration part of $a+b$.
\end{thm}

We first prove a result which will be useful for all four families:

\begin{prop}\label{prop:d2-pattern}
If $x$ is an element defined on the Adams $E_2$ page with internal degree $t(x)$ divisible by 2 and Adams filtration at least 1, then
\[
d_2(\text{Sq}^0(x)) = h_0\text{Sq}^1(x).
\]
\end{prop}
\begin{proof}
By Bruner's formula, the first possible contribution to a possible differential on $\text{Sq}^0(x)$ is
\[
\text{Sq}^{i+r-1}d_rx = \text{Sq}^{0+2-1}d_2x = \text{Sq}^1d_2x
\]
since we are considering the case $i=0,r=2$. This contribution lies in Adams filtration
\[
s(\text{Sq}^1d_2x) = 1+2+s(x) = 3+s(x)
\]
so contributes to a $d_3$ differential if it is nonzero at all.

For the second possible contribution to a differential on $\text{Sq}^0(x)$, we must compute $v(t-i) = v(t)$. To do so, note that the degree $t = t(x)$ is divisible by 2, by assumption. Then $t+1$ is odd, so its 2-adic valuation is 0. As a result, we have
\[
v(t) = 8\cdot0+2^0 = 1.
\]
Since $k = s - 0 = s$ is assumed to be at least 1, this gives either $v = k$ or $v < k$ with $v\leq 10$. In either case, the second possible contribution to the differential on $\text{Sq}^0(x)$ is
\[
a\text{Sq}^{1+0}(x) = a\text{Sq}^1(x)
\]
with $a$ detecting a generator of the image of $J$ in $\pi_{1-1}S^0 = \pi_0S^0$. Any such generator is detected by $h_0$, so this term becomes
\[
h_0\text{Sq}^1(x).
\]
Note that this element lies in Adams filtration
\[
s(h_0\text{Sq}^1(x)) = 2+s(x)
\]
so lies in strictly lower Adams filtration than the first possible contribution
\[
\text{Sq}^1d_2x.
\]
As a result, Bruner's formula gives
\[
d_2(\text{Sq}^0(x)) = h_0\text{Sq}^1(x),
\]
as claimed.
\end{proof}

\begin{cor}\label{cor:d2-pattern}
Given a family $\{x_j\}_{j\geq0}$ on the Adams $E_2$ page with
\[
x_j = (\text{Sq}^0)^j x_0,
\]
the formula
\[
d_2(x_{j+1}) = h_0\text{Sq}^1(x_j)
\]
holds for all $j\geq1$.
\end{cor}
\begin{proof}
\Cref{prop:d2-pattern} applies to all elements in filtration at least 1 and degree divisible by 2. Note that every element $x\in E_2$ in the image of $\text{Sq}^0$ has degree divisible by 2 since $\text{Sq}^0$ acts by squaring each component in the corresponding cobar representative. Then \Cref{prop:d2-pattern} applies to every member of $\{x_j\}_{j\geq1}$ in Adams filtration at least 1. Since the only element on the 0-line has trivial $\text{Sq}^i$ for all $i\geq0$, the formula holds vacuously there, hence the lack of a restriction on Adams filtration.
\end{proof}

The following $d_2$ differentials are now read off from Example 1.23(6),(8)-(10) of \cite{BX}, all valid for $j\geq0$:
\[
\begin{split}
d_2(f_{j+2}) &= 0 \\
d_2(p_{j+2}) &= h_0h_{j+1}p_{j+1}' \\
d_2(D_3(j+2)) &= h_0K_{j+1} \\
d_2(p_{j+2}') &= h_0T_{j+1}
\end{split}
\]

Three of these differentials provide the differentials claimed in \Cref{thm:main}, and the only differential remaining to be proven is
\[
d_3(f_{j+2}) = h_1h_{j+4}Q_3(j)
\]
for $j\geq0$. This will be accomplished by the following two propositions.

\begin{prop}
There is an Adams $d_3$ differential
\[
d_3(f_{j+2}) = h_1y_{j+1},
\]
which may or may not be trivial, for every $j\geq0$.
\end{prop}
\begin{proof}
Bruner's formula gives
\[
\begin{split}
d_\ast(\text{Sq}^1c_{j+1}) &= \text{Sq}^{1+2-1}d_2c_{j+1} \\
&= \text{Sq}^2(h_0f_j)
\end{split}
\]
since $c_{j+1}$ has internal degree divisible by precisely $j+1$ factors of 2. Using the Cartan formula for algebraic Steenrod operations, this yields
\[
\begin{split}
d_\ast(\text{Sq}^1(c_{j+1})) &= \text{Sq}^2(h_0)\text{Sq}^0(f_j)+\text{Sq}^1(h_0)\text{Sq}^1(f_j)+\text{Sq}^0(h_0)\text{Sq}^2(f_j) \\
&= 0 + h_1^2\text{Sq}^1(f_j) + h_1\text{Sq}^2(f_j)
\end{split}
\]
for the differential. Using the algebraic Steenrod operations computed in \cite{Bruner-Nassau-Tilson}, we may identify $\text{Sq}^1(c_{j+1})$ with $f_{j+1}$, $\text{Sq}^1(f_j)$ with 0, and $\text{Sq}^2(f_j)$ with an element called $y_j$ to obtain
\[
d_3(f_{j+1}) = h_1y_j
\]
for all $j\geq1$. The first two members $f_0$ and $f_1$ of the $f_j$ family are known, respectively, to support a $d_2$ with target $h_0^2e_0$ and to be a permanent cycle, as shown in the Lin-Wang-Xu charts (\cite{LWX}, \cite{Weinan}). Here we have changed $d_\ast$ to $d_3$ since this lands in the correct filtration of the target.
\end{proof}

We now use some automated tools to identify the $y_{j+1}$ family with a more familiar one.

\begin{prop}\label{prop:identifysq0y0}
The element $y_{j+1} = (\text{Sq}^0)^{j+1}(y_0)$ can be identified with $h_{j+4}Q_3(j)$ for all $j\geq0$.
\end{prop}
\begin{proof}
First note that $y_0$ is the name given in \cite{Bruner-Nassau-Tilson} to the element $6_{16}$ in the Adams resolution for the sphere, produced by Bob Bruner's Ext resolver (\cite{BrunerExt}). Technically, it is the element \textit{dual to} $6_{16}$, but we suppress this terminology here for ease of reading.

In order to identify $(\text{Sq}^0)^j(y_0)$ with a more common name, we will use both Bruner's resolver and the SSeq software developed and maintained by Chatham, Chua, and Beauvais-Feisthauer (\cite{BFCC24}). In \Cref{sec:appendix}, we provide the commands needed to reproduce our results and give more detailed arguments alongside the commands.

%It now remains to confirm that $6_{63}$ is dual to $h_4Q_3(0)$. To this end, we look to stem 36, where there is a class $5_{13}$ in Bruner's notation. In this bidegree, there are classes $h_0^4h_5$, $n_0$, and $h_0^4h_5+n_0$. Of these, only $n_0$ does not support multiplication by $h_0$. Since $5_{13}$ does not support multiplication by $h_0$, we conclude that $5_{13}$ corresponds to $n_0$. Translating $5_{13}$ into Chatham-Chua-Beauvais-Feisthauer notation, we obtain $x_{31,5,1}$, which maps to $x_{67,5,1}$ under $\text{Sq}^0$. Translating this back into Bruner's notation, we obtain the class $5_{38}$. <-- this is overkill

Using both the Bruner and BFCC software, one finds that $\text{Sq}^0(y_0)$ is given by the element $6_{63}$ produced by Bruner's resolver. It then remains to confirm that $6_{63}$ is precisely the element more commonly called $h_4Q_3(0)$. To this end, we study the $(\text{stem},\text{filtration})$-bidegree $(67,5)$ to determine the name of $Q_3(0)$ given by Bruner's resolver.

The bidegree $(67,5)$ contains precisely three nonzero classes, by the charts produced by Isaksen-Wang-Xu (\cite{charts}). Namely, we have $Q_3(0)$, $n_1$, and their sum $Q_3(0)+n_1$. Of these three classes, only $Q_3(0)$ does not support multiplication by $h_3$ (\cite{charts}). Using Bruner's program, we find that the class $5_{39}$ does not support multiplication by $h_3$, so it must correspond to $Q_3(0)$. Bruner's program also shows that multiplying $5_{39}$ by $h_4$ produces the class $6_{63}$, so $6_{63}$ does indeed correspond to $h_4Q_3(0)$.

Then we have shown that $\text{Sq}^0(y_0) = h_4Q_3(0)$, so we have
\[
\begin{split}
d_3(f_{j+2}) &= h_1y_{j+1} \\
&= h_1(\text{Sq}^0)^{j+1}(y_0) \\
&= h_1(\text{Sq}^0)^j(h_4Q_3(0)) \\
&= h_1h_{j+4}Q_3(j)
\end{split}
\]
for $j\geq 0$, as claimed.
\end{proof}

\section{The $p_j$ Family}

We prove the following theorem:
\begin{thm}
There are Adams differentials
\[
d_2(p_{j+1}) = h_0h_jp_j' \neq 0
\]
for $j\geq 3$.
\end{thm}

Recall that the form of the differential was shown to hold in \Cref{sec:differentials}, so all that remains is to show that the target is nonzero on the Adams $E_2$ page. The proof proceeds in two steps, following the program of Burklund-Xu. First we show, for $j\geq2$, that the element $q_0\cdot h_jp_j'$ is a permanent cycle in the algebraic Burklund-Xu spectral sequence and detects a nonzero class in the abutment, an associated graded of the Cartan-Eilenberg $E_2$ page. We then show that $q_0\cdot h_jp_j'$ is a permanent cycle in that spectral sequence and detects a nonzero class in the abutment, an associated graded of the Adams $E_2$ page. Since $h_jp_j'$ detects the classical Adams element $h_{j+1}p_{j+1}'$ under the doubling-up isomorphism
\[
\text{Ext}_\mathcal{A}\xrightarrow{\cong}\text{Ext}_\mathcal{P},
\]
\cite{Benson} shows that $q_0\cdot h_jp_j'$ must detect $h_0h_{j+1}p_{j+1}'$. As a result, we will have shown that $h_0h_jp_j'$ is nonzero on the Adams $E_2$ page for $j\geq3$.

\subsection{Burklund-Xu Differentials}

We first show the following:
\begin{prop}
The Burklund-Xu element $q_0\cdot h_jp_j'$ is a permanent cycle for $j\geq1$.
\end{prop}
\begin{proof}
The Burklund-Xu degree of $q_0\cdot h_jp_j'$ is 0. Since Burklund-Xu differentials strictly lower filtration, this class cannot support any differential in that spectral sequence.
\end{proof}

Then the class $q_0\cdot h_jp_j'$ is defined in the abutment of the Burklund-Xu spectral sequence, so it remains to show that it cannot receive any differential. We use the following simple proposition to simplify this argument.

\begin{prop}\label{prop:deg-pj}
If $y\in\text{Ext}_\mathcal{P}^4$ has internal degree $t(y)$ divisible by $2^2$, then any Burklund-Xu element $q_a\cdot y$ in the same internal degree as $q_0\cdot h_jp_j'$ must not support a differential.
\end{prop}
\begin{proof}
If $q_a\cdot y$ lies in the same internal degree as $q_0\cdot h_jp_j'$, then the equality
\[
t(q_a\cdot y) = 2^{a+1}-1+t(y) = 1+2^{j+2}+2^{j+4}+2^{j+7} = t(q_0\cdot h_jp_j')
\]
must hold. Adding 1 to both sides and dividing by 2 gives
\[
2^a+2^{-1}t(y) = 1+2^{j+1}+2^{j+3}+2^{j+6},
\]
and, by assumption, the term $2^{-1}t(y)$ is still divisible by 2. Then the only possibility is that $a=0$, meaning that
\[
q_a\cdot y = q_0\cdot y
\]
is the source of the differential. Again, classes of this form cannot support Burklund-Xu differentials since they already live in the lowest possible filtration in that spectral sequence. Then no such differential can exist.
\end{proof}

We now show that $q_0\cdot h_jp_j'$ cannot be the target of any Burklund-Xu differential:
\begin{prop}
There is no Burklund-Xu differential of the form
\[
d_r(x) = q_0\cdot h_jp_j'
\]
for any $j\geq1$ and Burklund-Xu element $x$.
\end{prop}
\begin{proof}
Recall that we can write $x$ as a linear combination of elements of the form $q_a\cdot y$ for $y\in\text{Ext}^4_\mathcal{P}$. As a result, we first collect all such classes having the same internal degree as $q_0\cdot h_jp_j'$ and then show that no linear combination of these classes can support a differential as in the proposition statement. We begin with the case where $y$ is indecomposable (see \cite{Chen}, \cite{Lin} for a full presentation of the classical Adams 4-line). Of the possible classes, we rule out $y\in\{d_b,f_b,g_b\}$ using \autoref{prop:deg-pj}.

We now consider $y = e_b$ so that the equality
\[
t(q_a\cdot e_b) = 2^{a+1}-1+2^{b+1}+2^{b+3}+2^{b+5} = 1+2^{j+2}+2^{j+4}+2^{j+7} = t(q_0\cdot h_jp_j')
\]
holds. Adding 1 to both sides and dividing by 2 gives
\[
2^a+2^b+2^{b+2}+2^{b+4} = 1+2^{j+1}+2^{j+3}+2^{j+6},
\]
and since we must avoid $a = 0$ by the reasoning above, we conclude that $b=0$. After subtracting 1 from both sides, this gives
\[
2^a+2^2+2^4 = 2^{j+1}+2^{j+3}+2^{j+6},
\]
but there is no choice of $a$ and $j$ compatible with this relation. Then $y\neq e_b$.

We now consider $y = p_b$ so that the equality
\[
2^a+2^b+2^{b+2}+2^{b+5} = 2^{j+1}+2^{j+3}+2^{j+6}
\]
holds, where we have repeated the same first step of adding 1 and dividing by 2. Again we must have $b=0$. Subtracting 1 from both sides, we obtain
\[
2^a+2^2+2^5 = 2^{j+1}+2^{j+3}+2^{j+6},
\]
but there is no choice of $a$ and $j$ compatible with this relation. Then $y\neq p_b$.

The same logic applies to $y = D_3(b)$, which has internal degree $2^{b+1}+2^{b+7}$, so we now focus on the final indecomposable $y = p_b'$ so that the equality
\[
2^a+2^b+2^{b+3}+2^{b+6} = 1+2^{j+1}+2^{j+3}+2^{j+6}
\]
holds. We must have $b=0$. Then
\[
2^a+2^6+2^3 = 2^{j+1}+2^{j+3}+2^{j+6},
\]
but this holds only if $j=0$, which is below the range being considered. Then $y\neq p_b'$, and we have ruled out all indecomposable $y$.

The decomposable $y$ are of the form $h_\ell c_m$ and $h_bh_ch_dh_e$. We first consider $y = h_\ell c_m$ so that the equality
\[
2^a+2^\ell+2^m+2^{m+1}+2^{m+3} = 1+2^{j+1}+2^{j+3}+2^{j+6}
\]
holds after adding 1 and dividing by 2. We must have $\ell = 0$ or $m = 0$. In the first case, we obtain
\[
2^a+2^{m+3}+2^{m+1}+2^m = 2^{j+1}+2^{j+3}+2^{j+6}
\]
after subtracting 1. There is no compatible choice of $a,m,j$, so we rule out this possibility. If $m=0$, then we obtain
\[
2^a+2^\ell+2^1+2^3 = 2^{j+1}+2^{j+3}+2^{j+6}
\]
after subtracting 1. This forces $a = \ell = j+5 = 5$, so
\[
q_a\cdot h_\ell c_m = q_5\cdot h_5c_0.
\]
However, this class is precisely $d_1(q_6\cdot c_0)$ so cannot support any differentials.

The only remaining classes $y\in\text{Ext}^4_\mathcal{P}$ are of the form $y = h_bh_ch_dh_e$ so that the equality
\[
2^a+2^b+2^c+2^d+2^e = 1+2^{j+1}+2^{j+3}+2^{j+6}
\]
holds after adding 1 and dividing by 2. Once again, we must have $a\neq 0$. Without loss of generality, we can take
\[
b\leq c\leq d\leq e,
\]
so we arrive at $b=0$. This gives
\[
2^a+2^c+2^d+2^e = 2^{j+1}+2^{j+3}+2^{j+6}
\]
after subtracting 1. We now consider all the ways in which $a,c,d,e$ can match to give the correct 2-adic expansion. The number of such combinations is reduced by the assumption on the ordering of $c,d,e$.

If $a = c = j$, then
\[
q_a\cdot h_bh_ch_dh_e = q_j\cdot h_0h_jh_{j+3}h_{j+6}.
\]
However, this class is precisely $d_1(q_{j+1}\cdot h_0h_{j+3}h_{j+6})$ so cannot support any differentials.

If $a = d = j+2$, then
\[
q_a\cdot h_bh_ch_dh_e = q_{j+2}\cdot h_0h_{j+1}h_{j+2}h_{j+6}.
\]
However, the relation $h_{j+1}h_{j+2} = 0$ forces this class to vanish identically, so it cannot support any differentials.

Finally, if $a = e = j+5$, then
\[
q_a\cdot h_bh_ch_dh_e = q_{j+5}\cdot h_0h_{j+1}h_{j+3}h_{j+5}.
\]
However, this class is precisely $d_1(q_{j+6}\cdot h_0h_{j+1}h_{j+3})$ so cannot support any differentials.

This rules out all potential sources of Burklund-Xu differentials hitting $q_0\cdot h_jp_j'$.
\end{proof}

\subsection{Cartan-Eilenberg Differentials}

We have now shown that $q_0\cdot h_jp_j'$ survives to detect a nonzero class on the Cartan-Eilenberg $E_2$ page for $j\geq1$. It remains to show that this class, which we will denote by the same name, survives to detect a nonzero class on the Adams $E_2$ page. We first show that it cannot support any Cartan-Eilenberg differential:
\begin{prop}
The class
\[
q_0\cdot h_jp_j'
\]
does not support any Cartan-Eilenberg differential for $j\geq3$.
\end{prop}
\begin{proof}
The Cartan-Eilenberg filtration of $q_0\cdot h_jp_j'$ is 1. Since Cartan-Eilenberg differentials lower filtration by at least 2, this class cannot possibly support such a differential.
\end{proof}

We now show that $q_0\cdot h_jp_j'$ cannot be the target of any Cartan-Eilenberg differential:
\begin{prop}
\label{CESS p}
The class
\[
q_0\cdot h_jp_j'
\]
is not the target of any Cartan-Eilenberg differential for $j\geq3$.
\end{prop}
\begin{proof}
If a class does support a differential hitting $q_0\cdot h_jp_j'$, then this differential source must have Cartan-Eilenberg filtration $k=3$ or $k=5$. If $k=5$, then the class must be $q_0^5$ since this is the only class surviving to the Cartan-Eilenberg $E_2$ page in filtration 5. This class cannot support any differentials since it must survive to detect $h_0^5$ on the Adams $E_2$ page.

We now rule out any class in filtration $k=3$. Such a class must be of the form
\[
q_aq_bq_c\cdot h_dh_e
\]
since the 2-line $\text{Ext}^2_\mathcal{P}$ is generated by products of $h_i$ classes. If one carries out analysis similar to that used in studying the Burklund-Xu spectral sequence, then one finds the following possibilities for such a class, given that it must have the same total degree as $q_0\cdot h_jp_j'$:
\[
\begin{array}{c c c c c c c c c}
q_1q_{j}^2\cdot h_{j+3}h_{j+6} &&q_1q_{j+5}^2\cdot h_{j+1}h_{j+3} && q_1q_{j+2}^2\cdot h_{j+1}h_{j+6} && q_{j+1}q_{j+3}q_{j+6}\cdot h_0^2 &&\\ \\

q_1q_{j}q_{j+3}\cdot h_{j}h_{j+6} &&q_1q_{j+3}q_{j+5}\cdot h_{j+1}h_{j+5} &&q_1q_{j+1}q_{j+6}\cdot h_{j+2}^2 && q_0q_{j+3}q_{j+6}\cdot h_0h_{j+1} &&\\ \\

q_1q_{j}q_{j+6}\cdot h_{j}h_{j+3} &&q_1q_{j+1}q_{j+5}\cdot h_{j+3}h_{j+5} &&q_{j+1}q_{j+2}^2\cdot h_1h_{j+5} && q_0q_{j+1}q_{j+6}\cdot h_0h_{j+3} &&\\ \\

q_1q_{j+3}q_{j+6}\cdot h_{j}^2 &&q_1q_{j+1}q_{j+3}\cdot h_{j+5}^2 &&q_{j+2}^2q_{j+6}\cdot h_1h_{j+1} && q_0q_{j+1}q_{j+3}\cdot h_0h_{j+6} &&\\ \\

q_{j}^2q_{j+3}\cdot h_1h_{j+6} &&q_{j+1}q_{j+5}^2\cdot h_1h_{j+3} && q_1q_{j+1}q_{j+2}\cdot h_{j+2}h_{j+6} && q_0^2q_{j+6}\cdot h_{j+1}h_{j+3} &&\\ \\

q_{j}^2q_{j+6}\cdot h_1h_{j+3} && q_{j+1}q_{j+5}^2\cdot h_1h_{j+3} && q_1q_{j+2}q_{j+6} \cdot h_{j+1}h_{j+2}=0 && q_0^2q_{j+3}\cdot h_{j+1}h_{j+6} &&\\ \\

q_{j}q_{j+3}q_{j+6}\cdot h_1h_{j} && q_{j+1}q_{j+3}q_{j+5}\cdot h_1h_{j+5} && q_{j+1}q_{j+2}q_{j+6}\cdot h_1h_{j+2} && q_0^2q_{j+1}\cdot h_{j+3}h_{j+6} &&
\end{array}
\]
To be clear, we are using the Burklund-Xu names for these classes since they are not known to survive to the Cartan-Eilenberg $E_2$ page. In fact, we will now show that no linear combination of any of these classes will survive. These classes participate in precisely the following Burklund-Xu $d_1$ and $d_2$ differentials: %where we use $\longrightarrow$ to indicate such a differential:
\begin{enumerate}
    \item $d_2(q_{j+3}q_{j+5}^2\cdot h_1h_{j+1}) = q_{j+1}q_{j+5}^2\cdot h_1h_{j+2}^2 + q_{j+3}q_{j+4}^2\cdot h_1h_{j+1}h_{j+5}$
    \item $d_1(q_{j+1}q_{j+5}^2\cdot h_1h_{j+3}) = q_{j}q_{j+5}^2\cdot h_1h_{j}h_{j+3}$
    \item $d_1(q_{j+2}^2q_{j+6}\cdot h_1h_{j+1}) = q_{j+2}^2q_{j+5}\cdot h_1h_{j+1}h_{j+5}$
    \item $d_1(q_{j+1}q_{j+3}q_{j+6}\cdot h_0^2) = q_{j}q_{j+3}q_{j+6}\cdot h_0^2h_{j} + q_{j+1}q_{j+2}q_{j+6}\cdot h_0^2h_{j+2} + q_{j+1}q_{j+3}q_{j+5}\cdot h_0^2h_{j+5}$
    \item $d_1(q_{j}q_{j+3}q_{j+6}\cdot h_1h_{j}) = q_{j}q_{j+3}q_{j+5}\cdot h_1h_{j}h_{j+5} + q_jq_{j+2}q_{j+6}\cdot h_1h_jh_{j+2}$
    \item $d_1(q_{j+1}q_{j+3}q_{j+5}\cdot h_1h_{j+5}) = q_{j}q_{j+3}q_{j+5}\cdot h_1h_{j}h_{j+5} + q_{j+1}q_{j+2}q_{j+5}\cdot h_1h_{j+2}h_{j+5}$
    \item $d_1(q_{j+1}q_{j+2}q_{j+6}\cdot h_1h_{j+2}) = q_jq_{j+2}q_{j+6}\cdot h_1h_jh_{j+2} + q_{j+1}q_{j+2}q_{j+5}\cdot h_1h_{j+2}h_{j+5}$
    \item $d_1(q_{j}^2q_{j+6}\cdot h_1h_{j+3}) = q_{j}^2q_{j+5}\cdot h_1h_{j+3}h_{j+5}$
    \item $d_1(q_{j+1}q_{j+2}^2\cdot h_1h_{j+6}) = q_{j}q_{j+2}^2\cdot h_1h_{j}h_{j+6}$
    \item $d_1(q_{j}^2q_{j+3}\cdot h_1h_{j+6}) = q_{j}^2q_{j+2}\cdot h_1h_{j+2}h_{j+6}$
    \item $d_1(q_1q_{j+3}q_{j+6}\cdot h_{j}^2) = q_1q_{j+3}q_{j+5}\cdot h_{j}^2h_{j+5} + q_0q_{j+3}q_{j+6}\cdot h_0h_{j}^2 + q_1q_{j+2}q_{j+6}\cdot h_j^2h_{j+2}$
    \item $d_1(q_1q_{j+3}q_{j+5}\cdot h_{j+1}h_{j+5}) = q_0q_{j+3}q_{j+5}\cdot h_0h_{j+1}h_{j+5}$
    \item $d_1(q_0q_{j+3}q_{j+6}\cdot h_0h_{j+1}) = q_0q_{j+3}q_{j+5}\cdot h_0h_{j+1}h_{j+5}$
    \item $d_1(q_1q_{j+1}q_{j+6}\cdot h_{j+2}^2) = q_1q_{j+1}q_{j+5}\cdot h_{j+2}^2h_{j+5}$
    \item $d_1(q_1q_{j}q_{j+6}\cdot h_{j}h_{j+3}) = q_1q_{j}q_{j+5}\cdot h_{j}h_{j+3}h_{j+5} + q_0q_{j}q_{j+6}\cdot h_0h_{j}h_{j+3}$
    \item $d_1(q_1q_{j+1}q_{j+5}\cdot h_{j+3}h_{j+5}) = q_1q_{j}q_{j+5}\cdot h_{j}h_{j+3}h_{j+5} + q_0q_{j+1}q_{j+5}\cdot h_0h_{j+3}h_{j+5}$
    \item $d_1(q_0q_{j+1}q_{j+6}\cdot h_0h_{j+3}) = q_0q_{j}q_{j+6}\cdot h_0h_{j}h_{j+3} + q_0q_{j+1}q_{j+5}\cdot h_0h_{j+3}h_{j+5}$
    \item $d_1(q_1q_{j+1}q_{j+3}\cdot h_{j+5}^2) = q_1q_{j}q_{j+3}\cdot h_{j}h_{j+5}^2 + q_1q_{j+1}q_{j+2}\cdot h_{j+2}h_{j+5}^2 + q_0q_{j+1}q_{j+3}\cdot h_0h_{j+5}^2$
    \item $d_1(q_1q_{j+2}^2\cdot h_{j+1}h_{j+6}) = q_0q_{j+2}^2\cdot h_0h_{j+1}h_{j+6}$
    \item $d_1(q_1q_{j}q_{j+3}\cdot h_{j}h_{j+6}) = q_0q_{j}q_{j+3}\cdot h_0h_{j}h_{j+6} + q_1q_{j}q_{j+2}\cdot h_jh_{j+2}h_{j+6}$
    \item $d_1(q_1q_{j+1}q_{j+2}\cdot h_{j+2}h_{j+6}) = q_0q_{j+1}q_{j+2}\cdot h_0h_{j+2}h_{j+6} + q_1q_jq_{j+2}\cdot h_jh_{j+2}h_{j+6}$
    \item $d_1(q_0q_{j+1}q_{j+3}\cdot h_0h_{j+6}) = q_0q_{j}q_{j+3}\cdot h_0h_{j}h_{j+6} + q_0q_{j+1}q_{j+2}\cdot h_0h_{j+2}h_{j+6}$
    \item $d_1(q_1q_{j}^2\cdot h_{j+3}h_{j+6}) = q_0q_{j}^2\cdot h_0h_{j+3}h_{j+6}$
    \item $d_2(q_0^2q_{j+3}\cdot h_{j+1}h_{j+6}) = q_0^2q_{j+1}\cdot h_{j+2}^2h_{j+6}$
    \item $d_1(q_0^2q_{j+1}\cdot h_{j+3}h_{j+6}) = q_0^2q_{j}\cdot h_{j}h_{j+3}h_{j+6}$
    \item $d_1(q_1q_{j+5}^2\cdot h_{j+1}h_{j+3}) = q_0q_{j+5}^2\cdot h_0h_{j+1}h_{j+3}$
    \item $d_1(q_{j+1}q_{j+2}^2\cdot h_1h_{j+5})=q_jq_{j+2}^2\cdot h_1h_{j}h_{j+5}$
\end{enumerate}

We have arranged these differentials in groups so that those with common terms appear in the same group. Note that all targets of these differentials are nontrivial for $j\geq1$ since we are studying the Burklund-Xu $E_1$ page, which is subject only to relations in $\text{Ext}_\mathcal{P}$. Inspecting the list of differentials, there are four linear combinations which form Burklund-Xu $d_1$ cycles:
\[
\begin{array}{c}
q_0q_{j+3}q_{j+6}\cdot h_0h_{j+1}+q_1q_{j+3}q_{j+5}\cdot h_{j+1}h_{j+5} \\q_0q_{j+1}q_{j+6}\cdot h_0h_{j+3}+q_1q_{j}q_{j+6}\cdot h_{j}h_{j+3}+q_1q_{j+1}q_{j+5}\cdot h_{j+3}h_{j+5} \\
q_0q_{j+1}q_{j+3}\cdot h_{0}h_{j+6}+q_1q_{j+1}q_{j+2}\cdot h_{j+2}h_{j+6}+q_1q_jq_{j+3}\cdot h_jh_{j+6} \\ q_{j}q_{j+3}q_{j+6}\cdot h_1h_{j} + q_{j+1}q_{j+3}q_{j+5}\cdot h_1h_{j+5} +q_{j+1}q_{j+2}q_{j+6}\cdot h_1h_{j+2} 
\end{array}
\]
These linear combinations are, respectively, the targets of Burklund-Xu $d_1$ differentials on the following classes:
\begin{enumerate}
    \item $d_1(q_1q_{j+3}q_{j+6}\cdot h_{j+1}) = q_0q_{j+3}q_{j+6}\cdot h_0h_{j+1} + q_1q_{j+3}q_{j+5}\cdot h_{j+1}h_{j+5}$
    \item $d_1(q_1q_{j+1}q_{j+6}\cdot h_{j+3}) = q_0q_{j+1}q_{j+6}\cdot h_0h_{j+3} + q_1q_{j}q_{j+6}\cdot h_{j}h_{j+3} + q_1q_{j+1}q_{j+5}\cdot h_{j+3}h_{j+5}$
    \item $d_1(q_1q_{j+1}q_{j+3}\cdot h_{j+6}) = q_0q_{j+1}q_{j+3}\cdot h_{0}h_{j+6} + q_1q_{j+1}q_{j+2}\cdot h_{j+2}h_{j+6} + q_1q_jq_{j+3}\cdot h_jh_{j+6}$
    \item $d_1(q_{j+1}q_{j+3}q_{j+6}\cdot h_1) = q_{j}q_{j+3}q_{j+6}\cdot h_1h_{j} + q_{j+1}q_{j+3}q_{j+5}\cdot h_1h_{j+5} + q_{j+1}q_{j+2}q_{j+6}\cdot h_1h_{j+2}$
\end{enumerate}
\begin{comment}
    \[
\begin{array}{c c c}
q_1q_{j+3}q_{j+6}\cdot h_{j+1} && q_1q_{j+1}q_{j+5}\cdot h_{j+3} \\ \\
q_1q_{j+1}q_{j+3}\cdot h_{j+6} &&  q_{j+1}q_{j+3}q_{j+6}\cdot h_1
\end{array}
\]
\end{comment}

For $j\geq 3$, the classes above are distinct and nonzero.
\begin{comment}
for $j=3$, the following classes and differentials changes:%check d2 targets are not killed by d1
\begin{enumerate}
    \item $d_2(q_3q_{7}^2\cdot h_1h_{4})=q_3q_{6}^2\cdot h_1h_{4}h_{7}$
    \item $q_2q_4q_5\cdot h_1h_2=0$
    \item $d_1(q_3q_4q_7\cdot h_1h_7)=q_3^2q_7\cdot h_1h_3h_7$
    \item $d_2(q_3^3\cdot h_1h_{8})=q_1q_3^2\cdot h_2^2h_8+q_1q_2^2\cdot h_1h_3h_8$
\end{enumerate}
\end{comment}
As a result, starting with our original list of possible $q_aq_bq_c\cdot h_dh_e$ classes above, the only Burklund-Xu $d_1$ cycles arising from this list are in fact boundaries in the Burklund-Xu spectral sequence. Therefore, no linear combination of these possible classes survives to detect a nontrivial class in the Cartan-Eilenberg spectral sequence. The claim follows.
\end{proof}

\begin{figure}[H]
\centering
\begingroup
\setlength{\tabcolsep}{6pt}
\renewcommand{\arraystretch}{1.0}
\begin{tabular}{|wc{1.2cm}||wc{2.5cm}|p{5cm}|wc{4cm}|} \hline
      $w$ & $s=1$ & \multicolumn{1}{p{5cm}|}{\centering $s=2$} & $s=3$ \\\hline\hline
% r = 3j+13
    $3j+13$ &  & \tikzmark{part1-figure_p-p1L}$(1)\; q_{j+3}q_{j+5}^2\cdot h_1h_{j+1}$\tikzmark{part1-figure_p-p1R} &  \\\hline
% r = 3j+12
    $3j+12$ &  &  &  \\\hline
% r = 3j+11
    $3j+11$ &  & \tikzmark{part1-figure_p-p4L}$(2)\; q_{j+1}q_{j+5}^2\cdot h_1h_{j+3}$\tikzmark{part1-figure_p-p4R} & \tikzmark{part1-figure_p-p2L}$q_{j+1}q_{j+5}^2\cdot h_1h_{j+2}^2$\tikzmark{part1-figure_p-p2R} \\
     &  &  & \tikzmark{part1-figure_p-p3L}$q_{j+3}q_{j+4}^2\cdot h_1h_{j+1}h_{j+5}$\tikzmark{part1-figure_p-p3R} \\\hline
% r = 3j+10
    $3j+10$ & \tikzmark{part1-figure_p-p60L}$q_{j+1}q_{j+3}q_{j+6}\cdot h_1$\tikzmark{part1-figure_p-p60R} & \tikzmark{part1-figure_p-p6L}$(3)\; q_{j+2}^2q_{j+6}\cdot h_1h_{j+1}$\tikzmark{part1-figure_p-p6R} & \tikzmark{part1-figure_p-p5L}$q_jq_{j+5}^2\cdot h_1h_jh_{j+3}$\tikzmark{part1-figure_p-p5R} \\
     &  & \tikzmark{part1-figure_p-p8L}$(4)\; q_{j+1}q_{j+3}q_{j+6}\cdot h_0^2$\tikzmark{part1-figure_p-p8R} &  \\\hline
% r = 3j+9
    $3j+9$ &  & \tikzmark{part1-figure_p-p12L}$(5)\; q_jq_{j+3}q_{j+6}\cdot h_1h_j$\tikzmark{part1-figure_p-p12R} & \tikzmark{part1-figure_p-p7L}$q_{j+2}^2q_{j+5}\cdot h_1h_{j+1}h_{j+5}$\tikzmark{part1-figure_p-p7R} \\
     &  & \tikzmark{part1-figure_p-p15L}$(6)\; q_{j+1}q_{j+3}q_{j+5}\cdot h_1h_{j+5}$\tikzmark{part1-figure_p-p15R} & \tikzmark{part1-figure_p-p9L}$q_jq_{j+3}q_{j+6}\cdot h_0^2h_j$\tikzmark{part1-figure_p-p9R} \\
     &  & \tikzmark{part1-figure_p-p17L}$(7)\; q_{j+1}q_{j+2}q_{j+6}\cdot h_1h_{j+2}$\tikzmark{part1-figure_p-p17R} & \tikzmark{part1-figure_p-p10L}$q_{j+1}q_{j+2}q_{j+6}\cdot h_0^2h_{j+2}$\tikzmark{part1-figure_p-p10R} \\
     &  &  & \tikzmark{part1-figure_p-p11L}$q_{j+1}q_{j+3}q_{j+5}\cdot h_0^2h_{j+5}$\tikzmark{part1-figure_p-p11R} \\\hline
% r = 3j+8
    $3j+8$ &  &  & \tikzmark{part1-figure_p-p13L}$q_jq_{j+3}q_{j+5}\cdot h_1h_jh_{j+5}$\tikzmark{part1-figure_p-p13R} \\
     &  &  & \tikzmark{part1-figure_p-p14L}$q_jq_{j+2}q_{j+6}\cdot h_1h_jh_{j+2}$\tikzmark{part1-figure_p-p14R} \\
     &  &  & \tikzmark{part1-figure_p-p16L}$q_{j+1}q_{j+2}q_{j+5}\cdot h_1h_{j+2}h_{j+5}$\tikzmark{part1-figure_p-p16R} \\\hline
% r = 3j+7
    $3j+7$ &  &  &  \\\hline
% r = 3j+6
    $3j+6$ &  & \tikzmark{part1-figure_p-p18L}$(8)\; q_j^2q_{j+6}\cdot h_1h_{j+3}$\tikzmark{part1-figure_p-p18R} &  \\\hline
% r = 3j+5
    $3j+5$ &  & \tikzmark{part1-figure_p-p20L}$(9)\; q_{j+1}q_{j+2}^2\cdot h_1h_{j+6}$\tikzmark{part1-figure_p-p20R} & \tikzmark{part1-figure_p-p19L}$q_j^2q_{j+5}\cdot h_1h_{j+3}h_{j+5}$\tikzmark{part1-figure_p-p19R} \\
     &  & \tikzmark{part1-figure_p-p63L}(27) $q_{j+1}q_{j+2}^2\cdot h_1h_{j+5}$\tikzmark{part1-figure_p-p63R} &  \\\hline
% r = 3j+4
    $3j+4$ &  &  & \tikzmark{part1-figure_p-p21L}$q_jq_{j+2}^2\cdot h_1h_jh_{j+6}$\tikzmark{part1-figure_p-p21R} \\
     &  &  & \tikzmark{part1-figure_p-p64L}$q_jq_{j+2}^2\cdot h_1h_jh_{j+5}$\tikzmark{part1-figure_p-p64R} \\\hline
% r = 3j+3
    $3j+3$ &  & \tikzmark{part1-figure_p-p22L}$(10)\; q_j^2q_{j+3}\cdot h_1h_{j+6}$\tikzmark{part1-figure_p-p22R} &  \\\hline
% r = 3j+2
    $3j+2$ &  &  & \tikzmark{part1-figure_p-p23L}$q_j^2q_{j+2}\cdot h_1h_{j+2}h_{j+6}$\tikzmark{part1-figure_p-p23R} \\\hline
    $\vdots$ & & & \\\hline
\end{tabular}
\begin{tikzpicture}[overlay, remember picture, shorten >=.5pt, shorten <=.5pt, transform canvas={yshift=.25\baselineskip}]

    \draw [->] ({pic cs:part1-figure_p-p1R}) [blue,bend left=8] to ({pic cs:part1-figure_p-p2L});
    \draw [->] ({pic cs:part1-figure_p-p1R}) [blue,bend right=8] to ({pic cs:part1-figure_p-p3L});
    \draw [->] ({pic cs:part1-figure_p-p4R}) [red] to ({pic cs:part1-figure_p-p5L});
    \draw [->] ({pic cs:part1-figure_p-p6R}) [red] to ({pic cs:part1-figure_p-p7L});
    \draw [->] ({pic cs:part1-figure_p-p8R}) [red,bend left=16] to ({pic cs:part1-figure_p-p9L});
    \draw [->] ({pic cs:part1-figure_p-p8R}) [red] to ({pic cs:part1-figure_p-p10L});
    \draw [->] ({pic cs:part1-figure_p-p8R}) [red,bend left=0] to ({pic cs:part1-figure_p-p11L});
    \draw [->] ({pic cs:part1-figure_p-p12R}) [red,out=0,in=160] to ({pic cs:part1-figure_p-p13L});
    \draw [->] ({pic cs:part1-figure_p-p12R}) [red,out=0,in=140] to ({pic cs:part1-figure_p-p14L});
    \draw [->] ({pic cs:part1-figure_p-p15R}) [red,bend left=4] to ({pic cs:part1-figure_p-p13L});
    \draw [->] ({pic cs:part1-figure_p-p15R}) [red,bend right=4] to ({pic cs:part1-figure_p-p16L});
    \draw [->] ({pic cs:part1-figure_p-p17R}) [red,bend left=4] to ({pic cs:part1-figure_p-p14L});
    \draw [->] ({pic cs:part1-figure_p-p17R}) [red,bend right=12] to ({pic cs:part1-figure_p-p16L});
    \draw [->] ({pic cs:part1-figure_p-p18R}) [red] to ({pic cs:part1-figure_p-p19L});
    \draw [->] ({pic cs:part1-figure_p-p20R}) [red] to ({pic cs:part1-figure_p-p21L});
    \draw [->] ({pic cs:part1-figure_p-p22R}) [red] to ({pic cs:part1-figure_p-p23L});
    \draw [->] ({pic cs:part1-figure_p-p60R}) [red,bend left=16] to ({pic cs:part1-figure_p-p12L});
    \draw [->] ({pic cs:part1-figure_p-p60R}) [red] to ({pic cs:part1-figure_p-p15L});
    \draw [->] ({pic cs:part1-figure_p-p60R}) [red,bend left=0] to ({pic cs:part1-figure_p-p17L});
    \draw [->] ({pic cs:part1-figure_p-p63R}) [red] to ({pic cs:part1-figure_p-p64L});
\end{tikzpicture}
\par\endgroup
\end{figure}

\begin{figure}[H]
\centering
\begingroup
\setlength{\tabcolsep}{6pt}
\renewcommand{\arraystretch}{1.0}
\begin{tabular}{|wc{1.2cm}||wc{2.5cm}|p{5cm}|wc{4cm}|} \hline
      $w$ & $s=1$ & \multicolumn{1}{p{5cm}|}{\centering $s=2$} & $s=3$ \\\hline\hline
% r = 2j+11
    $2j+11$ &  & \tikzmark{part2-figure_p-p61L}(26) $q_1q_{j+5}^2\cdot h_{j+1}h_{j+3}$\tikzmark{part2-figure_p-p61R} &  \\\hline
% r = 2j+10
    $2j+10$ & \tikzmark{part2-figure_p-p57L}$q_1q_{j+3}q_{j+6}\cdot h_{j+1}$\tikzmark{part2-figure_p-p57R} & \tikzmark{part2-figure_p-p24L}$(11)\; q_1q_{j+3}q_{j+6}\cdot h_j^2$\tikzmark{part2-figure_p-p24R} & \tikzmark{part2-figure_p-p62L}$q_0q_{j+5}^2\cdot h_0h_{j+1}h_{j+3}$\tikzmark{part2-figure_p-p62R} \\\hline
% r = 2j+9
    $2j+9$ &  & \tikzmark{part2-figure_p-p28L}$(12)\; q_1q_{j+3}q_{j+5}\cdot h_{j+1}h_{j+5}$\tikzmark{part2-figure_p-p28R} & \tikzmark{part2-figure_p-p25L}$q_1q_{j+3}q_{j+5}\cdot h_j^2h_{j+5}$\tikzmark{part2-figure_p-p25R} \\
     &  & \tikzmark{part2-figure_p-p30L}$(13)\; q_0q_{j+3}q_{j+6}\cdot h_0h_{j+1}$\tikzmark{part2-figure_p-p30R} & \tikzmark{part2-figure_p-p26L}$q_0q_{j+3}q_{j+6}\cdot h_0h_j^2$\tikzmark{part2-figure_p-p26R} \\
     &  &  & \tikzmark{part2-figure_p-p27L}$q_1q_{j+2}q_{j+6}\cdot h_j^2h_{j+2}$\tikzmark{part2-figure_p-p27R} \\\hline
% r = 2j+8
    $2j+8$ & \tikzmark{part2-figure_p-p58L}$q_1q_{j+1}q_{j+6}\cdot h_{j+3}$\tikzmark{part2-figure_p-p58R} & \tikzmark{part2-figure_p-p31L}$(14)\; q_1q_{j+1}q_{j+6}\cdot h_{j+2}^2$\tikzmark{part2-figure_p-p31R} & \tikzmark{part2-figure_p-p29L}$q_0q_{j+3}q_{j+5}\cdot h_0h_{j+1}h_{j+5}$\tikzmark{part2-figure_p-p29R} \\\hline
% r = 2j+7
    $2j+7$ &  & \tikzmark{part2-figure_p-p33L}$(15)\; q_1q_jq_{j+6}\cdot h_jh_{j+3}$\tikzmark{part2-figure_p-p33R} & \tikzmark{part2-figure_p-p32L}$q_1q_{j+1}q_{j+5}\cdot h_{j+2}^2h_{j+5}$\tikzmark{part2-figure_p-p32R} \\
     &  & \tikzmark{part2-figure_p-p36L}$(16)\; q_1q_{j+1}q_{j+5}\cdot h_{j+3}h_{j+5}$\tikzmark{part2-figure_p-p36R} &  \\
     &  & \tikzmark{part2-figure_p-p38L}$(17)\; q_0q_{j+1}q_{j+6}\cdot h_0h_{j+3}$\tikzmark{part2-figure_p-p38R} &  \\\hline
% r = 2j+6
    $2j+6$ &  &  & \tikzmark{part2-figure_p-p34L}$q_1q_jq_{j+5}\cdot h_jh_{j+3}h_{j+5}$\tikzmark{part2-figure_p-p34R} \\
     &  &  & \tikzmark{part2-figure_p-p35L}$q_0q_jq_{j+6}\cdot h_0h_jh_{j+3}$\tikzmark{part2-figure_p-p35R} \\
     &  &  & \tikzmark{part2-figure_p-p37L}$q_0q_{j+1}q_{j+5}\cdot h_0h_{j+3}h_{j+5}$\tikzmark{part2-figure_p-p37R} \\\hline
% r = 2j+5
    $2j+5$ & \tikzmark{part2-figure_p-p59L}$q_1q_{j+1}q_{j+3}\cdot h_{j+6}$\tikzmark{part2-figure_p-p59R} & \tikzmark{part2-figure_p-p39L}$(18)\; q_1q_{j+1}q_{j+3}\cdot h_{j+5}^2$\tikzmark{part2-figure_p-p39R} &  \\
     &  & \tikzmark{part2-figure_p-p43L}$(19)\; q_1q_{j+2}^2\cdot h_{j+1}h_{j+6}$\tikzmark{part2-figure_p-p43R} &  \\\hline
% r = 2j+4
    $2j+4$ &  & \tikzmark{part2-figure_p-p45L}$(20)\; q_1q_jq_{j+3}\cdot h_jh_{j+6}$\tikzmark{part2-figure_p-p45R} & \tikzmark{part2-figure_p-p40L}$q_1q_jq_{j+3}\cdot h_jh_{j+5}^2$\tikzmark{part2-figure_p-p40R} \\
     &  & \tikzmark{part2-figure_p-p48L}$(21)\; q_1q_{j+1}q_{j+2}\cdot h_{j+2}h_{j+6}$\tikzmark{part2-figure_p-p48R} & \tikzmark{part2-figure_p-p41L}$q_1q_{j+1}q_{j+2}\cdot h_{j+2}h_{j+5}^2$\tikzmark{part2-figure_p-p41R} \\
     &  & \tikzmark{part2-figure_p-p50L}$(22)\; q_0q_{j+1}q_{j+3}\cdot h_0h_{j+6}$\tikzmark{part2-figure_p-p50R} & \tikzmark{part2-figure_p-p42L}$q_0q_{j+1}q_{j+3}\cdot h_0h_{j+5}^2$\tikzmark{part2-figure_p-p42R} \\
     &  &  & \tikzmark{part2-figure_p-p44L}$q_0q_{j+2}^2\cdot h_0h_{j+1}h_{j+6}$\tikzmark{part2-figure_p-p44R} \\\hline
% r = 2j+3
    $2j+3$ &  &  & \tikzmark{part2-figure_p-p46L}$q_0q_jq_{j+3}\cdot h_0h_jh_{j+6}$\tikzmark{part2-figure_p-p46R} \\
     &  &  & \tikzmark{part2-figure_p-p47L}$q_1q_jq_{j+2}\cdot h_jh_{j+2}h_{j+6}$\tikzmark{part2-figure_p-p47R} \\
     &  &  & \tikzmark{part2-figure_p-p49L}$q_0q_{j+1}q_{j+2}\cdot h_0h_{j+2}h_{j+6}$\tikzmark{part2-figure_p-p49R} \\\hline
% r = 2j+2
    $2j+2$ &  &  &  \\\hline
% r = 2j+1
    $2j+1$ &  & \tikzmark{part2-figure_p-p51L}$(23)\; q_1q_j^2\cdot h_{j+3}h_{j+6}$\tikzmark{part2-figure_p-p51R} &  \\\hline
% r = 2j
    $2j$ &  &  & \tikzmark{part2-figure_p-p52L}$q_0q_j^2\cdot h_0h_{j+3}h_{j+6}$\tikzmark{part2-figure_p-p52R} \\\hline
    $\vdots$ & & & \\\hline
% r = j+3
    $j+3$ &  & \tikzmark{part2-figure_p-p53L}$(24)\; q_0^2q_{j+3}\cdot h_{j+1}h_{j+6}$\tikzmark{part2-figure_p-p53R} &  \\\hline
% r = j+2
    $j+2$ &  &  &  \\\hline
% r = j+1
    $j+1$ &  & \tikzmark{part2-figure_p-p55L}$(25)\; q_0^2q_{j+1}\cdot h_{j+3}h_{j+6}$\tikzmark{part2-figure_p-p55R} & \tikzmark{part2-figure_p-p54L}$q_0^2q_{j+1}\cdot h_{j+2}^2h_{j+6}$\tikzmark{part2-figure_p-p54R} \\\hline
% r = j
    $j$ &  &  & \tikzmark{part2-figure_p-p56L}$q_0^2q_j\cdot h_jh_{j+3}h_{j+6}$\tikzmark{part2-figure_p-p56R} \\\hline
    \end{tabular}
\begin{tikzpicture}[overlay, remember picture, shorten >=.5pt, shorten <=.5pt, transform canvas={yshift=.25\baselineskip}]
    \draw [->] ({pic cs:part2-figure_p-p24R}) [red,bend left=16] to ({pic cs:part2-figure_p-p25L});
    \draw [->] ({pic cs:part2-figure_p-p24R}) [red,bend left=8] to ({pic cs:part2-figure_p-p26L});
    \draw [->] ({pic cs:part2-figure_p-p24R}) [red,bend left=15] to ({pic cs:part2-figure_p-p27L});
    \draw [->] ({pic cs:part2-figure_p-p28R}) [red,bend left=4] to ({pic cs:part2-figure_p-p29L});
    \draw [->] ({pic cs:part2-figure_p-p30R}) [red,bend right=4] to ({pic cs:part2-figure_p-p29L});
    \draw [->] ({pic cs:part2-figure_p-p31R}) [red] to ({pic cs:part2-figure_p-p32L});
    \draw [->] ({pic cs:part2-figure_p-p33R}) [red,bend left=15] to ({pic cs:part2-figure_p-p34L});
    \draw [->] ({pic cs:part2-figure_p-p33R}) [red,out=0,in=160] to ({pic cs:part2-figure_p-p35L});
    \draw [->] ({pic cs:part2-figure_p-p36R}) [red,bend left=4] to ({pic cs:part2-figure_p-p34L});
    \draw [->] ({pic cs:part2-figure_p-p36R}) [red,bend right=4] to ({pic cs:part2-figure_p-p37L});
    \draw [->] ({pic cs:part2-figure_p-p38R}) [red,bend left=4] to ({pic cs:part2-figure_p-p35L});
    \draw [->] ({pic cs:part2-figure_p-p38R}) [red,bend right=12] to ({pic cs:part2-figure_p-p37L});
    \draw [->] ({pic cs:part2-figure_p-p39R}) [red,bend left=16] to ({pic cs:part2-figure_p-p40L});
    \draw [->] ({pic cs:part2-figure_p-p39R}) [red] to ({pic cs:part2-figure_p-p41L});
    \draw [->] ({pic cs:part2-figure_p-p39R}) [red,bend right=16] to ({pic cs:part2-figure_p-p42L});
    \draw [->] ({pic cs:part2-figure_p-p43R}) [red] to ({pic cs:part2-figure_p-p44L});
    \draw [->] ({pic cs:part2-figure_p-p45R}) [red,out=0,in=160] to ({pic cs:part2-figure_p-p46L});
    \draw [->] ({pic cs:part2-figure_p-p45R}) [red,out=0,in=140] to ({pic cs:part2-figure_p-p47L});
    \draw [->] ({pic cs:part2-figure_p-p48R}) [red,bend left=4] to ({pic cs:part2-figure_p-p47L});
    \draw [->] ({pic cs:part2-figure_p-p48R}) [red,bend right=4] to ({pic cs:part2-figure_p-p49L});
    \draw [->] ({pic cs:part2-figure_p-p50R}) [red,bend left=4] to ({pic cs:part2-figure_p-p46L});
    \draw [->] ({pic cs:part2-figure_p-p50R}) [red,bend right=12] to ({pic cs:part2-figure_p-p49L});
    \draw [->] ({pic cs:part2-figure_p-p51R}) [red] to ({pic cs:part2-figure_p-p52L});
    \draw [->] ({pic cs:part2-figure_p-p53R}) [blue] to ({pic cs:part2-figure_p-p54L});
    \draw [->] ({pic cs:part2-figure_p-p55R}) [red] to ({pic cs:part2-figure_p-p56L});
    \draw [->] ({pic cs:part2-figure_p-p57R}) [red,bend left=8] to ({pic cs:part2-figure_p-p28L});
    \draw [->] ({pic cs:part2-figure_p-p57R}) [red,bend right=8] to ({pic cs:part2-figure_p-p30L});
    \draw [->] ({pic cs:part2-figure_p-p58R}) [red,bend left=16] to ({pic cs:part2-figure_p-p33L});
    \draw [->] ({pic cs:part2-figure_p-p58R}) [red] to ({pic cs:part2-figure_p-p36L});
    \draw [->] ({pic cs:part2-figure_p-p58R}) [red,bend left=0] to ({pic cs:part2-figure_p-p38L});
    \draw [->] ({pic cs:part2-figure_p-p59R}) [red,bend left=16] to ({pic cs:part2-figure_p-p45L});
    \draw [->] ({pic cs:part2-figure_p-p59R}) [red] to ({pic cs:part2-figure_p-p48L});
    \draw [->] ({pic cs:part2-figure_p-p59R}) [red,bend left=0] to ({pic cs:part2-figure_p-p50L});
    \draw [->] ({pic cs:part2-figure_p-p61R}) [red] to ({pic cs:part2-figure_p-p62L});
\end{tikzpicture}
\par\endgroup
\caption{The BXSS in degree $(s,k,t) = (2,3,2^{j+2}\cdot 37+1)$ for \cref{CESS p}, where $j\geq 3$.
The vertical axis is the BX filtration and the horizontal axis is the $s$-degree. Red arrows indicate $d_1$-differentials, and blue arrows indicate $d_2$-differentials.}
\label{fig:figure_p}
\label{fig:figure_p-part2}
\end{figure}

\section{The $D_3(j)$ Family}
In this section we prove the following theorem:
\begin{thm}
There are Adams differentials
\[
d_2(D_3(j+1)) = h_0K_j \neq 0
\]
for all $j\geq1$.
\end{thm}

The proof proceeds in three steps, following the program of Burklund-Xu. First we show, for $j\geq1$, that the element $q_0\cdot K_j$ is a permanent cycle in the Burklund-Xu spectral sequence and detects a nonzero class in the abutment, an associated graded of the Cartan-Eilenberg $E_2$ page. We then show that $q_0\cdot K_j$ is a permanent cycle in that spectral sequence and detects a nonzero class in the abutment, an associated graded of the Adams $E_2$ page. Since $K_j$ detects the classical Adams element $K_{j+1}$ under the doubling-up isomorphism
\[
\text{Ext}_\mathcal{A} \xrightarrow{\cong} \text{Ext}_\mathcal{P},
\]
\cite{Benson} shows that $q_0\cdot K_j$ must detect $h_0K_{j+1}$. Since we do this for all $j\geq1$, this proves the theorem for $j\geq2$. The remaining case of $j = 1$, namely
\[
d_2(D_3(2)) = h_0K_1,
\]
is handled directly at the end of the section.

\subsection{Burklund-Xu Differentials}
We first show the following:
\begin{prop}
The Burklund-Xu element $q_0\cdot K_j$ is a permanent cycle for $j\geq0$.
\end{prop}
\begin{proof}
The Burklund-Xu filtration of $q_0\cdot K_j$ is 0. Since Burklund-Xu differentials strictly lower filtration, this class cannot support any differential in that spectral sequence.
\end{proof}

Then the class $q_0\cdot K_j$ is defined in the abutment of the Burklund-Xu spectral sequence for all $j\geq0$, so it remains to show that it cannot receive any differential. We use the following simple proposition to simplify this argument:
\begin{prop}\label{prop:deg-Dj}
If $y\in\text{Ext}^4_\mathcal{P}$ has internal degree $t(y)$ divisible by $2^2$, then any Burklund-Xu element $q_a\cdot y$ in the same internal degree as $q_0\cdot K_j$ must not support a differential.
\end{prop}
\begin{proof}
If $q_a\cdot y$ lies in the same internal degree as $q_0\cdot K_j$, then the equality
\[
t(q_a\cdot y) = 2^{a+1}-1+t(y) = 1+2^{j+2}+2^{j+8} = t(q_0\cdot K_j)
\]
must hold. Adding 1 to both sides and dividing by 2 gives
\[
2^a + 2^{-1}t(y) = 1+2^{j+2}+2^{j+8},
\]
and, by assumption, the term $2^{-1}t(y)$ is still divisible by 2. Then the only possibility is that $a=0$, meaning that
\[
q_a\cdot y = q_0\cdot y
\]
is the element in question.
\end{proof}

We now show that $q_0\cdot K_j$ cannot be the target of any Burklund-Xu differential:

\begin{prop}
There is no Burklund-Xu differential of the form
\[
d_r(x) = q_0\cdot K_j
\]
for any $j\geq0$ and Burklund-Xu element $x$.
\end{prop}
\begin{proof}
Recall that we can write $x$ as a linear combination of elements of the form $q_a\cdot y$ for $y\in\text{Ext}^4_\mathcal{P}$. As a result, we first collect all classes $q_a\cdot y$ having the same internal degree as $q_0\cdot K_j$ and then show that no linear combination of these classes can support a differential as in the proposition statement. We begin with the case where $y$ is indecomposable (see \cite{Chen}, \cite{Lin} for a full presentation of the classical Adams 4-line). Of the possible cases, we find that if $y\in\{d_b,f_b,g_b\}$, then \autoref{prop:deg-Dj} implies that $q_a\cdot y$ = $q_0\cdot y$. A class of this form cannot support any Burklund-Xu differential since it already lies in the lowest possible filtration, so we discard these classes.

We now consider $y = e_b$ so that the equality
\[
t(q_a\cdot e_b) = 2^{a+1}-1+2^{b+1}+2^{b+3}+2^{b+5} = 1+2^{j+2}+2^{j+8} = t(q_0\cdot K_j)
\]
holds. Adding 1 to both sides and dividing by 2 gives
\[
2^a+2^b+2^{b+2}+2^{b+4} = 1+2^{j+1}+2^{j+7},
\]
and since we must avoid $a=0$ by the reasoning above, we conclude that $b=0$. After subtracting 1 from both sides, this gives
\[
2^a+2^2+2^4 = 2^{j+1}+2^{j+7},
\]
but there is no choice of $a$ and $j$ compatible with this relation. Then $y\neq e_b$.

We now consider $y = p_b$ so that the equality
\[
t(q_a\cdot p_b) = 2^{a+1}-1+2^{b+1}+2^{b+3}+2^{b+6} = 1+2^{j+2}+2^{j+8} = t(q_0\cdot K_j)
\]
holds. Adding 1 to both sides and dividing by 2 gives
\[
2^a+2^b+2^{b+2}+2^{b+5} = 1+2^{j+1}+2^{j+7},
\]
and since we must avoid $a=0$ by the reasoning above, we conclude that $b=0$. After subtracting 1 from both sides, this gives
\[
2^a+2^2+2^5 = 2^{j+1}+2^{j+7},
\]
but there is no choice of $a$ and $j$ compatible with this relation. Then $y\neq p_b$.

We now consider $y = D_3(b)$ so that the equality
\[
t(q_a\cdot D_3(b)) = 2^{a+1}-1+2^{b+1}+2^{b+7} = 1+2^{j+2}+2^{j+8} = t(q_0\cdot K_j)
\]
holds. Adding 1 to both sides and dividing by 2 gives
\[
2^a+2^b+2^{b+6} = 1+2^{j+1}+2^{j+7},
\]
and since we must avoid $a=0$ by the reasoning above, we conclude that $b=0$. After subtracting 1 from both sides, this gives
\[
2^a+2^6 = 2^{j+1}+2^{j+7}.
\]
In order to avoid $j=-1$, which would make $q_0\cdot K_j$ undefined, we must have $j=5$ and $a=12$ so that
\[
q_{12}\cdot D_3(0)
\]
lies in the same internal degree as $q_0\cdot K_5$. We need to consider linear combinations of this class with other classes which we may discover below.

We now consider $y = p_b'$ so that the equality
\[
t(q_a\cdot p_b') = 2^{a+1}-1+2^{b+1}+2^{b+4}+2^{b+7} = 1+2^{j+2}+2^{j+8} = t(q_0\cdot K_j)
\]
holds. Adding 1 to both sides and dividing by 2 gives
\[
2^a+2^b+2^{b+3}+2^{b+6} = 1+2^{j+1}+2^{j+7},
\]
and since we must avoid $a=0$ by the reasoning above, we conclude that $b=0$. After subtracting 1 from both sides, this gives
\[
2^a+2^3+2^6 = 2^{j+1}+2^{j+7},
\]
but there is no choice of $a$ and $j$ compatible with this relation. Then $y\neq p_b'$, and we have exhausted all indecomposable $y$. So far we have only $q_{12}\cdot D_3(0)$ to consider.

We now consider decomposable $y$, which are of the form $h_\ell c_m$ and $h_bh_ch_dh_e$. We first consider $y = h_\ell c_m$ so that the equality
\[
t(q_a\cdot h_\ell c_m) = 2^{a+1}-1+2^{\ell+1}+2^{m+1}+2^{m+2}+2^{m+4} = 1+2^{j+2}+2^{j+8} = t(q_0\cdot K_j)
\]
holds. Adding 1 to both sides and dividing by 2 gives
\[
2^a+2^\ell+2^m+2^{m+1}+2^{m+3} = 1+2^{j+1}+2^{j+7},
\]
and since we must avoid $a=0$ by the reasoning above, we conclude that $\ell = 0$ or $m = 0$. We cannot have $\ell = 0$, as there would be no choice of $a$, $m$, and $j$ compatible with this relation, so it must be that $m = 0$. Then we obtain
\[
2^a+2^\ell+2^1+2^3 = 2^{j+1}+2^{j+7}
\]
after setting $m=0$ and subtracting 1 from the resulting equation. There is no choice of $a$, $\ell$, and $j$ compatible with this relation, so $y \neq h_\ell c_m$.

The only remaining classes $y\in\text{Ext}^4_\mathcal{P}$ are of the form $y = h_bh_ch_dh_e$ so that the equality
\[
t(q_a\cdot h_bh_ch_dh_e) = 2^{a+1}-1+2^{b+1}+2^{c+1}+2^{d+1}+2^{e+1} = 1+2^{j+2}+2^{j+8} = t(q_0\cdot K_j)
\]
holds. Adding 1 to both sides and dividing by 2 gives
\[
2^a+2^b+2^c+2^d+2^e = 1+2^{j+1}+2^{j+7},
\]
and since we must avoid $a=0$ by the reasoning above, we conclude that $b=0$ since we can take $b\leq c\leq d\leq e$ without loss of generality. After subtracting 1 from both sides, this gives
\[
2^a+2^c+2^d+2^e = 2^{j+1}+2^{j+7}.
\]
Again using that $c\leq d\leq e$, we have the following possibilities:
\[
\begin{array}{c c c}
a = c = j,\text{ }d = e = j+6 &&a = j+1,\text{ }c = d = j+5,\text{ }e = j+6 \\ \\
c = j+1,\text{ }a = d = j+5,\text{ }e = j+6 && c = j+1,\text{ }d = e = j+5,\text{ }a = j+6
\end{array}
\]
These give, respectively, the following elements:
\[
\begin{array}{c c c}
q_j\cdot h_0h_jh_{j+6}^2 &&q_{j+1}\cdot h_0h_{j+5}^2h_{j+6} = 0 \\ \\
q_{j+5}\cdot h_0h_{j+1}h_{j+5}h_{j+6} = 0 &&q_{j+6}\cdot h_0h_{j+1}h_{j+5}^2
\end{array}
\]
This concludes the search for elements in the same internal degree as $q_0\cdot K_j$, so we collect all the nonzero classes we have found:
\[
\begin{array}{c c c c c}
q_{12}\cdot D_3(0) &&q_j\cdot h_0h_jh_{j+6}^2 &&q_{j+6}\cdot h_0h_{j+1}h_{j+5}^2
\end{array}
\]
Note that the middle class is trivial since there is a Burklund-Xu $d_1$ differential
\[
d_1(q_{j+1}\cdot h_0h_{j+6}^2) = q_j\cdot h_0h_jh_{j+6}^2
\]
for all $j\neq 1$ and a $d_2$ differential
\[
d_2(q_2\cdot h_0h_7^2)=q_0\cdot h_1^2h_7^2
\]
for $j=1$. Of the remaining two classes, we note that
\[
d_1(q_{12}\cdot D_3(0)) = q_{11}\cdot h_{11}D_3(0) \neq 0
\]
and
\[
d_1(q_{j+6}\cdot h_0h_{j+1}h_{j+5}^2) = q_{j+5}\cdot h_0h_{j+1}h_{j+5}^3 = q_{j+5}\cdot h_0h_{j+1}h_{j+4}^2h_{j+6} \neq 0
\]
cannot be equal for any choice of $j$. As a result, neither class survives long enough to participate in a Burklund-Xu differential with target $q_0\cdot K_j$.

This rules out all potential sources of Burklund-Xu differentials hitting $q_0\cdot K_j$.
\end{proof}

\subsection{Cartan-Eilenberg Differentials}
\begin{prop}
\label{CESS D_3}
    For $j\geq 1$ the class
    \[ q_0 \cdot K_j\]
    survives in the Cartan-Eilenberg spectral sequence.
\end{prop}
\begin{proof}
    $q_0 \cdot K_j$ has Cartan-Eilenberg filtration 1, so it cannot support a Cartan-Eilenberg differential, which lower the Cartan-Eilenberg filtration by 2.
    
    The source of differentials that could hit $q_0\cdot K_j$ must have CE filtration $k=3$ or $k=5$.\\
    If $k=5$, then the possible classes are 
    \begin{enumerate}
        \item $q_0^3q_{j+1}q_{j+7}$
        \item $q_0q_1q_j^2q_{j+7}$
        \item $q_0q_1q_{j+1}q_{j+6}^2$
    \end{enumerate} in degree $(s,k,t) = (0,5,2^{j+2}+2^{j+8}+1)$. However, they all support non-trivial Burklund-Xu differentials:
    \begin{enumerate}
        \item $d_1(q_0^3q_{j+1}q_{j+7}) = q_0^3q_jq_{j+7}\cdot h_{j}+q_0^3q_{j+1}q_{j+6}\cdot h_{j+6}\not=0$
        \item $d_1(q_0q_1q_j^2q_{j+7})=q_0^2q_j^2q_{j+7}h_0\not=0$
        \item $d_1(q_0q_1q_{j+1}q_{j+6}^2) = q_0^2q_{j+1}q_{j+6}^2\cdot h_0\not=0$
    \end{enumerate}
    and does not survive to the CESS.\\
    If $k=3$, then the possible classes in degree $(s,k,t) = (2,3,2^{j+2}+2^{j+8}+1)$ are the following (all distinct and nonzero for $j\ge 2$):
    \begin{enumerate}
        \item $q_{j+1}q_{j+6}^2\cdot h_0^2$
        \item $q_{j}q_{j+6}^2\cdot h_1h_{j}$
        \item $q_{j}^2q_{j+6}\cdot h_1h_{j+6}$
        \item $q_j^2q_{j+7}\cdot h_0^2$
        \item $q_1q_{j+6}^2\cdot h_{j}^2$
        \item $q_0q_{j+6}^2\cdot h_0h_{j+1}$
        \item $q_1q_{j}q_{j+6} \cdot h_{j}h_{j+6}$
        \item $q_0q_{j+1}q_{j+6}\cdot h_0h_{j+6}$
        \item $q_0q_jq_{j+7}\cdot h_0h_j$
        \item $q_1q_{j}^2\cdot h_{j+6}^2$
        \item $q_0q_j^2\cdot h_0h_{j+7}$
        \item $q_0^2q_{j+7}\cdot h_j^2$
        \item $q_0^2q_{j+1}\cdot h_{j+6}^2$
        \item $q_0^2q_{j+6}\cdot h_{j+1}h_{j+6}$ 
        \item $q_0^2q_{j}\cdot h_{j}h_{j+7}$
    \end{enumerate}
    The last two class $q_0q_{j+1}q_{j+7} \cdot h_{j+1}h_{j+7}$ are killed by Burklund-Xu differentials:
    \[d_1(q_0^2q_{j+7}\cdot h_{j+1}) = q_0^2q_{j+6}\cdot h_{j+1}h_{j+6}.\]
    \[d_1(q_0^2q_{j+1}\cdot h_{j+7})=q_0^2q_{j}\cdot h_{j}h_{j+7}\]
    Other classes all support non-trivial Burklund-Xu differentials: %\textbf{(address (1), (2), (3), (7), (8), (9), (11))}
    \begin{enumerate}
        \item $d_1(q_{j+1}q_{j+6}^2\cdot h_0^2)=q_{j}q_{j+6}^2\cdot h_0^2h_j$ %d_2(q_2q_7^2\cdot h_0) = q_2q_6^2\cdot h_0^2h_7 \neq 0 \in E_2
        \item $d_2(q_jq_{j+6}^2\cdot h_1h_j) = q_jq_{j+5}^2\cdot h_1h_jh_{j+6}$ %source is d_1(q_{j+1}q_{j+6}^2\cdot h_1)
        \item $d_2(q_j^2q_{j+6}\cdot h_1h_{j+6}) = q_{j-1}^2q_{j+6}\cdot h_1h_jh_{j+6}$ %source is d_1(q_j^2 q_{j+7}\cdot h_1)
        \item $d_1(q_j^2q_{j+7}\cdot h_0^2)=q_j^2q_{j+6}\cdot h_0^2h_{j+6}$ %fine
        \item $d_1(q_1q_{j+6}^2\cdot h_j^2) = q_0q_{j+6}^2\cdot h_j^3$
        \item $d_2(q_0q_{j+6}^2\cdot h_0h_{j+1}) = q_0q_{j+5}^2\cdot h_0h_{j+1}h_{j+6}$ %d_1(q_1q_{j+6}^2\cdot h_{j+1})=q_0q_{j+6}^2\cdot h_0h_{j+1}
        \item $d_1(q_1q_{j}q_{j+6} \cdot h_{j}h_{j+6}) = q_0q_{j}q_{j+6} \cdot h_0h_{j}h_{j+6}$
        \item $d_1(q_0q_{j+1}q_{j+6}\cdot h_0h_{j+6})=q_0q_{j}q_{j+6} \cdot h_0h_{j}h_{j+6}$
        \item $d_1(q_0q_jq_{j+7}\cdot h_0h_j)=q_0q_{j}q_{j+6} \cdot h_0h_{j}h_{j+6}$
        \item $d_1(q_1q_j^2\cdot h_{j+6}^2) = q_0q_j^2\cdot h_0h_{j+6}^2$
        \item $d_2(q_0q_j^2\cdot h_0h_{j+7})=q_0q_{j-1}^2\cdot h_0h_jh_{j+7}$ %source is d_1 boundary
        \item $d_1(q_0^2q_{j+7}\cdot h_j^2) = q_0^2q_{j+6}\cdot h_j^2h_{j+6}$
        \item $d_1(q_0^2q_{j+1}\cdot h_{j+6}^2) = q_0^2q_{j}\cdot h_jh_{j+6}^2$
    \end{enumerate}
    Note there are linear relations among targets of these differentials. In particular, we should account for
    \begin{enumerate}
        \item $d_1(q_1q_{j}q_{j+6} \cdot h_{j}h_{j+6}+q_0q_{j+1}q_{j+6}\cdot h_0h_{j+6}) = 0$
        \item $d_1(q_0q_{j+1}q_{j+6}\cdot h_0h_{j+6}+q_0q_jq_{j+7}\cdot h_0h_j)=0$
        \item $d_1(q_1q_{j}q_{j+6} \cdot h_{j}h_{j+6}+q_0q_jq_{j+7}\cdot h_0h_j)=0$
    \end{enumerate}
    These classes are indeed killed by Burklund-Xu differentials:
    \begin{enumerate}
        \item $d_1(q_1q_{j+1}q_{j+6}\cdot h_{j+6}) = q_1q_{j}q_{j+6} \cdot h_{j}h_{j+6}+q_0q_{j+1}q_{j+6}\cdot h_0h_{j+6}$
        \item $d_1(q_0q_{j+1}q_{j+7}\cdot h_0)=q_0q_{j+1}q_{j+6}\cdot h_0h_{j+6}+q_0q_jq_{j+7}\cdot h_0h_j$
        \item $d_1(q_1q_jq_{j+7}\cdot h_j)=q_1q_{j}q_{j+6} \cdot h_{j}h_{j+6}+q_0q_jq_{j+7}\cdot h_0h_j$
    \end{enumerate}

\begin{figure}[H]
\centering
\begingroup
\setlength{\tabcolsep}{6pt}
\renewcommand{\arraystretch}{1.1}
\begin{tabular}{|wc{1.2cm}||wc{2.5cm}|p{5cm}|wc{4cm}|} \hline
      $w$ & $s=1$ & \multicolumn{1}{p{5cm}|}{\centering $s=2$} & $s=3$ \\\hline\hline
% r = 3j+13
    $3j+13$ &  & $(1)\; q_{j+1}q_{j+6}^2\cdot h_0^2$\tikzmark{part1-figure_D_3-n1R} &  \\\hline
% r = 3j+12
    $3j+12$ &  & $(2)\; q_jq_{j+6}^2\cdot h_1h_j$\tikzmark{part1-figure_D_3-n2R} & \tikzmark{part1-figure_D_3-t1L}$q_jq_{j+6}^2\cdot h_0^2h_j$ \\\hline
% r = 3j+11
    $3j+11$ &  &  &  \\\hline
% r = 3j+10
    $3j+10$ &  &  & \tikzmark{part1-figure_D_3-t2L}$q_jq_{j+5}^2\cdot h_1h_jh_{j+6}$ \\\hline
% r = 3j+9
    $3j+9$ &  &  &  \\\hline
% r = 3j+8
    $3j+8$ &  &  &  \\\hline
% r = 3j+7
    $3j+7$ &  & $(4)\; q_j^2q_{j+7}\cdot h_0^2$\tikzmark{part1-figure_D_3-n4R} &  \\\hline
% r = 3j+6
    $3j+6$ &  & $(3)\; q_j^2q_{j+6}\cdot h_1h_{j+6}$\tikzmark{part1-figure_D_3-n3R} & \tikzmark{part1-figure_D_3-t4L}$q_j^2q_{j+6}\cdot h_0^2h_{j+6}$ \\\hline
% r = 3j+5
    $3j+5$ &  &  &  \\\hline
% r = 3j+4
    $3j+4$ &  &  & \tikzmark{part1-figure_D_3-t3L}$q_{j-1}^2q_{j+6}\cdot h_1h_jh_{j+6}$ \\\hline
% r = \vdots
    $\vdots$ &  &  &  \\\hline
% r = 2j+13
    $2j+13$ &  & $(5)\; q_1q_{j+6}^2\cdot h_j^2$\tikzmark{part1-figure_D_3-n5R} &  \\\hline
% r = 2j+12
    $2j+12$ &  & $(6)\; q_0q_{j+6}^2\cdot h_0h_{j+1}$\tikzmark{part1-figure_D_3-n6R} & \tikzmark{part1-figure_D_3-t5L}$q_0q_{j+6}^2\cdot h_j^3$ \\\hline
% r = 2j+11
    $2j+11$ &  &  &  \\\hline
% r = 2j+10
    $2j+10$ &  &  & \tikzmark{part1-figure_D_3-t6L}$q_0q_{j+5}^2\cdot h_0h_{j+1}h_{j+6}$ \\\hline
% r = 2j+9
    $2j+9$ &  &  &  \\\hline
\end{tabular}
\begin{tikzpicture}[overlay, remember picture, shorten >=.5pt, shorten <=.5pt, transform canvas={yshift=.25\baselineskip}]

    \draw [->] ({pic cs:part1-figure_D_3-n1R}) [red] to ({pic cs:part1-figure_D_3-t1L});
    \draw [->] ({pic cs:part1-figure_D_3-n2R}) [blue] to ({pic cs:part1-figure_D_3-t2L});
    \draw [->] ({pic cs:part1-figure_D_3-n3R}) [blue] to ({pic cs:part1-figure_D_3-t3L});
    \draw [->] ({pic cs:part1-figure_D_3-n4R}) [red] to ({pic cs:part1-figure_D_3-t4L});
    \draw [->] ({pic cs:part1-figure_D_3-n5R}) [red] to ({pic cs:part1-figure_D_3-t5L});
    \draw [->] ({pic cs:part1-figure_D_3-n6R}) [blue] to ({pic cs:part1-figure_D_3-t6L});
\end{tikzpicture}
\par\endgroup
\end{figure}

\begin{figure}[H]
\centering
\begingroup
\setlength{\tabcolsep}{6pt}
\renewcommand{\arraystretch}{1.1}
\begin{tabular}{|wc{1.2cm}||wc{2.5cm}|p{5cm}|wc{4cm}|} \hline
      $w$ & $s=1$ & \multicolumn{1}{p{5cm}|}{\centering $s=2$} & $s=3$ \\\hline\hline
% r = 2j+8
    $2j+8$ & $q_1q_{j+1}q_{j+6}\cdot h_{j+6}$\tikzmark{part2-figure_D_3-b1R} &  &  \\
     & $q_0q_{j+1}q_{j+7}\cdot h_0$\tikzmark{part2-figure_D_3-b2R} &  &  \\
     & $q_1q_jq_{j+7}\cdot h_j$\tikzmark{part2-figure_D_3-b3R} &  &  \\\hline
% r = 2j+7
    $2j+7$ &  & \tikzmark{part2-figure_D_3-n7L}$(7)\; q_1q_jq_{j+6}\cdot h_jh_{j+6}$\tikzmark{part2-figure_D_3-n7R} &  \\
     &  & \tikzmark{part2-figure_D_3-n8L}$(8)\; q_0q_{j+1}q_{j+6}\cdot h_0h_{j+6}$\tikzmark{part2-figure_D_3-n8R} &  \\
     &  & \tikzmark{part2-figure_D_3-n9L}$(9)\; q_0q_jq_{j+7}\cdot h_0h_j$\tikzmark{part2-figure_D_3-n9R} &  \\\hline
% r = 2j+6
    $2j+6$ &  &  & \tikzmark{part2-figure_D_3-t7L}$q_0q_jq_{j+6}\cdot h_0h_jh_{j+6}$ \\\hline
% r = \vdots
    $\vdots$ &  &  &  \\\hline
% r = 2j+1
    $2j+1$ &  & $(10)\; q_1q_j^2\cdot h_{j+6}^2$\tikzmark{part2-figure_D_3-n10R} &  \\\hline
% r = 2j
    $2j$ &  & $(11)\; q_0q_j^2\cdot h_0h_{j+7}$\tikzmark{part2-figure_D_3-n11R} & \tikzmark{part2-figure_D_3-t10L}$q_0q_j^2\cdot h_0h_{j+6}^2$ \\\hline
% r = 2j-1
    $2j-1$ &  &  &  \\\hline
% r = 2j-2
    $2j-2$ &  &  & \tikzmark{part2-figure_D_3-t11L}$q_0q_{j-1}^2\cdot h_0h_jh_{j+7}$ \\\hline
% r = \vdots
    $\vdots$ &  &  &  \\\hline
% r = j+7
    $j+7$ & $q_0^2q_{j+7}\cdot h_{j+1}$\tikzmark{part2-figure_D_3-a1R} & $(12)\; q_0^2q_{j+7}\cdot h_j^2$\tikzmark{part2-figure_D_3-n12R} &  \\\hline
% r = j+6
    $j+6$ &  & \tikzmark{part2-figure_D_3-t14L}$q_0^2q_{j+6}\cdot h_{j+1}h_{j+6}$ & \tikzmark{part2-figure_D_3-t12L}$q_0^2q_{j+6}\cdot h_j^2h_{j+6}$ \\\hline
% r = \vdots
    $\vdots$ &  &  &  \\\hline
% r = j+1
    $j+1$ & $q_0^2q_{j+1}\cdot h_{j+7}$\tikzmark{part2-figure_D_3-a2R} & $(13)\; q_0^2q_{j+1}\cdot h_{j+6}^2$\tikzmark{part2-figure_D_3-n13R} &  \\\hline
% r = j
    $j$ &  & \tikzmark{part2-figure_D_3-t15L}$q_0^2q_j\cdot h_jh_{j+7}$ & \tikzmark{part2-figure_D_3-t13L}$q_0^2q_j\cdot h_jh_{j+6}^2$ \\\hline
    \end{tabular}
\begin{tikzpicture}[overlay, remember picture, shorten >=.5pt, shorten <=.5pt, transform canvas={yshift=.25\baselineskip}]
    \draw [->] ({pic cs:part2-figure_D_3-n7R}) [red,bend left=20] to ({pic cs:part2-figure_D_3-t7L});
    \draw [->] ({pic cs:part2-figure_D_3-n8R}) [red] to ({pic cs:part2-figure_D_3-t7L});
    \draw [->] ({pic cs:part2-figure_D_3-n9R}) [red,bend right=18] to ({pic cs:part2-figure_D_3-t7L});
    \draw [->] ({pic cs:part2-figure_D_3-n10R}) [red] to ({pic cs:part2-figure_D_3-t10L});
    \draw [->] ({pic cs:part2-figure_D_3-n11R}) [blue] to ({pic cs:part2-figure_D_3-t11L});
    \draw [->] ({pic cs:part2-figure_D_3-n12R}) [red] to ({pic cs:part2-figure_D_3-t12L});
    \draw [->] ({pic cs:part2-figure_D_3-n13R}) [red] to ({pic cs:part2-figure_D_3-t13L});
    \draw [->] ({pic cs:part2-figure_D_3-a1R}) [red] to ({pic cs:part2-figure_D_3-t14L});
    \draw [->] ({pic cs:part2-figure_D_3-a2R}) [red] to ({pic cs:part2-figure_D_3-t15L});
    \draw [->] ({pic cs:part2-figure_D_3-b1R}) [red,out=-25,in=135] to ({pic cs:part2-figure_D_3-n7L});
    \draw [->] ({pic cs:part2-figure_D_3-b1R}) [red,out=-85,in=150] to ({pic cs:part2-figure_D_3-n8L});
    \draw [->] ({pic cs:part2-figure_D_3-b2R}) [red] to ({pic cs:part2-figure_D_3-n8L});
    \draw [->] ({pic cs:part2-figure_D_3-b2R}) [red,bend right=18] to ({pic cs:part2-figure_D_3-n9L});
    \draw [->] ({pic cs:part2-figure_D_3-b3R}) [red,bend right=24] to ({pic cs:part2-figure_D_3-n7L});
    \draw [->] ({pic cs:part2-figure_D_3-b3R}) [red] to ({pic cs:part2-figure_D_3-n9L});
\end{tikzpicture}
\par\endgroup
\caption{The BXSS in degree $(s,k,t) = (2,3,2^{j+2}\cdot 65+1)$ for \cref{CESS D_3}, where $j\geq 2$.
The vertical axis is the BX filtration and the horizontal axis is the $s$-degree. Red arrows indicate $d_1$-differentials, and blue arrows indicate $d_2$-differentials.}
\label{fig:D3-differentials}
\label{fig:D3-differentials-part2}
\end{figure}

    If $j=1$, the following classes and differentials changes:
    \[(1) d_2(q_2q_7^2\cdot h_0^2) = q_2q_6^2\cdot h_0^2h_7 \neq 0\]
    \[(7) \text{ is the same as } (3)\]
    \[(8) d_1(q_0q_2q_8\cdot h_0)= q_0q_2q_7\cdot h_0h_7\]
    \[(9) q_0q_1q_8\cdot h_0h_1=0\]
    \[(11) d_1(q_1^3\cdot h_8)=q_0q_1^2\cdot h_0h_8\]

\begin{figure}[H]
\centering
\begingroup
\setlength{\tabcolsep}{6pt}
\renewcommand{\arraystretch}{1.1}
\begin{tabular}{|wc{1.2cm}||wc{2.5cm}|p{5cm}|wc{4cm}|} \hline
      $w$ & $s=1$ & \multicolumn{1}{p{5cm}|}{\centering $s=2$} & $s=3$ \\\hline\hline
% r = 16
    $16$ &  & $(1)\; q_2q_7^2\cdot h_0^2$\tikzmark{part1-figure_D_3_j1-v1R} &  \\\hline
% r = 15
    $15$ &  & $(2,5)\; q_1q_7^2\cdot h_1^2$\tikzmark{part1-figure_D_3_j1-v3R} &  \\\hline
% r = 14
    $14$ &  & $(6)\; q_0q_7^2\cdot h_0h_2$\tikzmark{part1-figure_D_3_j1-v10R} & \tikzmark{part1-figure_D_3_j1-v2L}$q_2q_6^2\cdot h_0^2h_7$ \\
     &  &  & \tikzmark{part1-figure_D_3_j1-v9L}$q_0q_7^2\cdot h_1^3$ \\\hline
% r = 13
    $13$ &  &  & \tikzmark{part1-figure_D_3_j1-v4L}$q_1q_6^2\cdot h_1^2h_7$ \\\hline
% r = 12
    $12$ &  &  & \tikzmark{part1-figure_D_3_j1-v11L}$q_0q_6^2\cdot h_0h_2h_7$ \\\hline
% r = 11
    $11$ &  &  &  \\\hline
% r = 10
    $10$ & $q_1q_2q_7\cdot h_7$\tikzmark{part1-figure_D_3_j1-v25R} & $(4)\; q_1^2q_8\cdot h_0^2$\tikzmark{part1-figure_D_3_j1-v7R} &  \\
     & $(8)\; q_0q_2q_8\cdot h_0$\tikzmark{part1-figure_D_3_j1-v12R} &  &  \\
     & $q_1^2q_8\cdot h_1$\tikzmark{part1-figure_D_3_j1-v26R} &  &  \\\hline
% r = 9
    $9$ &  & \tikzmark{part1-figure_D_3_j1-v5L}$(3,7)\; q_1^2q_7\cdot h_1h_7$\tikzmark{part1-figure_D_3_j1-v5R} & \tikzmark{part1-figure_D_3_j1-v8L}$q_1^2q_7\cdot h_0^2h_7$ \\
     &  & \tikzmark{part1-figure_D_3_j1-v13L}$q_0q_2q_7\cdot h_0h_7$ &  \\\hline
% r = 8
    $8$ & $q_0^2q_8\cdot h_2$\tikzmark{part1-figure_D_3_j1-v21R} & $(12)\; q_0^2q_8\cdot h_1^2$\tikzmark{part1-figure_D_3_j1-v18R} &  \\\hline
% r = 7
    $7$ &  & \tikzmark{part1-figure_D_3_j1-v22L}$q_0^2q_7\cdot h_2h_7$ & \tikzmark{part1-figure_D_3_j1-v6L}$q_0^2q_7\cdot h_1^2h_7$ \\\hline
% r = 6
    $6$ &  &  &  \\\hline
% r = 5
    $5$ &  &  &  \\\hline
% r = 4
    $4$ &  &  &  \\\hline
% r = 3
    $3$ & $(11)\; q_1^3\cdot h_8$\tikzmark{part1-figure_D_3_j1-v16R} & $(10)\; q_1^3\cdot h_7^2$\tikzmark{part1-figure_D_3_j1-v14R} &  \\\hline
% r = 2
    $2$ & $q_0^2q_2\cdot h_8$\tikzmark{part1-figure_D_3_j1-v23R} & \tikzmark{part1-figure_D_3_j1-v17L}$q_0q_1^2\cdot h_0h_8$ & \tikzmark{part1-figure_D_3_j1-v15L}$q_0q_1^2\cdot h_0h_7^2$ \\
     &  & $(13)\; q_0^2q_2\cdot h_7^2$\tikzmark{part1-figure_D_3_j1-v19R} &  \\\hline
% r = 1
    $1$ &  & \tikzmark{part1-figure_D_3_j1-v24L}$q_0^2q_1\cdot h_1h_8$ & \tikzmark{part1-figure_D_3_j1-v20L}$q_0^2q_1\cdot h_1h_7^2$ \\\hline
    \end{tabular}
\begin{tikzpicture}[overlay, remember picture, shorten >=.5pt, shorten <=.5pt, transform canvas={yshift=.25\baselineskip}]

    \draw [->] ({pic cs:part1-figure_D_3_j1-v1R}) [blue] to ({pic cs:part1-figure_D_3_j1-v2L});
    \draw [->] ({pic cs:part1-figure_D_3_j1-v7R}) [red] to ({pic cs:part1-figure_D_3_j1-v8L});
    \draw [->] ({pic cs:part1-figure_D_3_j1-v3R}) [red] to ({pic cs:part1-figure_D_3_j1-v9L});
    \draw [->] ({pic cs:part1-figure_D_3_j1-v10R}) [blue] to ({pic cs:part1-figure_D_3_j1-v11L});
    \draw [->] ({pic cs:part1-figure_D_3_j1-v12R}) [red] to ({pic cs:part1-figure_D_3_j1-v13L});
    \draw [->] ({pic cs:part1-figure_D_3_j1-v14R}) [red] to ({pic cs:part1-figure_D_3_j1-v15L});
    \draw [->] ({pic cs:part1-figure_D_3_j1-v16R}) [red] to ({pic cs:part1-figure_D_3_j1-v17L});
    \draw [->] ({pic cs:part1-figure_D_3_j1-v18R}) [red] to ({pic cs:part1-figure_D_3_j1-v6L});
    \draw [->] ({pic cs:part1-figure_D_3_j1-v19R}) [red] to ({pic cs:part1-figure_D_3_j1-v20L});
    \draw [->] ({pic cs:part1-figure_D_3_j1-v21R}) [red] to ({pic cs:part1-figure_D_3_j1-v22L});
    \draw [->] ({pic cs:part1-figure_D_3_j1-v23R}) [red] to ({pic cs:part1-figure_D_3_j1-v24L});
    \draw [->] ({pic cs:part1-figure_D_3_j1-v25R}) [red] to ({pic cs:part1-figure_D_3_j1-v5L});
    \draw [->] ({pic cs:part1-figure_D_3_j1-v25R}) [red,bend right=12] to ({pic cs:part1-figure_D_3_j1-v13L});
    \draw [->] ({pic cs:part1-figure_D_3_j1-v26R}) [red,bend right=22] to ({pic cs:part1-figure_D_3_j1-v5L});
\end{tikzpicture}
\par\endgroup
\caption{The BXSS in degree $(s,k,t) = (2,3,2^{3}\cdot 65+1)$ for \cref{CESS D_3}.
The vertical axis is the BX filtration and the horizontal axis is the $s$-degree. Red arrows indicate $d_1$-differentials, and blue arrows indicate $d_2$-differentials.}
\label{fig:D3-j1-differentials}
\end{figure}

\end{proof}

\subsection{The $j=1$ Case}
We now show that the differential
\[
d_2(D_3(2)) = h_0K_1
\]
is nonzero. To see this, note that the charts produced by Lin-Wang-Xu (\cite{LWX}, \cite{Weinan}) show that the bidegree containing $K_1$, namely stem 255 and Adams filtration 5, is generated as an $\mathbb{F}_2$-vector space by
\[
\begin{array}{c c c c c}
h_0^4h_8 &&\text{and} &&x_{255,5},
\end{array}
\]
where $x_{255,5}$ is one of $K_1$ or $h_0^4h_8+K_1$. We also find that
\[
h_0\cdot h_0^4h_8 \neq h_0\cdot x_{255,5}
\]
with both $h_0$-extensions being nontrivial.

If $x_{255,5} = K_1$, then the claim follows immediately. Otherwise, using the inequalities described above, we have that
\[
h_0x_{255,5} = h_0^5h_8+h_0K_1
\]
is not equal to either of 0 or $h_0^5h_8$. If $h_0K_1 = 0$, this leads to a contradiction. Then we find that, regardless of the precise value of $x_{255,5}$, the target of the differential
\[
d_2(D_3(2)) = h_0K_1
\]
is nontrivial, as claimed.

\section{The $p_j'$ Family}

We prove the following theorem:
\begin{thm}
There are Adams differentials
\[
d_2(p_{j+2}') = h_0T_{j+1} \neq 0
\]
for $j\geq2$.
\end{thm}

Recall that the form of the differential was shown to hold in \Cref{sec:differentials}, so all that remains is to show that the target is nonzero on the Adams $E_2$ page. The proof proceeds in two steps, following the program of Burklund-Xu. First we show, for $j\geq2$, that the element $q_0\cdot T_j$ is a permanent cycle in the Burklund-Xu spectral sequence and detects a nonzero class in the abutment, an associated graded of the Cartan-Eilenberg $E_2$ page. We then show that $q_0\cdot T_j$ is a permanent cycle in that spectral sequence and detects a nonzero class in the abutment, an associated graded of the Adams $E_2$ page. Since $T_j$ detects the classical Adams element $T_{j+1}$ under the doubling-up isomorphism
\[
\text{Ext}_\mathcal{A}\xrightarrow{\cong}\text{Ext}_\mathcal{P},
\]
\cite{Benson} shows that $q_0\cdot T_j$ must detect $h_0T_{j+1}$. As a result, we will have shown that $h_0T_j$ is nonzero on the Adams $E_2$ page for $j\geq3$.

\subsection{Burklund-Xu Differentials}

We first show the following:
\begin{prop}
The Burklund-Xu element $q_0\cdot T_j$ is a permanent cycle for $j\geq0$.
\end{prop}
\begin{proof}
The Burklund-Xu degree of $q_0\cdot T_j$ is 0. Since Burklund-Xu differentials strictly lower filtration, this class cannot support any differential in that spectral sequence.
\end{proof}

Then the class $q_0\cdot T_j$ is defined in the abutment of the Burklund-Xu spectral sequence for all $j\geq0$, so it remains to show that the $j\geq2$ members of the family cannot receive any differential. We use the following simple proposition to simplify this argument:

\begin{prop}\label{prop:deg-pj'}
If $y\in\text{Ext}^4_\mathcal{P}$ has internal degree $t(y)$ divisible by $2^2$, then any Burklund-Xu element $q_a\cdot y$ in the same internal degree as $q_0\cdot T_j$ cannot support a differential.
\end{prop}
\begin{proof}
If $q_a\cdot y$ lies in the same internal degree as $q_0\cdot T_j$, then the equality
\[
t(q_a\cdot y) = 2^{a+1}-1+t(y) = 1+2^{j+2}+2^{j+5}+2^{j+8} = t(q_0\cdot T_j)
\]
must hold. Adding 1 to both sides and dividing by 2 gives
\[
2^a + 2^{-1}t(y) = 1+2^{j+1}+2^{j+4}+2^{j+7},
\]
and, by assumption, the term $2^{-1}t(y)$ is still divisible by 2. Then the only possibility is that $a=0$, meaning that
\[
q_a\cdot y = q_0\cdot y
\]
is the element in question.
\end{proof}

We now show that $q_0\cdot T_j$ cannot be the target of any Burklund-Xu differential:
\begin{prop}
There is no Burklund-Xu differential of the form
\[
d_r(x) = q_0\cdot T_j
\]
for any $j\geq2$ and Burklund-Xu element $x$.
\end{prop}
\begin{proof}
Recall that we can write $x$ as a linear combination of elements of the form $q_a\cdot y$ for $y\in\text{Ext}_\mathcal{P}^4$. As a result, we first collect all classes $q_a\cdot y$ having the same internal degree as $q_0\cdot T_j$ and then show that no linear combination of these classes can support a differential as in the proposition statement.

We begin with the case where $y$ is indecomposable (see \cite{Lin}, \cite{Chen} for a full presentation of the classical Adams 4-line). Of the possible cases, we find that if $y\in \{d_b,f_b,g_b\}$, then \Cref{prop:deg-pj'} implies that $q_a\cdot y = q_0\cdot y$. A class of this form cannot support any Burklund-Xu differential since it already lies in the lowest possible filtration, so we discard these classes.

We now consider $y = e_b$ so that the equality
\[
t(q_a\cdot e_b) = 2^{a+1}-1+2^{b+1}+2^{b+3}+2^{b+5} = 2^0+2^{j+2}+2^{j+5}+2^{j+8} = t(q_0\cdot T_j)
\]
holds. Adding 1 to both sides and dividing by 2 gives
\[
2^a+2^b+2^{b+2}+2^{b+4} = 2^0+2^{j+1}+2^{j+4}+2^{j+7},
\]
and since we must avoid $a=0$ by the reasoning above, we conclude that $b=0$. After subtracting 1 from both sides, this gives
\[
2^a+2^2+2^4 = 2^{j+1}+2^{j+4}+2^{j+7},
\]
but there is no choice of $a$ and $j$ compatible with this relation. Then $y\neq e_b$.

We now consider $y = p_b$ so that the equality
\[
t(q_a\cdot p_b) = 2^{a+1}-1+2^{b+1}+2^{b+3}+2^{b+6} = 2^0+2^{j+2}+2^{j+5}+2^{j+8} = t(q_0\cdot T_j)
\]
holds. Adding 1 to both sides and dividing by 2 gives
\[
2^a+2^b+2^{b+2}+2^{b+5} = 2^0+2^{j+1}+2^{j+4}+2^{j+7},
\]
and since we must avoid $a=0$ by the reasoning above, we conclude that $b=0$. After subtracting 1 from both sides, this gives
\[
2^a+2^2+2^5 = 2^{j+1}+2^{j+4}+2^{j+7},
\]
which is possible only if $j=1$ and $a = j+7 = 8$. However, we are considering only $j\geq2$, so we discard this case. Then $y \neq p_b$.

We no consider $y = D_3(b)$ so that the equality
\[
t(q_a\cdot D_3(b)) = 2^{a+1}-1+2^{b+1}+2^{b+7} = 2^0+2^{j+2}+2^{j+5}+2^{j+8} = t(q_0\cdot T_j)
\]
holds. Adding 1 to both sides and dividing by 2 gives
\[
2^a+2^b+2^{b+6} = 2^0+2^{j+1}+2^{j+4}+2^{j+7},
\]
and since we must avoid $a=0$ by the reasoning above, we conclude that $b=0$. After subtracting 1 from both sides, this gives
\[
2^a+2^6 = 2^{j+1}+2^{j+4}+2^{j+7},
\]
but there simply are not enough powers of 2 on the left-hand side to match the right-hand side. Then $y \neq D_3(b)$.

We now consider $y = p_b'$ so that the equality
\[
t(q_a\cdot p_b') = 2^{a+1}-1+2^{b+1}+2^{b+4}+2^{b+7} = 2^0+2^{j+2}+2^{j+5}+2^{j+8} = t(q_0\cdot T_j)
\]
holds. Adding 1 to both sides and dividing by 2 gives
\[
2^a+2^b+2^{b+3}+2^{b+6} = 2^0+2^{j+1}+2^{j+4}+2^{j+7},
\]
and since we must avoid $a=0$ by the reasoning above, we conclude that $b=0$. After subtracting 1 from both sides, this gives
\[
2^a+2^3+2^6 = 2^{j+1}+2^{j+4}+2^{j+7},
\]
which is possible only if $j=2$ and $a = j+7 = 9$. Then we will consider below the possibility that $q_a\cdot p_b' = q_9\cdot p_0'$, and we have conclude our search for indecomposable $y$.

The decomposable $y$ are of the form $h_\ell c_m$ and $h_bh_ch_dh_e$. We first consider $y = h_\ell c_m$ so that the equality
\[
2^a+2^\ell+2^m+2^{m+1}+2^{m+3} = 2^0+2^{j+1}+2^{j+4}+2^{j+7}
\]
holds after adding 1 and dividing by 2. In order to avoid $a=0$, we must have $\ell=0$ or $m=0$. If $\ell=0$, then we obtain
\[
2^a+2^m+2^{m+1}+2^{m+3} = 2^{j+1}+2^{j+4}+2^{j+7}
\]
after subtracting 1. Then $a\in\{m,m+1,m+3\}$. If $a=m$, then we obtain
\[
2^{m+2}+2^{m+3} = 2^{j+1}+2^{j+4}+2^{j+7},
\]
but there simply are not enough powers of 2 on the left-hand side to match the right-hand side. If $a=m+1$, then we obtain
\[
2^m+2^{m+2}+2^{m+3} = 2^{j+1}+2^{j+4}+2^{j+7},
\]
but then the smallest two powers of 2 on the left differ by two factors of 2, whereas the smallest two powers on the right differ by three factors of 2. Finally, if $a=m+3$, then we obtain
\[
2^m+2^{m+1}+2^{m+4} = 2^{j+1}+2^{j+4}+2^{j+7},
\]
but just as in the $a=m+1$ case, we see that the powers of 2 on the left cannot possibly match those on the right.

If instead $m=0$, then we obtain
\[
2^a+2^\ell+2^1+2^3 = 2^{j+1}+2^{j+4}+2^{j+7}
\]
after subtracting 1. Since we must avoid $a=0$ and have already explored the possibility that $\ell=0$, it must be that $a,\ell\geq1$. Then the smallest power on the right, $2^{j+1}$, must match the smallest power on the left, $2^1$, so that $j=0$. Since we consider only $j\geq2$, we disregard this case. Then $y\neq h_\ell c_m$.

Finally, we consider $y = h_bh_ch_dh_e$ so that the equality
\[
2^a+2^b+2^c+2^d+2^e = 2^0+2^{j+1}+2^{j+4}+2^{j+7}
\]
holds after adding 1 and dividing by 2. Since $\text{Ext}_\mathcal{P}$ is commutative, we may assume that $b\leq c\leq d\leq e$. Since we must avoid $a=0$, we then conclude that $b=0$ so that the equality
\[
2^a+2^c+2^d+2^e = 2^{j+1}+2^{j+4}+2^{j+7}
\]
holds after subtracting 1. In order for the four powers of 2 on the left to match the three powers of 2 on the right, two of $\{a,c,d,e\}$ must be equal. By the assumption that $c\leq d\leq e$, it then suffices to consider the following cases:
\[
\begin{array}{c c c}
a=\ell=j, m=j+4, n=j+7 &&\ell=j+1, a=m=j+3, n=j+7 \\
&\ell=j+1, m=j+4, a=n=j+6 &
\end{array}
\]
These cases correspond, respectively, to the following elements:
\[
\begin{array}{c c c}
q_j\cdot h_0h_jh_{j+4}h_{j+7} &&q_{j+3}\cdot h_0h_{j+1}h_{j+3}h_{j+7} \\
&q_{j+6}\cdot h_0h_{j+1}h_{j+4}h_{j+6} &
\end{array}
\]
These elements are, respectively, the result of evaluating $d_1$ on
\[
\begin{array}{c c c}
q_{j+1}\cdot h_0h_{j+4}h_{j+7} &&q_{j+4}\cdot h_0h_{j+1}h_{j+7} \\
&q_{j+7}\cdot h_0h_{j+1}h_{j+4}, &
\end{array}
\]
so we discard them. Then $y\neq h_bh_ch_dh_e$.

In summary, we have collected only the class
\[
q_9\cdot p_0',
\]
but it is easily seen to support the differential
\[
d_1(q_9\cdot p_0') = q_8\cdot h_8p_0',
\]
by the presentation in \cite{Chen}. Then no Burklund-Xu class can possibly support a differential with target
\[
q_0\cdot T_j
\]
for any $j\geq2$. %add j=1 case
\end{proof}

\subsection{Cartan-Eilenberg Differentials}

\begin{prop}
\label{CESS p prime}
The class 
\[q_0\cdot T_j\]
survives in the Catan-Eilenberg spectral sequence for $j\ge 1$.
\end{prop}
\begin{proof}
    $q_0\cdot T_j$ has CE filtration 1, so it cannot support a CE differential, which lower the CE filtration by 2.\\
    \\
    The source of differentials that could hit $q_o\cdot T_j$ must have CE filtration $k=3$ or $k=5$.\\
    If $k=5$, then the only possible class in degree $(s,k,t)= (0,5,2^{j+8}+2^{j+5}+2^{j+2}+1)$ is $q_0q_1q_{j+1}q_{j+4}q_{j+7}$, which does not survive to the Cartan-Eilenberg spectral sequence because it supports a Burklund-Xu differential:
    \[d_1(q_0q_1q_{j+1}q_{j+4}q_{j+7}) = q_0^2q_{j+1}q_{j+4}q_{j+7}\cdot h_0+q_0q_1q_{j}q_{j+4}q_{j+7}h_j+q_0q_1q_{j+1}q_{j+3}q_{j+7}h_{j+3}+q_0q_1q_{j+1}q_{j+4}q_{j+6}h_{j+6}\]
    If $k=3$, then the possible classes in degree $(2,3,2^{j+8}+2^{j+5}+2^{j+2}+1)$ all support non-trivial Burklund-Xu differentials:
    \begin{enumerate}
    \item $d_1(q_{j+4}q_{j+6}^2\cdot h_1h_{j+1}) = q_{j+3}q_{j+6}^2\cdot h_1h_{j+1}h_{j+3}$
    \item $d_1(q_{j+1}q_{j+6}^2\cdot h_1h_{j+4}) = q_jq_{j+6}^2\cdot h_1h_jh_{j+4}$
    \item $d_1(q_{j+3}^2q_{j+7}\cdot h_1h_{j+1})=q_{j+3}^2q_{j+6}\cdot h_1h_{j+1}h_{j+6}$
    \item $d_1(q_{j+1}q_{j+4}q_{j+7}\cdot h_0^2)=q_jq_{j+4}q_{j+7}\cdot h_0^2h_j+q_{j+1}q_{j+3}q_{j+7}\cdot h_0^2h_{j+3}+q_{j+1}q_{j+4}q_{j+6}\cdot h_0^2h_{j+6}$
    \item $d_1(q_{j+1}q_{j+4}q_{j+6}\cdot h_1h_{j+6}) = q_jq_{j+4}q_{j+6}\cdot h_1h_jh_{j+6} + q_{j+1}q_{j+3}q_{j+6}\cdot h_1h_{j+3}h_{j+6}$
    \item $d_1(q_{j+1}q_{j+3}q_{j+7}\cdot h_1h_{j+3}) = q_jq_{j+3}q_{j+7}\cdot h_1h_jh_{j+3}+q_{j+1}q_{j+3}q_{j+6}\cdot h_1h_{j+3}h_{j+6}$
    \item $d_1(q_jq_{j+4}q_{j+7}\cdot h_1h_j) = q_jq_{j+4}q_{j+6}\cdot h_1h_jh_{j+6}+q_jq_{j+3}q_{j+7}\cdot h_1h_jh_{j+3}$
    \item $d_1(q_{j+1}q_{j+3}^2\cdot h_1h_{j+7}) = q_jq_{j+3}^2\cdot h_1h_jh_{j+7}$
    \item $d_1(q_j^2q_{j+7}\cdot h_1h_{j+4}) = q_j^2q_{j+6}\cdot h_1h_{j+4}h_{j+6}$
    \item $d_1(q_j^2q_{j+4}\cdot h_1h_{j+7})=q_j^2q_{j+3}\cdot h_1h_{j+3}h_{j+7}$
    \item $d_1(q_1q_{j+6}^2\cdot h_{j+1}h_{j+4}) = q_0q_{j+6}^2\cdot h_0h_{j+1}h_{j+4}$
    \item $d_1(q_1q_{j+4}q_{j+7}\cdot h_{j}^2) = q_1q_{j+3}q_{j+7}\cdot h_{j}^2h_{j+3}+q_1q_{j+4}q_{j+6}\cdot h_{j}^2h_{j+6}+q_0q_{j+4}q_{j+7}\cdot h_0 h_{j}^2$
    \item $d_1(q_1q_{j+4}q_{j+6}\cdot h_{j+1}h_{j+6}) = q_1q_{j+3}q_{j+6}\cdot h_{j+1}h_{j+3}h_{j+6}+q_0q_{j+4}q_{j+6}\cdot h_0h_{j+1}h_{j+6}$
    \item $d_1(q_1q_{j+3}q_{j+7}\cdot h_{j+1}h_{j+3}) = q_1q_{j+3}q_{j+6}\cdot h_{j+1}h_{j+3}h_{j+6}+q_0q_{j+3}q_{j+7}\cdot h_0h_{j+1}h_{j+3}$
    \item $d_1(q_0q_{j+4}q_{j+7}\cdot h_0h_{j+1})=q_0q_{j+4}q_{j+6}\cdot h_0h_{j+1}h_{j+6}+q_0q_{j+3}q_{j+7}\cdot h_0h_{j+1}h_{j+3}$
    \item $d_1(q_1q_{j+1}q_{j+7}\cdot h_{j+3}^2) = q_1q_jq_{j+7}\cdot h_jh_{j+3}^2+q_1q_{j+1}q_{j+6}\cdot h_{j+3}^2h_{j+6}+q_0q_{j+1}q_{j+7}\cdot h_0h_{j+3}^2$
    \item $d_1(q_1q_{j+1}q_{j+6}\cdot h_{j+4}h_{j+6}) = q_1q_jq_{j+6}\cdot h_jh_{j+4}h_{j+6}+q_0q_{j+1}q_{j+6}\cdot h_0h_{j+4}h_{j+6}$
    \item $d_1(q_1q_jq_{j+7}\cdot h_jh_{j+4}) = q_1q_jq_{j+6}\cdot h_jh_{j+4}h_{j+6}+q_0q_jq_{j+7}\cdot h_0h_jh_{j+4}$
    \item $d_1(q_0q_{j+1}q_{j+7}\cdot h_0h_{j+4})=q_0q_{j+1}q_{j+6}\cdot h_0h_{j+4}h_{j+6}+q_0q_jq_{j+7}\cdot h_0h_jh_{j+4}$
    \item $d_1(q_1q_{j+3}^2\cdot h_{j+1}h_{j+7}) = q_0q_{j+3}^2\cdot h_0h_{j+1}h_{j+7}$
    \item $d_1(q_1q_{j+1}q_{j+5}\cdot h_{j+6}^2) = q_1q_jq_{j+5}\cdot h_jh_{j+6}^2+q_1q_{j+1}q_{j+4}\cdot h_{j+4}h_{j+6}^2+q_0q_{j+1}q_{j+5}\cdot h_0h_{j+6}^2$
    \item $d_1(q_1q_{j+1}q_{j+3}\cdot h_{j+3}h_{j+7}) = q_1q_jq_{j+3}\cdot h_jh_{j+3}h_{j+7}+q_0q_{j+1}q_{j+3}\cdot h_0h_{j+3}h_{j+7}$
    \item $d_1(q_1q_jq_{j+4}\cdot h_jh_{j+7}) = q_1q_jq_{j+3}\cdot h_jh_{j+3}h_{j+7}+q_0q_jq_{j+4}\cdot h_0h_jh_{j+7}$
    \item $d_1(q_0q_{j+1}q_{j+4}\cdot h_0h_{j+7})=q_0q_{j+1}q_{j+3}\cdot h_0h_{j+3}h_{j+7}+q_0q_jq_{j+4}\cdot h_0h_jh_{j+7}$
    \item $d_1(q_1q_j^2\cdot h_{j+4}h_{j+7}) = q_0q_j^2\cdot h_0h_{j+4}h_{j+7}$
    \item $d_1(q_0^2q_{j+7}\cdot h_{j+1}h_{j+4})=q_0^2q_{j+6}\cdot h_{j+1}h_{j+4}h_{j+6}$
    \item $d_1(q_0^2q_{j+4}\cdot h_{j+1}h_{j+7})=q_0^2q_{j+3}\cdot h_{j+1}h_{j+3}h_{j+7}$
    \item $d_1(q_0^2q_{j+1}\cdot h_{j+4}h_{j+7})=q_0^2q_{j}\cdot h_{j}h_{j+4}h_{j+7}$
    \end{enumerate}
    Note there are linear relations among targets of these differentials, in particular, we should account for
    \begin{enumerate}
        \item $d_1(q_{j+1}q_{j+4}q_{j+6}\cdot h_1h_{j+6}+q_{j+1}q_{j+3}q_{j+7}\cdot h_1h_{j+3}+q_jq_{j+4}q_{j+7}\cdot h_1h_j) = 0$
        \item $d_1(q_1q_{j+4}q_{j+6}\cdot h_{j+1}h_{j+6}+q_1q_{j+3}q_{j+7}\cdot h_{j+1}h_{j+3}+q_0q_{j+4}q_{j+7}\cdot h_0h_{j+1}) = 0$
        \item $d_1(q_1q_{j+1}q_{j+6}\cdot h_{j+4}h_{j+6}+q_1q_jq_{j+7}\cdot h_jh_{j+4}+q_0q_{j+1}q_{j+7}\cdot h_0h_{j+4})= 0$
        \item $d_1(q_1q_{j+1}q_{j+3}\cdot h_{j+3}h_{j+7}+q_1q_jq_{j+5}\cdot h_jh_{j+7}+q_0q_{j+1}q_{j+4}\cdot h_0h_{j+7})=0$
    \end{enumerate}
    by finding differentials killing the sources. Indeed, we have
    \begin{enumerate}
        \item $d_1(q_{j+1}q_{j+4}q_{j+7}\cdot h_1) = q_{j+1}q_{j+4}q_{j+6}\cdot h_1h_{j+6}+q_{j+1}q_{j+3}q_{j+7}\cdot h_1h_{j+3}+q_jq_{j+4}q_{j+7}\cdot h_1h_j$
        \item $d_1(q_1q_{j+4}q_{j+7}\cdot h_{j+1}) = q_1q_{j+4}q_{j+6}\cdot h_{j+1}h_{j+6}+q_1q_{j+3}q_{j+7}\cdot h_{j+1}h_{j+3}+q_0q_{j+4}q_{j+7}\cdot h_0h_{j+1}$
        \item $d_1(q_1q_{j+1}q_{j+7}\cdot h_{j+4}) = q_1q_{j+1}q_{j+6}\cdot h_{j+4}h_{j+6}+q_1q_jq_{j+7}\cdot h_jh_{j+4}+q_0q_{j+1}q_{j+7}\cdot h_0h_{j+4}$
        \item $d_1(q_1q_{j+1}q_{j+4}\cdot h_{j+7}) = q_1q_{j+1}q_{j+3}\cdot h_{j+3}h_{j+7}+q_1q_jq_{j+5}\cdot h_jh_{j+7}+q_0q_{j+1}q_{j+4}\cdot h_0h_{j+7}$
    \end{enumerate}
    For $j\geq 3$, all the classes above are distinct and nonzero.

\begin{figure}[H]
\centering
\begingroup
\setlength{\tabcolsep}{6pt}
\renewcommand{\arraystretch}{1.1}
\begin{tabular}{|wc{1.2cm}||wc{2.5cm}|p{5cm}|wc{4cm}|} \hline
      $w$ & $s=1$ & \multicolumn{1}{p{5cm}|}{\centering $s=2$} & $s=3$ \\\hline\hline
% r = 3j+16
    $3j+16$ &  & \tikzmark{part1-figure_p_prime-p1L}$(1)\; q_{j+4}q_{j+6}^2\cdot h_1h_{j+1}$\tikzmark{part1-figure_p_prime-p1R} &  \\\hline
% r = 3j+15
    $3j+15$ &  &  & \tikzmark{part1-figure_p_prime-p2L}$q_{j+3}q_{j+6}^2\cdot h_1h_{j+1}h_{j+3}$\tikzmark{part1-figure_p_prime-p2R} \\\hline
% r = 3j+14
    $3j+14$ &  &  &  \\\hline
% r = 3j+13
    $3j+13$ &  & \tikzmark{part1-figure_p_prime-p3L}$(2)\; q_{j+1}q_{j+6}^2\cdot h_1h_{j+4}$\tikzmark{part1-figure_p_prime-p3R} &  \\
     &  & \tikzmark{part1-figure_p_prime-p5L}$(3)\; q_{j+3}^2q_{j+7}\cdot h_1h_{j+1}$\tikzmark{part1-figure_p_prime-p5R} &  \\\hline
% r = 3j+12
    $3j+12$ & \tikzmark{part1-figure_p_prime-p65L}$q_{j+1}q_{j+4}q_{j+7}\cdot h_1$\tikzmark{part1-figure_p_prime-p65R} & \tikzmark{part1-figure_p_prime-p7L}$(4)\; q_{j+1}q_{j+4}q_{j+7}\cdot h_0^2$\tikzmark{part1-figure_p_prime-p7R} & \tikzmark{part1-figure_p_prime-p4L}$q_jq_{j+6}^2\cdot h_1h_jh_{j+4}$\tikzmark{part1-figure_p_prime-p4R} \\
     &  &  & \tikzmark{part1-figure_p_prime-p6L}$q_{j+3}^2q_{j+6}\cdot h_1h_{j+1}h_{j+6}$\tikzmark{part1-figure_p_prime-p6R} \\\hline
% r = 3j+11
    $3j+11$ &  & \tikzmark{part1-figure_p_prime-p11L}$(5)\; q_{j+1}q_{j+4}q_{j+6}\cdot h_1h_{j+6}$\tikzmark{part1-figure_p_prime-p11R} & \tikzmark{part1-figure_p_prime-p8L}$q_jq_{j+4}q_{j+7}\cdot h_0^2h_j$\tikzmark{part1-figure_p_prime-p8R} \\
     &  & \tikzmark{part1-figure_p_prime-p14L}$(6)\; q_{j+1}q_{j+3}q_{j+7}\cdot h_1h_{j+3}$\tikzmark{part1-figure_p_prime-p14R} & \tikzmark{part1-figure_p_prime-p9L}$q_{j+1}q_{j+3}q_{j+7}\cdot h_0^2h_{j+3}$\tikzmark{part1-figure_p_prime-p9R} \\
     &  & \tikzmark{part1-figure_p_prime-p16L}$(7)\; q_jq_{j+4}q_{j+7}\cdot h_1h_j$\tikzmark{part1-figure_p_prime-p16R} & \tikzmark{part1-figure_p_prime-p10L}$q_{j+1}q_{j+4}q_{j+6}\cdot h_0^2h_{j+6}$\tikzmark{part1-figure_p_prime-p10R} \\\hline
% r = 3j+10
    $3j+10$ &  &  & \tikzmark{part1-figure_p_prime-p12L}$q_jq_{j+4}q_{j+6}\cdot h_1h_jh_{j+6}$\tikzmark{part1-figure_p_prime-p12R} \\
     &  &  & \tikzmark{part1-figure_p_prime-p13L}$q_{j+1}q_{j+3}q_{j+6}\cdot h_1h_{j+3}h_{j+6}$\tikzmark{part1-figure_p_prime-p13R} \\
     &  &  & \tikzmark{part1-figure_p_prime-p15L}$q_jq_{j+3}q_{j+7}\cdot h_1h_jh_{j+3}$\tikzmark{part1-figure_p_prime-p15R} \\\hline
% r = 3j+9
    $3j+9$ &  &  &  \\\hline
% r = 3j+8
    $3j+8$ &  &  &  \\\hline
% r = 3j+7
    $3j+7$ &  & \tikzmark{part1-figure_p_prime-p17L}$(8)\; q_{j+1}q_{j+3}^2\cdot h_1h_{j+7}$\tikzmark{part1-figure_p_prime-p17R} &  \\
     &  & \tikzmark{part1-figure_p_prime-p19L}$(9)\; q_j^2q_{j+7}\cdot h_1h_{j+4}$\tikzmark{part1-figure_p_prime-p19R} &  \\\hline
% r = 3j+6
    $3j+6$ &  &  & \tikzmark{part1-figure_p_prime-p18L}$q_jq_{j+3}^2\cdot h_1h_jh_{j+7}$\tikzmark{part1-figure_p_prime-p18R} \\
     &  &  & \tikzmark{part1-figure_p_prime-p20L}$q_j^2q_{j+6}\cdot h_1h_{j+4}h_{j+6}$\tikzmark{part1-figure_p_prime-p20R} \\\hline
% r = 3j+5
    $3j+5$ &  &  &  \\\hline
% r = 3j+4
    $3j+4$ &  & \tikzmark{part1-figure_p_prime-p21L}$(10)\; q_j^2q_{j+4}\cdot h_1h_{j+7}$\tikzmark{part1-figure_p_prime-p21R} &  \\\hline
% r = 3j+3
    $3j+3$ &  &  & \tikzmark{part1-figure_p_prime-p22L}$q_j^2q_{j+3}\cdot h_1h_{j+3}h_{j+7}$\tikzmark{part1-figure_p_prime-p22R} \\\hline
    $\vdots$ & & & \\\hline
% r = 2j+13
    $2j+13$ &  & \tikzmark{part1-figure_p_prime-p23L}$(11)\; q_1q_{j+6}^2\cdot h_{j+1}h_{j+4}$\tikzmark{part1-figure_p_prime-p23R} &  \\\hline
% r = 2j+12
    $2j+12$ & \tikzmark{part1-figure_p_prime-p66L}$q_1q_{j+4}q_{j+7}\cdot h_{j+1}$\tikzmark{part1-figure_p_prime-p66R} & \tikzmark{part1-figure_p_prime-p25L}$(12)\; q_1q_{j+4}q_{j+7}\cdot h_j^2$\tikzmark{part1-figure_p_prime-p25R} & \tikzmark{part1-figure_p_prime-p24L}$q_0q_{j+6}^2\cdot h_0h_{j+1}h_{j+4}$\tikzmark{part1-figure_p_prime-p24R} \\\hline
% r = 2j+11
    $2j+11$ &  & \tikzmark{part1-figure_p_prime-p29L}$(13)\; q_1q_{j+4}q_{j+6}\cdot h_{j+1}h_{j+6}$\tikzmark{part1-figure_p_prime-p29R} & \tikzmark{part1-figure_p_prime-p26L}$q_1q_{j+3}q_{j+7}\cdot h_j^2h_{j+3}$\tikzmark{part1-figure_p_prime-p26R} \\
     &  & \tikzmark{part1-figure_p_prime-p32L}$(14)\; q_1q_{j+3}q_{j+7}\cdot h_{j+1}h_{j+3}$\tikzmark{part1-figure_p_prime-p32R} & \tikzmark{part1-figure_p_prime-p27L}$q_1q_{j+4}q_{j+6}\cdot h_j^2h_{j+6}$\tikzmark{part1-figure_p_prime-p27R} \\
     &  & \tikzmark{part1-figure_p_prime-p34L}$(15)\; q_0q_{j+4}q_{j+7}\cdot h_0h_{j+1}$\tikzmark{part1-figure_p_prime-p34R} & \tikzmark{part1-figure_p_prime-p28L}$q_0q_{j+4}q_{j+7}\cdot h_0h_j^2$\tikzmark{part1-figure_p_prime-p28R} \\\hline
% r = 2j+10
    $2j+10$ &  &  & \tikzmark{part1-figure_p_prime-p30L}$q_1q_{j+3}q_{j+6}\cdot h_{j+1}h_{j+3}h_{j+6}$\tikzmark{part1-figure_p_prime-p30R} \\
     &  &  & \tikzmark{part1-figure_p_prime-p31L}$q_0q_{j+4}q_{j+6}\cdot h_0h_{j+1}h_{j+6}$\tikzmark{part1-figure_p_prime-p31R} \\
     &  &  & \tikzmark{part1-figure_p_prime-p33L}$q_0q_{j+3}q_{j+7}\cdot h_0h_{j+1}h_{j+3}$\tikzmark{part1-figure_p_prime-p33R} \\\hline
\end{tabular}
\begin{tikzpicture}[overlay, remember picture, shorten >=.5pt, shorten <=.5pt, transform canvas={yshift=.25\baselineskip}]

    \draw [->] ({pic cs:part1-figure_p_prime-p1R}) [red] to ({pic cs:part1-figure_p_prime-p2L});
    \draw [->] ({pic cs:part1-figure_p_prime-p3R}) [red] to ({pic cs:part1-figure_p_prime-p4L});
    \draw [->] ({pic cs:part1-figure_p_prime-p5R}) [red] to ({pic cs:part1-figure_p_prime-p6L});
    \draw [->] ({pic cs:part1-figure_p_prime-p7R}) [red,bend left=16] to ({pic cs:part1-figure_p_prime-p8L});
    \draw [->] ({pic cs:part1-figure_p_prime-p7R}) [red] to ({pic cs:part1-figure_p_prime-p9L});
    \draw [->] ({pic cs:part1-figure_p_prime-p7R}) [red,bend left=0] to ({pic cs:part1-figure_p_prime-p10L});
    \draw [->] ({pic cs:part1-figure_p_prime-p11R}) [red,bend left=12] to ({pic cs:part1-figure_p_prime-p12L});
    \draw [->] ({pic cs:part1-figure_p_prime-p11R}) [red,bend right=4] to ({pic cs:part1-figure_p_prime-p13L});
    \draw [->] ({pic cs:part1-figure_p_prime-p14R}) [red,bend left=4] to ({pic cs:part1-figure_p_prime-p13L});
    \draw [->] ({pic cs:part1-figure_p_prime-p14R}) [red,bend right=15] to ({pic cs:part1-figure_p_prime-p15L});
    \draw [->] ({pic cs:part1-figure_p_prime-p16R}) [red,bend left=4] to ({pic cs:part1-figure_p_prime-p12L});
    \draw [->] ({pic cs:part1-figure_p_prime-p16R}) [red,bend right=12] to ({pic cs:part1-figure_p_prime-p15L});
    \draw [->] ({pic cs:part1-figure_p_prime-p17R}) [red] to ({pic cs:part1-figure_p_prime-p18L});
    \draw [->] ({pic cs:part1-figure_p_prime-p19R}) [red] to ({pic cs:part1-figure_p_prime-p20L});
    \draw [->] ({pic cs:part1-figure_p_prime-p21R}) [red] to ({pic cs:part1-figure_p_prime-p22L});
    \draw [->] ({pic cs:part1-figure_p_prime-p23R}) [red] to ({pic cs:part1-figure_p_prime-p24L});
    \draw [->] ({pic cs:part1-figure_p_prime-p25R}) [red,bend left=16] to ({pic cs:part1-figure_p_prime-p26L});
    \draw [->] ({pic cs:part1-figure_p_prime-p25R}) [red,bend left=8] to ({pic cs:part1-figure_p_prime-p27L});
    \draw [->] ({pic cs:part1-figure_p_prime-p25R}) [red,out=0,in=160] to ({pic cs:part1-figure_p_prime-p28L});
    \draw [->] ({pic cs:part1-figure_p_prime-p29R}) [red,bend left=12] to ({pic cs:part1-figure_p_prime-p30L});
    \draw [->] ({pic cs:part1-figure_p_prime-p29R}) [red,bend right=4] to ({pic cs:part1-figure_p_prime-p31L});
    \draw [->] ({pic cs:part1-figure_p_prime-p32R}) [red,bend left=4] to ({pic cs:part1-figure_p_prime-p30L});
    \draw [->] ({pic cs:part1-figure_p_prime-p32R}) [red,bend right=4] to ({pic cs:part1-figure_p_prime-p33L});
    \draw [->] ({pic cs:part1-figure_p_prime-p34R}) [red,bend left=4] to ({pic cs:part1-figure_p_prime-p31L});
    \draw [->] ({pic cs:part1-figure_p_prime-p34R}) [red,bend right=12] to ({pic cs:part1-figure_p_prime-p33L});
    \draw [->] ({pic cs:part1-figure_p_prime-p65R}) [red,bend left=16] to ({pic cs:part1-figure_p_prime-p11L});
    \draw [->] ({pic cs:part1-figure_p_prime-p65R}) [red] to ({pic cs:part1-figure_p_prime-p14L});
    \draw [->] ({pic cs:part1-figure_p_prime-p65R}) [red,bend left=0] to ({pic cs:part1-figure_p_prime-p16L});
    \draw [->] ({pic cs:part1-figure_p_prime-p66R}) [red,bend left=16] to ({pic cs:part1-figure_p_prime-p29L});
    \draw [->] ({pic cs:part1-figure_p_prime-p66R}) [red] to ({pic cs:part1-figure_p_prime-p32L});
    \draw [->] ({pic cs:part1-figure_p_prime-p66R}) [red,bend left=0] to ({pic cs:part1-figure_p_prime-p34L});
\end{tikzpicture}
\par\endgroup
\end{figure}

\begin{figure}[H]
\centering
\begingroup
\setlength{\tabcolsep}{6pt}
\renewcommand{\arraystretch}{1.1}
\begin{tabular}{|wc{1.2cm}||wc{2.5cm}|p{5cm}|wc{4cm}|} \hline
      $w$ & $s=1$ & \multicolumn{1}{p{5cm}|}{\centering $s=2$} & $s=3$ \\\hline\hline
% r = 2j+9
    $2j+9$ & \tikzmark{part2-figure_p_prime-p67L}$q_1q_{j+1}q_{j+7}\cdot h_{j+4}$\tikzmark{part2-figure_p_prime-p67R} & \tikzmark{part2-figure_p_prime-p35L}$(16)\; q_1q_{j+1}q_{j+7}\cdot h_{j+3}^2$\tikzmark{part2-figure_p_prime-p35R} &  \\\hline
% r = 2j+8
    $2j+8$ &  & \tikzmark{part2-figure_p_prime-p39L}$(17)\; q_1q_{j+1}q_{j+6}\cdot h_{j+4}h_{j+6}$\tikzmark{part2-figure_p_prime-p39R} & \tikzmark{part2-figure_p_prime-p36L}$q_1q_jq_{j+7}\cdot h_jh_{j+3}^2$\tikzmark{part2-figure_p_prime-p36R} \\
     &  & \tikzmark{part2-figure_p_prime-p42L}$(18)\; q_1q_jq_{j+7}\cdot h_jh_{j+4}$\tikzmark{part2-figure_p_prime-p42R} & \tikzmark{part2-figure_p_prime-p37L}$q_1q_{j+1}q_{j+6}\cdot h_{j+3}^2h_{j+6}$\tikzmark{part2-figure_p_prime-p37R} \\
     &  & \tikzmark{part2-figure_p_prime-p44L}$(19)\; q_0q_{j+1}q_{j+7}\cdot h_0h_{j+4}$\tikzmark{part2-figure_p_prime-p44R} & \tikzmark{part2-figure_p_prime-p38L}$q_0q_{j+1}q_{j+7}\cdot h_0h_{j+3}^2$\tikzmark{part2-figure_p_prime-p38R} \\\hline
% r = 2j+7
    $2j+7$ &  & \tikzmark{part2-figure_p_prime-p45L}$(20)\; q_1q_{j+3}^2\cdot h_{j+1}h_{j+7}$\tikzmark{part2-figure_p_prime-p45R} & \tikzmark{part2-figure_p_prime-p40L}$q_1q_jq_{j+6}\cdot h_jh_{j+4}h_{j+6}$\tikzmark{part2-figure_p_prime-p40R} \\
     &  & \tikzmark{part2-figure_p_prime-p47L}$(21)\; q_1q_{j+1}q_{j+5}\cdot h_{j+6}^2$\tikzmark{part2-figure_p_prime-p47R} & \tikzmark{part2-figure_p_prime-p41L}$q_0q_{j+1}q_{j+6}\cdot h_0h_{j+4}h_{j+6}$\tikzmark{part2-figure_p_prime-p41R} \\
     &  &  & \tikzmark{part2-figure_p_prime-p43L}$q_0q_jq_{j+7}\cdot h_0h_jh_{j+4}$\tikzmark{part2-figure_p_prime-p43R} \\\hline
% r = 2j+6
    $2j+6$ & \tikzmark{part2-figure_p_prime-p68L}$q_1q_{j+1}q_{j+4}\cdot h_{j+7}$\tikzmark{part2-figure_p_prime-p68R} &  & \tikzmark{part2-figure_p_prime-p46L}$q_0q_{j+3}^2\cdot h_0h_{j+1}h_{j+7}$\tikzmark{part2-figure_p_prime-p46R} \\
     &  &  & \tikzmark{part2-figure_p_prime-p48L}$q_1q_jq_{j+5}\cdot h_jh_{j+6}^2$\tikzmark{part2-figure_p_prime-p48R} \\
     &  &  & \tikzmark{part2-figure_p_prime-p49L}$q_1q_{j+1}q_{j+4}\cdot h_{j+4}h_{j+6}^2$\tikzmark{part2-figure_p_prime-p49R} \\
     &  &  & \tikzmark{part2-figure_p_prime-p50L}$q_0q_{j+1}q_{j+5}\cdot h_0h_{j+6}^2$\tikzmark{part2-figure_p_prime-p50R} \\\hline
% r = 2j+5
    $2j+5$ &  & \tikzmark{part2-figure_p_prime-p51L}$(22)\; q_1q_{j+1}q_{j+3}\cdot h_{j+3}h_{j+7}$\tikzmark{part2-figure_p_prime-p51R} &  \\
     &  & \tikzmark{part2-figure_p_prime-p54L}$(23)\; q_1q_jq_{j+4}\cdot h_jh_{j+7}$\tikzmark{part2-figure_p_prime-p54R} &  \\
     &  & \tikzmark{part2-figure_p_prime-p56L}$(24)\; q_0q_{j+1}q_{j+4}\cdot h_0h_{j+7}$\tikzmark{part2-figure_p_prime-p56R} &  \\\hline
% r = 2j+4
    $2j+4$ &  &  & \tikzmark{part2-figure_p_prime-p52L}$q_1q_jq_{j+3}\cdot h_jh_{j+3}h_{j+7}$\tikzmark{part2-figure_p_prime-p52R} \\
     &  &  & \tikzmark{part2-figure_p_prime-p53L}$q_0q_{j+1}q_{j+3}\cdot h_0h_{j+3}h_{j+7}$\tikzmark{part2-figure_p_prime-p53R} \\
     &  &  & \tikzmark{part2-figure_p_prime-p55L}$q_0q_jq_{j+4}\cdot h_0h_jh_{j+7}$\tikzmark{part2-figure_p_prime-p55R} \\\hline
% r = 2j+3
    $2j+3$ &  &  &  \\\hline
% r = 2j+2
    $2j+2$ &  &  &  \\\hline
% r = 2j+1
    $2j+1$ &  & \tikzmark{part2-figure_p_prime-p57L}$(25)\; q_1q_j^2\cdot h_{j+4}h_{j+7}$\tikzmark{part2-figure_p_prime-p57R} &  \\\hline
% r = 2j
    $2j$ &  &  & \tikzmark{part2-figure_p_prime-p58L}$q_0q_j^2\cdot h_0h_{j+4}h_{j+7}$\tikzmark{part2-figure_p_prime-p58R} \\\hline
    $\vdots$ & & & \\\hline
% r = j+7
    $j+7$ &  & \tikzmark{part2-figure_p_prime-p59L}$(26)\; q_0^2q_{j+7}\cdot h_{j+1}h_{j+4}$\tikzmark{part2-figure_p_prime-p59R} &  \\\hline
% r = j+6
    $j+6$ &  &  & \tikzmark{part2-figure_p_prime-p60L}$q_0^2q_{j+6}\cdot h_{j+1}h_{j+4}h_{j+6}$\tikzmark{part2-figure_p_prime-p60R} \\\hline
% r = j+5
    $j+5$ &  &  &  \\\hline
% r = j+4
    $j+4$ &  & \tikzmark{part2-figure_p_prime-p61L}$(27)\; q_0^2q_{j+4}\cdot h_{j+1}h_{j+7}$\tikzmark{part2-figure_p_prime-p61R} &  \\\hline
% r = j+3
    $j+3$ &  &  & \tikzmark{part2-figure_p_prime-p62L}$q_0^2q_{j+3}\cdot h_{j+1}h_{j+3}h_{j+7}$\tikzmark{part2-figure_p_prime-p62R} \\\hline
% r = j+2
    $j+2$ &  &  &  \\\hline
% r = j+1
    $j+1$ &  & \tikzmark{part2-figure_p_prime-p63L}$(28)\; q_0^2q_{j+1}\cdot h_{j+4}h_{j+7}$\tikzmark{part2-figure_p_prime-p63R} &  \\\hline
% r = j
    $j$ &  &  & \tikzmark{part2-figure_p_prime-p64L}$q_0^2q_j\cdot h_jh_{j+4}h_{j+7}$\tikzmark{part2-figure_p_prime-p64R} \\\hline
    \end{tabular}
\begin{tikzpicture}[overlay, remember picture, shorten >=.5pt, shorten <=.5pt, transform canvas={yshift=.25\baselineskip}]
    \draw [->] ({pic cs:part2-figure_p_prime-p35R}) [red,bend left=16] to ({pic cs:part2-figure_p_prime-p36L});
    \draw [->] ({pic cs:part2-figure_p_prime-p35R}) [red,bend left=8] to ({pic cs:part2-figure_p_prime-p37L});
    \draw [->] ({pic cs:part2-figure_p_prime-p35R}) [red,bend left=15] to ({pic cs:part2-figure_p_prime-p38L});
    \draw [->] ({pic cs:part2-figure_p_prime-p39R}) [red,bend left=12] to ({pic cs:part2-figure_p_prime-p40L});
    \draw [->] ({pic cs:part2-figure_p_prime-p39R}) [red,bend right=4] to ({pic cs:part2-figure_p_prime-p41L});
    \draw [->] ({pic cs:part2-figure_p_prime-p42R}) [red,bend left=4] to ({pic cs:part2-figure_p_prime-p40L});
    \draw [->] ({pic cs:part2-figure_p_prime-p42R}) [red,out=-25,in=160] to ({pic cs:part2-figure_p_prime-p43L});
    \draw [->] ({pic cs:part2-figure_p_prime-p44R}) [red,bend left=4] to ({pic cs:part2-figure_p_prime-p41L});
    \draw [->] ({pic cs:part2-figure_p_prime-p44R}) [red,bend right=12] to ({pic cs:part2-figure_p_prime-p43L});
    \draw [->] ({pic cs:part2-figure_p_prime-p45R}) [red] to ({pic cs:part2-figure_p_prime-p46L});
    \draw [->] ({pic cs:part2-figure_p_prime-p47R}) [red,bend left=16] to ({pic cs:part2-figure_p_prime-p48L});
    \draw [->] ({pic cs:part2-figure_p_prime-p47R}) [red] to ({pic cs:part2-figure_p_prime-p49L});
    \draw [->] ({pic cs:part2-figure_p_prime-p47R}) [red,bend right=16] to ({pic cs:part2-figure_p_prime-p50L});
    \draw [->] ({pic cs:part2-figure_p_prime-p51R}) [red,bend left=12] to ({pic cs:part2-figure_p_prime-p52L});
    \draw [->] ({pic cs:part2-figure_p_prime-p51R}) [red,bend right=4] to ({pic cs:part2-figure_p_prime-p53L});
    \draw [->] ({pic cs:part2-figure_p_prime-p54R}) [red,bend left=4] to ({pic cs:part2-figure_p_prime-p52L});
    \draw [->] ({pic cs:part2-figure_p_prime-p54R}) [red,out=-25,in=160] to ({pic cs:part2-figure_p_prime-p55L});
    \draw [->] ({pic cs:part2-figure_p_prime-p56R}) [red,bend left=4] to ({pic cs:part2-figure_p_prime-p53L});
    \draw [->] ({pic cs:part2-figure_p_prime-p56R}) [red,bend right=12] to ({pic cs:part2-figure_p_prime-p55L});
    \draw [->] ({pic cs:part2-figure_p_prime-p57R}) [red] to ({pic cs:part2-figure_p_prime-p58L});
    \draw [->] ({pic cs:part2-figure_p_prime-p59R}) [red] to ({pic cs:part2-figure_p_prime-p60L});
    \draw [->] ({pic cs:part2-figure_p_prime-p61R}) [red] to ({pic cs:part2-figure_p_prime-p62L});
    \draw [->] ({pic cs:part2-figure_p_prime-p63R}) [red] to ({pic cs:part2-figure_p_prime-p64L});
    \draw [->] ({pic cs:part2-figure_p_prime-p67R}) [red,bend left=16] to ({pic cs:part2-figure_p_prime-p39L});
    \draw [->] ({pic cs:part2-figure_p_prime-p67R}) [red] to ({pic cs:part2-figure_p_prime-p42L});
    \draw [->] ({pic cs:part2-figure_p_prime-p67R}) [red,bend left=0] to ({pic cs:part2-figure_p_prime-p44L});
    \draw [->] ({pic cs:part2-figure_p_prime-p68R}) [red,bend left=0] to ({pic cs:part2-figure_p_prime-p51L});
    \draw [->] ({pic cs:part2-figure_p_prime-p68R}) [red] to ({pic cs:part2-figure_p_prime-p54L});
    \draw [->] ({pic cs:part2-figure_p_prime-p68R}) [red,bend left=0] to ({pic cs:part2-figure_p_prime-p56L});
\end{tikzpicture}
\par\endgroup
\caption{The BXSS in degree $(s,k,t) = (2,3,2^{j+2}\cdot 73+1)$ for \cref{CESS p prime}, where $j\geq 3$.
The vertical axis is the BX filtration and the horizontal axis is the $s$-degree. Red arrows indicate $d_1$-differentials, and blue arrows indicate $d_2$-differentials.}
\label{fig:figure_p_prime}
\label{fig:figure_p_prime-part2}
\end{figure}

    If $j=2$, the following differentials change:
    \[(5) d_1(q_{3}q_{6}q_{8}\cdot h_1h_{8}) =  q_{3}q_{5}q_{8}\cdot h_1h_{5}h_{8}\]
    \[(6) d_1(q_3q_5q_9\cdot h_1h_5)=q_3q_5q_8\cdot h_1h_5h_8\]
    \[(8) d_2(q_3q_5^2\cdot h_1h_9)=q_1q_5^2\cdot h_2^2h_9+q_3q_4^2\cdot h_1h_5h_9\]
    \[(12) d_1(q_1q_{j+4}q_{j+7}\cdot h_{j}^2) = q_1q_{j+3}q_{j+7}\cdot h_{j}^2h_{j+3}+q_1q_{j+4}q_{j+6}\cdot h_{j}^2h_{j+6}\]
    To account for the new linear relation $d_1(q_{3}q_{6}q_{8}\cdot h_1h_{8}+q_3q_5q_9\cdot h_1h_5)=0$, we have $d_1(q_3q_6q_9\cdot h_1)=q_{3}q_{6}q_{8}\cdot h_1h_{8}+q_3q_5q_9\cdot h_1h_5$.

\begin{figure}[H]
\centering
\begingroup
\setlength{\tabcolsep}{6pt}
\renewcommand{\arraystretch}{1.1}
\begin{tabular}{|wc{1.2cm}||wc{2.5cm}|p{5cm}|wc{4cm}|} \hline
      $w$ & $s=1$ & \multicolumn{1}{p{5cm}|}{\centering $s=2$} & $s=3$ \\\hline\hline
% r = 22
    $22$ &  & \tikzmark{part1-figure_p_prime_j2-p1L}$(1)\; q_6q_8^2\cdot h_1h_3$\tikzmark{part1-figure_p_prime_j2-p1R} &  \\\hline
% r = 21
    $21$ &  &  & \tikzmark{part1-figure_p_prime_j2-p2L}$q_5q_8^2\cdot h_1h_3h_5$\tikzmark{part1-figure_p_prime_j2-p2R} \\\hline
% r = 20
    $20$ &  &  &  \\\hline
% r = 19
    $19$ &  & \tikzmark{part1-figure_p_prime_j2-p3L}$(2)\; q_3q_8^2\cdot h_1h_6$\tikzmark{part1-figure_p_prime_j2-p3R} &  \\
     &  & \tikzmark{part1-figure_p_prime_j2-p4L}$(3)\; q_5^2q_9\cdot h_1h_3$\tikzmark{part1-figure_p_prime_j2-p4R} &  \\\hline
% r = 18
    $18$ & \tikzmark{part1-figure_p_prime_j2-p61L}$q_3q_6q_9\cdot h_1$\tikzmark{part1-figure_p_prime_j2-p61R} & \tikzmark{part1-figure_p_prime_j2-p6L}$(4)\; q_3q_6q_9\cdot h_0^2$\tikzmark{part1-figure_p_prime_j2-p6R} & \tikzmark{part1-figure_p_prime_j2-p5L}$q_5^2q_8\cdot h_1h_3h_8$\tikzmark{part1-figure_p_prime_j2-p5R} \\\hline
% r = 17
    $17$ &  & \tikzmark{part1-figure_p_prime_j2-p10L}$(5)\; q_3q_6q_8\cdot h_1h_8$\tikzmark{part1-figure_p_prime_j2-p10R} & \tikzmark{part1-figure_p_prime_j2-p7L}$q_2q_6q_9\cdot h_0^2h_2$\tikzmark{part1-figure_p_prime_j2-p7R} \\
     &  & \tikzmark{part1-figure_p_prime_j2-p12L}$(6)\; q_3q_5q_9\cdot h_1h_5$\tikzmark{part1-figure_p_prime_j2-p12R} & \tikzmark{part1-figure_p_prime_j2-p8L}$q_3q_5q_9\cdot h_0^2h_5$\tikzmark{part1-figure_p_prime_j2-p8R} \\
     &  & \tikzmark{part1-figure_p_prime_j2-p20L}$(11)\; q_1q_8^2\cdot h_3h_6$\tikzmark{part1-figure_p_prime_j2-p20R} & \tikzmark{part1-figure_p_prime_j2-p9L}$q_3q_6q_8\cdot h_0^2h_8$\tikzmark{part1-figure_p_prime_j2-p9R} \\\hline
% r = 16
    $16$ & \tikzmark{part1-figure_p_prime_j2-p62L}$q_1q_6q_9\cdot h_3$\tikzmark{part1-figure_p_prime_j2-p62R} & \tikzmark{part1-figure_p_prime_j2-p22L}$(12)\; q_1q_6q_9\cdot h_2^2$\tikzmark{part1-figure_p_prime_j2-p22R} & \tikzmark{part1-figure_p_prime_j2-p11L}$q_3q_5q_8\cdot h_1h_5h_8$\tikzmark{part1-figure_p_prime_j2-p11R} \\
     &  &  & \tikzmark{part1-figure_p_prime_j2-p21L}$q_0q_8^2\cdot h_0h_3h_6$\tikzmark{part1-figure_p_prime_j2-p21R} \\\hline
% r = 15
    $15$ &  & \tikzmark{part1-figure_p_prime_j2-p25L}$(13)\; q_1q_6q_8\cdot h_3h_8$\tikzmark{part1-figure_p_prime_j2-p25R} & \tikzmark{part1-figure_p_prime_j2-p23L}$q_1q_5q_9\cdot h_2^2h_5$\tikzmark{part1-figure_p_prime_j2-p23R} \\
     &  & \tikzmark{part1-figure_p_prime_j2-p28L}$(14)\; q_1q_5q_9\cdot h_3h_5$\tikzmark{part1-figure_p_prime_j2-p28R} & \tikzmark{part1-figure_p_prime_j2-p24L}$q_1q_6q_8\cdot h_2^2h_8$\tikzmark{part1-figure_p_prime_j2-p24R} \\
     &  & \tikzmark{part1-figure_p_prime_j2-p30L}$(15)\; q_0q_6q_9\cdot h_0h_3$\tikzmark{part1-figure_p_prime_j2-p30R} &  \\\hline
% r = 14
    $14$ &  &  & \tikzmark{part1-figure_p_prime_j2-p26L}$q_1q_5q_8\cdot h_3h_5h_8$\tikzmark{part1-figure_p_prime_j2-p26R} \\
     &  &  & \tikzmark{part1-figure_p_prime_j2-p27L}$q_0q_6q_8\cdot h_0h_3h_8$\tikzmark{part1-figure_p_prime_j2-p27R} \\
     &  &  & \tikzmark{part1-figure_p_prime_j2-p29L}$q_0q_5q_9\cdot h_0h_3h_5$\tikzmark{part1-figure_p_prime_j2-p29R} \\\hline
\end{tabular}
\begin{tikzpicture}[overlay, remember picture, shorten >=.5pt, shorten <=.5pt, transform canvas={yshift=.25\baselineskip}]

    \draw [->] ({pic cs:part1-figure_p_prime_j2-p1R}) [red] to ({pic cs:part1-figure_p_prime_j2-p2L});
    \draw [->] ({pic cs:part1-figure_p_prime_j2-p4R}) [red] to ({pic cs:part1-figure_p_prime_j2-p5L});
    \draw [->] ({pic cs:part1-figure_p_prime_j2-p6R}) [red,bend left=16] to ({pic cs:part1-figure_p_prime_j2-p7L});
    \draw [->] ({pic cs:part1-figure_p_prime_j2-p6R}) [red] to ({pic cs:part1-figure_p_prime_j2-p8L});
    \draw [->] ({pic cs:part1-figure_p_prime_j2-p6R}) [red,bend right=16] to ({pic cs:part1-figure_p_prime_j2-p9L});
    \draw [->] ({pic cs:part1-figure_p_prime_j2-p10R}) [red,bend left=4] to ({pic cs:part1-figure_p_prime_j2-p11L});
    \draw [->] ({pic cs:part1-figure_p_prime_j2-p12R}) [red,bend right=4] to ({pic cs:part1-figure_p_prime_j2-p11L});
    \draw [->] ({pic cs:part1-figure_p_prime_j2-p20R}) [red] to ({pic cs:part1-figure_p_prime_j2-p21L});
    \draw [->] ({pic cs:part1-figure_p_prime_j2-p22R}) [red,bend left=8] to ({pic cs:part1-figure_p_prime_j2-p23L});
    \draw [->] ({pic cs:part1-figure_p_prime_j2-p22R}) [red,bend right=8] to ({pic cs:part1-figure_p_prime_j2-p24L});
    \draw [->] ({pic cs:part1-figure_p_prime_j2-p25R}) [red,bend left=12] to ({pic cs:part1-figure_p_prime_j2-p26L});
    \draw [->] ({pic cs:part1-figure_p_prime_j2-p25R}) [red,bend right=4] to ({pic cs:part1-figure_p_prime_j2-p27L});
    \draw [->] ({pic cs:part1-figure_p_prime_j2-p28R}) [red,bend left=4] to ({pic cs:part1-figure_p_prime_j2-p26L});
    \draw [->] ({pic cs:part1-figure_p_prime_j2-p28R}) [red,bend right=4] to ({pic cs:part1-figure_p_prime_j2-p29L});
    \draw [->] ({pic cs:part1-figure_p_prime_j2-p30R}) [red,bend left=4] to ({pic cs:part1-figure_p_prime_j2-p27L});
    \draw [->] ({pic cs:part1-figure_p_prime_j2-p30R}) [red,bend right=12] to ({pic cs:part1-figure_p_prime_j2-p29L});
    \draw [->] ({pic cs:part1-figure_p_prime_j2-p61R}) [red,bend left=8] to ({pic cs:part1-figure_p_prime_j2-p10L});
    \draw [->] ({pic cs:part1-figure_p_prime_j2-p61R}) [red,bend right=8] to ({pic cs:part1-figure_p_prime_j2-p12L});
    \draw [->] ({pic cs:part1-figure_p_prime_j2-p62R}) [red,bend left=16] to ({pic cs:part1-figure_p_prime_j2-p25L});
    \draw [->] ({pic cs:part1-figure_p_prime_j2-p62R}) [red] to ({pic cs:part1-figure_p_prime_j2-p28L});
    \draw [->] ({pic cs:part1-figure_p_prime_j2-p62R}) [red,bend right=16] to ({pic cs:part1-figure_p_prime_j2-p30L});
\end{tikzpicture}
\par\endgroup
\end{figure}

\begin{figure}[H]
\centering
\begingroup
\setlength{\tabcolsep}{6pt}
\renewcommand{\arraystretch}{1.1}
\begin{tabular}{|wc{1.2cm}||wc{2.5cm}|p{5cm}|wc{4cm}|} \hline
      $w$ & $s=1$ & \multicolumn{1}{p{5cm}|}{\centering $s=2$} & $s=3$ \\\hline\hline
% r = 13
    $13$ & \tikzmark{part2-figure_p_prime_j2-p63L}$q_1q_3q_9\cdot h_6$\tikzmark{part2-figure_p_prime_j2-p63R} & \tikzmark{part2-figure_p_prime_j2-p13L}$(8)\; q_3q_5^2\cdot h_1h_9$\tikzmark{part2-figure_p_prime_j2-p13R} &  \\
     &  & \tikzmark{part2-figure_p_prime_j2-p16L}$(9)\; q_2^2q_9\cdot h_1h_6$\tikzmark{part2-figure_p_prime_j2-p16R} &  \\
     &  & \tikzmark{part2-figure_p_prime_j2-p31L}$(16)\; q_1q_3q_9\cdot h_5^2$\tikzmark{part2-figure_p_prime_j2-p31R} &  \\\hline
% r = 12
    $12$ &  & \tikzmark{part2-figure_p_prime_j2-p35L}$(17)\; q_1q_3q_8\cdot h_6h_8$\tikzmark{part2-figure_p_prime_j2-p35R} & \tikzmark{part2-figure_p_prime_j2-p17L}$q_2^2q_8\cdot h_1h_6h_8$\tikzmark{part2-figure_p_prime_j2-p17R} \\
     &  & \tikzmark{part2-figure_p_prime_j2-p38L}$(18)\; q_1q_2q_9\cdot h_2h_6$\tikzmark{part2-figure_p_prime_j2-p38R} & \tikzmark{part2-figure_p_prime_j2-p32L}$q_1q_2q_9\cdot h_2h_5^2$\tikzmark{part2-figure_p_prime_j2-p32R} \\
     &  & \tikzmark{part2-figure_p_prime_j2-p40L}$(19)\; q_0q_3q_9\cdot h_0h_6$\tikzmark{part2-figure_p_prime_j2-p40R} & \tikzmark{part2-figure_p_prime_j2-p33L}$q_1q_3q_8\cdot h_5^2h_8$\tikzmark{part2-figure_p_prime_j2-p33R} \\
     &  &  & \tikzmark{part2-figure_p_prime_j2-p34L}$q_0q_3q_9\cdot h_0h_5^2$\tikzmark{part2-figure_p_prime_j2-p34R} \\\hline
% r = 11
    $11$ &  & \tikzmark{part2-figure_p_prime_j2-p41L}$(20)\; q_1q_5^2\cdot h_3h_9$\tikzmark{part2-figure_p_prime_j2-p41R} & \tikzmark{part2-figure_p_prime_j2-p14L}$q_1q_5^2\cdot h_2^2h_9$\tikzmark{part2-figure_p_prime_j2-p14R} \\
     &  & \tikzmark{part2-figure_p_prime_j2-p43L}$(21)\; q_1q_3q_7\cdot h_8^2$\tikzmark{part2-figure_p_prime_j2-p43R} & \tikzmark{part2-figure_p_prime_j2-p15L}$q_3q_4^2\cdot h_1h_5h_9$\tikzmark{part2-figure_p_prime_j2-p15R} \\
     &  &  & \tikzmark{part2-figure_p_prime_j2-p36L}$q_1q_2q_8\cdot h_2h_6h_8$\tikzmark{part2-figure_p_prime_j2-p36R} \\
     &  &  & \tikzmark{part2-figure_p_prime_j2-p37L}$q_0q_3q_8\cdot h_0h_6h_8$\tikzmark{part2-figure_p_prime_j2-p37R} \\
     &  &  & \tikzmark{part2-figure_p_prime_j2-p39L}$q_0q_2q_9\cdot h_0h_2h_6$\tikzmark{part2-figure_p_prime_j2-p39R} \\\hline
% r = 10
    $10$ & \tikzmark{part2-figure_p_prime_j2-p64L}$q_1q_3q_6\cdot h_9$\tikzmark{part2-figure_p_prime_j2-p64R} & \tikzmark{part2-figure_p_prime_j2-p18L}$(10)\; q_2^2q_6\cdot h_1h_9$\tikzmark{part2-figure_p_prime_j2-p18R} & \tikzmark{part2-figure_p_prime_j2-p42L}$q_0q_5^2\cdot h_0h_3h_9$\tikzmark{part2-figure_p_prime_j2-p42R} \\
     &  &  & \tikzmark{part2-figure_p_prime_j2-p44L}$q_1q_2q_7\cdot h_2h_8^2$\tikzmark{part2-figure_p_prime_j2-p44R} \\
     &  &  & \tikzmark{part2-figure_p_prime_j2-p45L}$q_1q_3q_6\cdot h_6h_8^2$\tikzmark{part2-figure_p_prime_j2-p45R} \\
     &  &  & \tikzmark{part2-figure_p_prime_j2-p46L}$q_0q_3q_7\cdot h_0h_8^2$\tikzmark{part2-figure_p_prime_j2-p46R} \\\hline
% r = 9
    $9$ &  & \tikzmark{part2-figure_p_prime_j2-p47L}$(22)\; q_1q_3q_5\cdot h_5h_9$\tikzmark{part2-figure_p_prime_j2-p47R} & \tikzmark{part2-figure_p_prime_j2-p19L}$q_2^2q_5\cdot h_1h_5h_9$\tikzmark{part2-figure_p_prime_j2-p19R} \\
     &  & \tikzmark{part2-figure_p_prime_j2-p50L}$(23)\; q_1q_2q_6\cdot h_2h_9$\tikzmark{part2-figure_p_prime_j2-p50R} &  \\
     &  & \tikzmark{part2-figure_p_prime_j2-p52L}$(24)\; q_0q_3q_6\cdot h_0h_9$\tikzmark{part2-figure_p_prime_j2-p52R} &  \\
     &  & \tikzmark{part2-figure_p_prime_j2-p55L}$(26)\; q_0^2q_9\cdot h_3h_6$\tikzmark{part2-figure_p_prime_j2-p55R} &  \\\hline
% r = 8
    $8$ &  &  & \tikzmark{part2-figure_p_prime_j2-p48L}$q_1q_2q_5\cdot h_2h_5h_9$\tikzmark{part2-figure_p_prime_j2-p48R} \\
     &  &  & \tikzmark{part2-figure_p_prime_j2-p49L}$q_0q_3q_5\cdot h_0h_5h_9$\tikzmark{part2-figure_p_prime_j2-p49R} \\
     &  &  & \tikzmark{part2-figure_p_prime_j2-p51L}$q_0q_2q_6\cdot h_0h_2h_9$\tikzmark{part2-figure_p_prime_j2-p51R} \\
     &  &  & \tikzmark{part2-figure_p_prime_j2-p56L}$q_0^2q_8\cdot h_3h_6h_8$\tikzmark{part2-figure_p_prime_j2-p56R} \\\hline
% r = 7
    $7$ &  &  &  \\\hline
% r = 6
    $6$ &  & \tikzmark{part2-figure_p_prime_j2-p57L}$(27)\; q_0^2q_6\cdot h_3h_9$\tikzmark{part2-figure_p_prime_j2-p57R} &  \\\hline
% r = 5
    $5$ &  & \tikzmark{part2-figure_p_prime_j2-p53L}$(25)\; q_1q_2^2\cdot h_6h_9$\tikzmark{part2-figure_p_prime_j2-p53R} & \tikzmark{part2-figure_p_prime_j2-p58L}$q_0^2q_5\cdot h_3h_5h_9$\tikzmark{part2-figure_p_prime_j2-p58R} \\\hline
% r = 4
    $4$ &  &  & \tikzmark{part2-figure_p_prime_j2-p54L}$q_0q_2^2\cdot h_0h_6h_9$\tikzmark{part2-figure_p_prime_j2-p54R} \\\hline
% r = 3
    $3$ &  & \tikzmark{part2-figure_p_prime_j2-p59L}$(28)\; q_0^2q_3\cdot h_6h_9$\tikzmark{part2-figure_p_prime_j2-p59R} &  \\\hline
% r = 2
    $2$ &  &  & \tikzmark{part2-figure_p_prime_j2-p60L}$q_0^2q_2\cdot h_2h_6h_9$\tikzmark{part2-figure_p_prime_j2-p60R} \\\hline
    \end{tabular}
\begin{tikzpicture}[overlay, remember picture, shorten >=.5pt, shorten <=.5pt, transform canvas={yshift=.25\baselineskip}]
    \draw [->] ({pic cs:part2-figure_p_prime_j2-p13R}) [blue,bend left=8] to ({pic cs:part2-figure_p_prime_j2-p14L});
    \draw [->] ({pic cs:part2-figure_p_prime_j2-p13R}) [blue,bend right=8] to ({pic cs:part2-figure_p_prime_j2-p15L});
    \draw [->] ({pic cs:part2-figure_p_prime_j2-p16R}) [red] to ({pic cs:part2-figure_p_prime_j2-p17L});
    \draw [->] ({pic cs:part2-figure_p_prime_j2-p18R}) [red] to ({pic cs:part2-figure_p_prime_j2-p19L});
    \draw [->] ({pic cs:part2-figure_p_prime_j2-p31R}) [red,bend left=16] to ({pic cs:part2-figure_p_prime_j2-p32L});
    \draw [->] ({pic cs:part2-figure_p_prime_j2-p31R}) [red] to ({pic cs:part2-figure_p_prime_j2-p33L});
    \draw [->] ({pic cs:part2-figure_p_prime_j2-p31R}) [red,bend left=0] to ({pic cs:part2-figure_p_prime_j2-p34L});
    \draw [->] ({pic cs:part2-figure_p_prime_j2-p35R}) [red,bend left=12] to ({pic cs:part2-figure_p_prime_j2-p36L});
    \draw [->] ({pic cs:part2-figure_p_prime_j2-p35R}) [red,bend right=4] to ({pic cs:part2-figure_p_prime_j2-p37L});
    \draw [->] ({pic cs:part2-figure_p_prime_j2-p38R}) [red,bend left=4] to ({pic cs:part2-figure_p_prime_j2-p36L});
    \draw [->] ({pic cs:part2-figure_p_prime_j2-p38R}) [red,bend right=4] to ({pic cs:part2-figure_p_prime_j2-p39L});
    \draw [->] ({pic cs:part2-figure_p_prime_j2-p40R}) [red,bend left=4] to ({pic cs:part2-figure_p_prime_j2-p37L});
    \draw [->] ({pic cs:part2-figure_p_prime_j2-p40R}) [red,bend right=12] to ({pic cs:part2-figure_p_prime_j2-p39L});
    \draw [->] ({pic cs:part2-figure_p_prime_j2-p41R}) [red] to ({pic cs:part2-figure_p_prime_j2-p42L});
    \draw [->] ({pic cs:part2-figure_p_prime_j2-p43R}) [red,bend left=16] to ({pic cs:part2-figure_p_prime_j2-p44L});
    \draw [->] ({pic cs:part2-figure_p_prime_j2-p43R}) [red] to ({pic cs:part2-figure_p_prime_j2-p45L});
    \draw [->] ({pic cs:part2-figure_p_prime_j2-p43R}) [red,bend right=16] to ({pic cs:part2-figure_p_prime_j2-p46L});
    \draw [->] ({pic cs:part2-figure_p_prime_j2-p47R}) [red,bend left=12] to ({pic cs:part2-figure_p_prime_j2-p48L});
    \draw [->] ({pic cs:part2-figure_p_prime_j2-p47R}) [red,bend right=4] to ({pic cs:part2-figure_p_prime_j2-p49L});
    \draw [->] ({pic cs:part2-figure_p_prime_j2-p50R}) [red,bend left=4] to ({pic cs:part2-figure_p_prime_j2-p48L});
    \draw [->] ({pic cs:part2-figure_p_prime_j2-p50R}) [red,bend right=4] to ({pic cs:part2-figure_p_prime_j2-p51L});
    \draw [->] ({pic cs:part2-figure_p_prime_j2-p52R}) [red,bend left=4] to ({pic cs:part2-figure_p_prime_j2-p49L});
    \draw [->] ({pic cs:part2-figure_p_prime_j2-p52R}) [red,bend right=12] to ({pic cs:part2-figure_p_prime_j2-p51L});
    \draw [->] ({pic cs:part2-figure_p_prime_j2-p53R}) [red] to ({pic cs:part2-figure_p_prime_j2-p54L});
    \draw [->] ({pic cs:part2-figure_p_prime_j2-p55R}) [red] to ({pic cs:part2-figure_p_prime_j2-p56L});
    \draw [->] ({pic cs:part2-figure_p_prime_j2-p57R}) [red] to ({pic cs:part2-figure_p_prime_j2-p58L});
    \draw [->] ({pic cs:part2-figure_p_prime_j2-p59R}) [red] to ({pic cs:part2-figure_p_prime_j2-p60L});
    \draw [->] ({pic cs:part2-figure_p_prime_j2-p63R}) [red,bend left=16] to ({pic cs:part2-figure_p_prime_j2-p35L});
    \draw [->] ({pic cs:part2-figure_p_prime_j2-p63R}) [red] to ({pic cs:part2-figure_p_prime_j2-p38L});
    \draw [->] ({pic cs:part2-figure_p_prime_j2-p63R}) [red,bend right=16] to ({pic cs:part2-figure_p_prime_j2-p40L});
    \draw [->] ({pic cs:part2-figure_p_prime_j2-p64R}) [red,bend left=16] to ({pic cs:part2-figure_p_prime_j2-p47L});
    \draw [->] ({pic cs:part2-figure_p_prime_j2-p64R}) [red] to ({pic cs:part2-figure_p_prime_j2-p50L});
    \draw [->] ({pic cs:part2-figure_p_prime_j2-p64R}) [red,bend right=16] to ({pic cs:part2-figure_p_prime_j2-p52L});
\end{tikzpicture}
\par\endgroup
\caption{The BXSS in degree $(s,k,t) = (2,3,2^{4}\cdot 73+1)$ for \cref{CESS p prime}.
The vertical axis is the BX filtration and the horizontal axis is the $s$-degree. Red arrows indicate $d_1$-differentials, and blue arrows indicate $d_2$-differentials.}
\label{fig:figure_p_prime_j2}
\label{fig:figure_p_prime_j2-part2}
\end{figure}
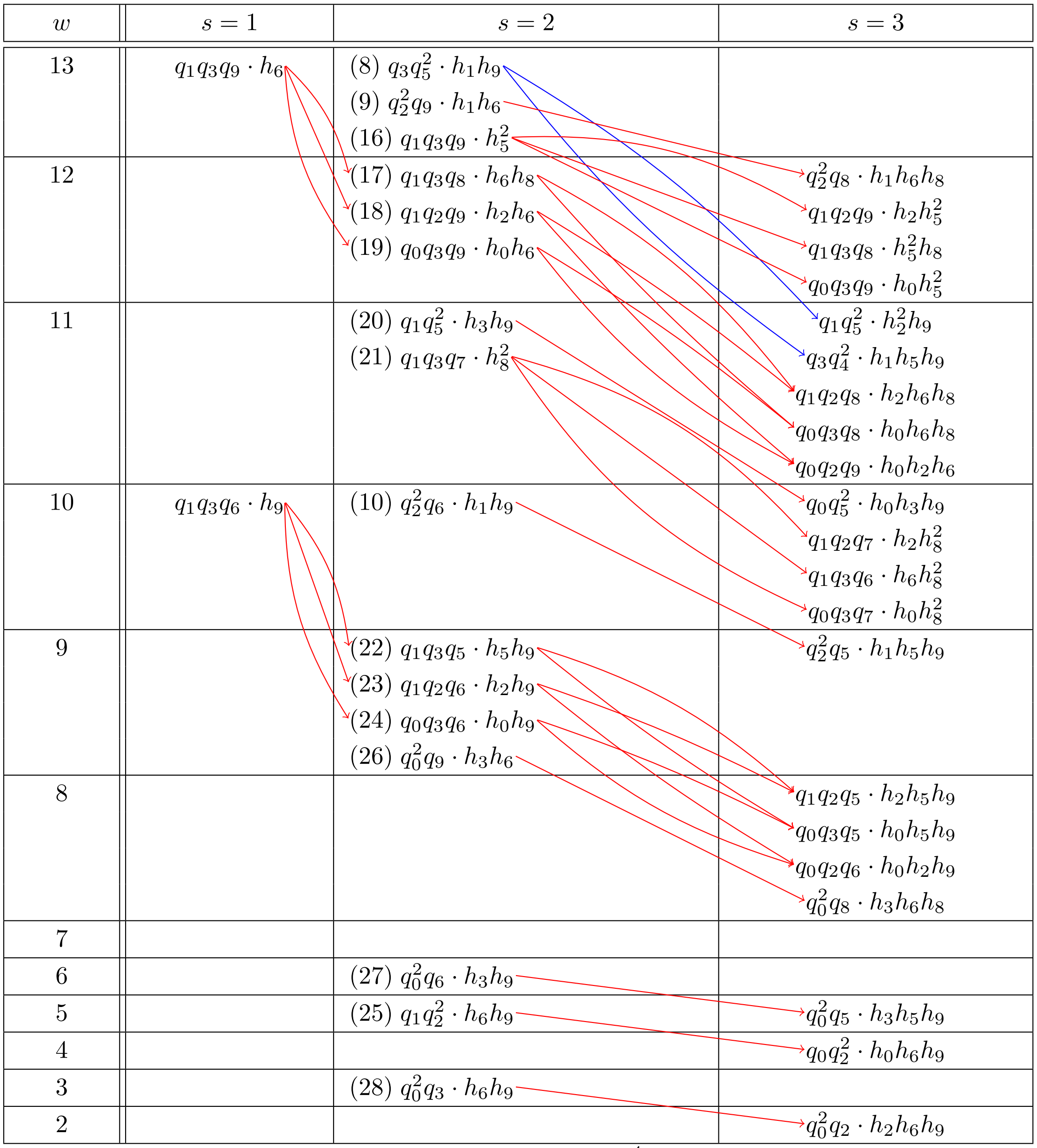

    If $j=1$, the following classes and differentials change:
    \[(3) q_4^2q_8\cdot h_1h_2=0\]
    \[(12)d_1(q_1q_5q_7\cdot h_1^2)=q_1q_4q_7\cdot h_1^2h_4+q_1q_5q_6\cdot h_1^2h_6\]
    \[\text{(18) is the same as (9)}\]

    \[(19) d_1(q_0q_2q_8\cdot h_0h_5)=q_0q_2q_7\cdot h_0h_5h_7\]
    \[\text{(23) is the same as (10)}\]

    \[(24)d_1(q_0q_2q_5\cdot h_0h_8)=q_0q_2q_4\cdot h_0h_4h_8\]

\usetikzlibrary{tikzmark}

\begin{figure}[H]
\centering
\begingroup
\setlength{\tabcolsep}{6pt}
\renewcommand{\arraystretch}{1.1}
\begin{tabular}{|wc{1.2cm}||wc{2.5cm}|p{5cm}|wc{4cm}|} \hline
      $w$ & $s=1$ & \multicolumn{1}{p{5cm}|}{\centering $s=2$} & $s=3$ \\\hline\hline
% r = 16
    $16$ &  & \tikzmark{part1-figure_p_prime_j1-p1L}$(2)\; q_2q_7^2\cdot h_1h_5$\tikzmark{part1-figure_p_prime_j1-p1R} &  \\\hline
% r = 15
    $15$ & \tikzmark{part1-figure_p_prime_j1-p53L}$q_2q_5q_8\cdot h_1$\tikzmark{part1-figure_p_prime_j1-p53R} & \tikzmark{part1-figure_p_prime_j1-p3L}$(4)\; q_2q_5q_8\cdot h_0^2$\tikzmark{part1-figure_p_prime_j1-p3R} & \tikzmark{part1-figure_p_prime_j1-p2L}$q_1q_7^2\cdot h_1^2h_5$\tikzmark{part1-figure_p_prime_j1-p2R} \\
     &  & \tikzmark{part1-figure_p_prime_j1-p18L}$(11)\; q_1q_7^2\cdot h_2h_5$\tikzmark{part1-figure_p_prime_j1-p18R} &  \\\hline
% r = 14
    $14$ & \tikzmark{part1-figure_p_prime_j1-p54L}$q_1q_5q_8\cdot h_2$\tikzmark{part1-figure_p_prime_j1-p54R} & \tikzmark{part1-figure_p_prime_j1-p6L}$(5)\; q_2q_5q_7\cdot h_1h_7$\tikzmark{part1-figure_p_prime_j1-p6R} & \tikzmark{part1-figure_p_prime_j1-p4L}$q_2q_4q_8\cdot h_0^2h_4$\tikzmark{part1-figure_p_prime_j1-p4R} \\
     &  & \tikzmark{part1-figure_p_prime_j1-p9L}$(6)\; q_2q_4q_8\cdot h_1h_4$\tikzmark{part1-figure_p_prime_j1-p9R} & \tikzmark{part1-figure_p_prime_j1-p5L}$q_2q_5q_7\cdot h_0^2h_7$\tikzmark{part1-figure_p_prime_j1-p5R} \\
     &  & \tikzmark{part1-figure_p_prime_j1-p11L}$(7)\; q_1q_5q_8\cdot h_1^2$\tikzmark{part1-figure_p_prime_j1-p11R} & \tikzmark{part1-figure_p_prime_j1-p19L}$q_0q_7^2\cdot h_0h_2h_5$\tikzmark{part1-figure_p_prime_j1-p19R} \\\hline
% r = 13
    $13$ &  & \tikzmark{part1-figure_p_prime_j1-p20L}$(12)\; q_1q_5q_7\cdot h_1^2$\tikzmark{part1-figure_p_prime_j1-p20R} & \tikzmark{part1-figure_p_prime_j1-p7L}$q_1q_5q_7\cdot h_1^2h_7$\tikzmark{part1-figure_p_prime_j1-p7R} \\
     &  & \tikzmark{part1-figure_p_prime_j1-p23L}$(13)\; q_1q_5q_7\cdot h_2h_7$\tikzmark{part1-figure_p_prime_j1-p23R} & \tikzmark{part1-figure_p_prime_j1-p8L}$q_2q_4q_7\cdot h_1h_4h_7$\tikzmark{part1-figure_p_prime_j1-p8R} \\
     &  & \tikzmark{part1-figure_p_prime_j1-p26L}$(14)\; q_1q_4q_8\cdot h_2h_4$\tikzmark{part1-figure_p_prime_j1-p26R} & \tikzmark{part1-figure_p_prime_j1-p10L}$q_1q_4q_8\cdot h_1^2h_4$\tikzmark{part1-figure_p_prime_j1-p10R} \\
     &  & \tikzmark{part1-figure_p_prime_j1-p28L}$(15)\; q_0q_5q_8\cdot h_0h_2$\tikzmark{part1-figure_p_prime_j1-p28R} &  \\\hline
% r = 12
    $12$ &  &  & \tikzmark{part1-figure_p_prime_j1-p21L}$q_1q_4q_7\cdot h_1^2h_4$\tikzmark{part1-figure_p_prime_j1-p21R} \\
     &  &  & \tikzmark{part1-figure_p_prime_j1-p22L}$q_1q_5q_6\cdot h_1^2h_6$\tikzmark{part1-figure_p_prime_j1-p22R} \\
     &  &  & \tikzmark{part1-figure_p_prime_j1-p24L}$q_1q_4q_7\cdot h_2h_4h_7$\tikzmark{part1-figure_p_prime_j1-p24R} \\
     &  &  & \tikzmark{part1-figure_p_prime_j1-p25L}$q_0q_5q_7\cdot h_0h_2h_7$\tikzmark{part1-figure_p_prime_j1-p25R} \\
     &  &  & \tikzmark{part1-figure_p_prime_j1-p27L}$q_0q_4q_8\cdot h_0h_2h_4$\tikzmark{part1-figure_p_prime_j1-p27R} \\\hline
\end{tabular}
\begin{tikzpicture}[overlay, remember picture, shorten >=.5pt, shorten <=.5pt, transform canvas={yshift=.25\baselineskip}]

    \draw [->] ({pic cs:part1-figure_p_prime_j1-p1R}) [red] to ({pic cs:part1-figure_p_prime_j1-p2L});
    \draw [->] ({pic cs:part1-figure_p_prime_j1-p3R}) [red,bend left=8] to ({pic cs:part1-figure_p_prime_j1-p4L});
    \draw [->] ({pic cs:part1-figure_p_prime_j1-p3R}) [red,bend right=8] to ({pic cs:part1-figure_p_prime_j1-p5L});
    \draw [->] ({pic cs:part1-figure_p_prime_j1-p6R}) [red,bend left=12] to ({pic cs:part1-figure_p_prime_j1-p7L});
    \draw [->] ({pic cs:part1-figure_p_prime_j1-p6R}) [red,bend right=4] to ({pic cs:part1-figure_p_prime_j1-p8L});
    \draw [->] ({pic cs:part1-figure_p_prime_j1-p9R}) [red,bend left=4] to ({pic cs:part1-figure_p_prime_j1-p8L});
    \draw [->] ({pic cs:part1-figure_p_prime_j1-p9R}) [red,bend right=4] to ({pic cs:part1-figure_p_prime_j1-p10L});
    \draw [->] ({pic cs:part1-figure_p_prime_j1-p11R}) [red,bend left=4] to ({pic cs:part1-figure_p_prime_j1-p7L});
    \draw [->] ({pic cs:part1-figure_p_prime_j1-p11R}) [red,bend right=12] to ({pic cs:part1-figure_p_prime_j1-p10L});
    \draw [->] ({pic cs:part1-figure_p_prime_j1-p18R}) [red] to ({pic cs:part1-figure_p_prime_j1-p19L});
    \draw [->] ({pic cs:part1-figure_p_prime_j1-p20R}) [red,bend left=8] to ({pic cs:part1-figure_p_prime_j1-p21L});
    \draw [->] ({pic cs:part1-figure_p_prime_j1-p20R}) [red,bend left=0] to ({pic cs:part1-figure_p_prime_j1-p22L});
    \draw [->] ({pic cs:part1-figure_p_prime_j1-p23R}) [red,bend left=12] to ({pic cs:part1-figure_p_prime_j1-p24L});
    \draw [->] ({pic cs:part1-figure_p_prime_j1-p23R}) [red,bend right=4] to ({pic cs:part1-figure_p_prime_j1-p25L});
    \draw [->] ({pic cs:part1-figure_p_prime_j1-p26R}) [red,bend left=4] to ({pic cs:part1-figure_p_prime_j1-p24L});
    \draw [->] ({pic cs:part1-figure_p_prime_j1-p26R}) [red,bend right=4] to ({pic cs:part1-figure_p_prime_j1-p27L});
    \draw [->] ({pic cs:part1-figure_p_prime_j1-p28R}) [red,bend left=4] to ({pic cs:part1-figure_p_prime_j1-p25L});
    \draw [->] ({pic cs:part1-figure_p_prime_j1-p28R}) [red,bend right=12] to ({pic cs:part1-figure_p_prime_j1-p27L});
    \draw [->] ({pic cs:part1-figure_p_prime_j1-p53R}) [red,bend left=16] to ({pic cs:part1-figure_p_prime_j1-p6L});
    \draw [->] ({pic cs:part1-figure_p_prime_j1-p53R}) [red] to ({pic cs:part1-figure_p_prime_j1-p9L});
    \draw [->] ({pic cs:part1-figure_p_prime_j1-p53R}) [red,bend right=16] to ({pic cs:part1-figure_p_prime_j1-p11L});
    \draw [->] ({pic cs:part1-figure_p_prime_j1-p54R}) [red,bend left=16] to ({pic cs:part1-figure_p_prime_j1-p23L});
    \draw [->] ({pic cs:part1-figure_p_prime_j1-p54R}) [red] to ({pic cs:part1-figure_p_prime_j1-p26L});
    \draw [->] ({pic cs:part1-figure_p_prime_j1-p54R}) [red,bend right=16] to ({pic cs:part1-figure_p_prime_j1-p28L});
\end{tikzpicture}
\par\endgroup
\end{figure}

\begin{figure}[H]
\centering
\begingroup
\setlength{\tabcolsep}{6pt}
\renewcommand{\arraystretch}{1.1}
\begin{tabular}{|wc{1.2cm}||wc{2.5cm}|p{5cm}|wc{4cm}|} \hline
      $w$ & $s=1$ & \multicolumn{1}{p{5cm}|}{\centering $s=2$} & $s=3$ \\\hline\hline
% r = 11
    $11$ & \tikzmark{part2-figure_p_prime_j1-p55L}$q_1q_2q_8\cdot h_5$\tikzmark{part2-figure_p_prime_j1-p55R} & \tikzmark{part2-figure_p_prime_j1-p29L}$(16)\; q_1q_2q_8\cdot h_4^2$\tikzmark{part2-figure_p_prime_j1-p29R} &  \\\hline
% r = 10
    $10$ &  & \tikzmark{part2-figure_p_prime_j1-p12L}$(8)\; q_2q_4^2\cdot h_1h_8$\tikzmark{part2-figure_p_prime_j1-p12R} & \tikzmark{part2-figure_p_prime_j1-p30L}$q_1^2q_8\cdot h_1h_4^2$\tikzmark{part2-figure_p_prime_j1-p30R} \\
     &  & \tikzmark{part2-figure_p_prime_j1-p14L}$(9,18)\; q_1^2q_8\cdot h_1h_5$\tikzmark{part2-figure_p_prime_j1-p14R} & \tikzmark{part2-figure_p_prime_j1-p31L}$q_1q_2q_7\cdot h_4^2h_7$\tikzmark{part2-figure_p_prime_j1-p31R} \\
     &  & \tikzmark{part2-figure_p_prime_j1-p33L}$(17)\; q_1q_2q_7\cdot h_5h_7$\tikzmark{part2-figure_p_prime_j1-p33R} & \tikzmark{part2-figure_p_prime_j1-p32L}$q_0q_2q_8\cdot h_0h_4^2$\tikzmark{part2-figure_p_prime_j1-p32R} \\
     &  & \tikzmark{part2-figure_p_prime_j1-p35L}$(19)\; q_0q_2q_8\cdot h_0h_5$\tikzmark{part2-figure_p_prime_j1-p35R} &  \\\hline
% r = 9
    $9$ &  & \tikzmark{part2-figure_p_prime_j1-p36L}$(20)\; q_1q_4^2\cdot h_2h_8$\tikzmark{part2-figure_p_prime_j1-p36R} & \tikzmark{part2-figure_p_prime_j1-p13L}$q_1q_4^2\cdot h_1^2h_8$\tikzmark{part2-figure_p_prime_j1-p13R} \\
     &  & \tikzmark{part2-figure_p_prime_j1-p38L}$(21)\; q_1q_2q_6\cdot h_7^2$\tikzmark{part2-figure_p_prime_j1-p38R} & \tikzmark{part2-figure_p_prime_j1-p15L}$q_1^2q_7\cdot h_1h_5h_7$\tikzmark{part2-figure_p_prime_j1-p15R} \\
     &  &  & \tikzmark{part2-figure_p_prime_j1-p34L}$q_0q_2q_7\cdot h_0h_5h_7$\tikzmark{part2-figure_p_prime_j1-p34R} \\\hline
% r = 8
    $8$ & \tikzmark{part2-figure_p_prime_j1-p56L}$q_1q_2q_5\cdot h_8$\tikzmark{part2-figure_p_prime_j1-p56R} & \tikzmark{part2-figure_p_prime_j1-p47L}$(26)\; q_0^2q_8\cdot h_2h_5$\tikzmark{part2-figure_p_prime_j1-p47R} & \tikzmark{part2-figure_p_prime_j1-p37L}$q_0q_4^2\cdot h_0h_2h_8$\tikzmark{part2-figure_p_prime_j1-p37R} \\
     &  &  & \tikzmark{part2-figure_p_prime_j1-p39L}$q_1^2q_6\cdot h_1h_7^2$\tikzmark{part2-figure_p_prime_j1-p39R} \\
     &  &  & \tikzmark{part2-figure_p_prime_j1-p40L}$q_1q_2q_5\cdot h_5h_7^2$\tikzmark{part2-figure_p_prime_j1-p40R} \\
     &  &  & \tikzmark{part2-figure_p_prime_j1-p41L}$q_0q_2q_6\cdot h_0h_7^2$\tikzmark{part2-figure_p_prime_j1-p41R} \\\hline
% r = 7
    $7$ &  & \tikzmark{part2-figure_p_prime_j1-p16L}$(10,23)\; q_1^2q_5\cdot h_1h_8$\tikzmark{part2-figure_p_prime_j1-p16R} & \tikzmark{part2-figure_p_prime_j1-p48L}$q_0^2q_7\cdot h_2h_5h_7$\tikzmark{part2-figure_p_prime_j1-p48R} \\
     &  & \tikzmark{part2-figure_p_prime_j1-p42L}$(22)\; q_1q_2q_4\cdot h_4h_8$\tikzmark{part2-figure_p_prime_j1-p42R} &  \\
     &  & \tikzmark{part2-figure_p_prime_j1-p44L}$(24)\; q_0q_2q_5\cdot h_0h_8$\tikzmark{part2-figure_p_prime_j1-p44R} &  \\\hline
% r = 6
    $6$ &  &  & \tikzmark{part2-figure_p_prime_j1-p17L}$q_1^2q_4\cdot h_1h_4h_8$\tikzmark{part2-figure_p_prime_j1-p17R} \\
     &  &  & \tikzmark{part2-figure_p_prime_j1-p43L}$q_0q_2q_4\cdot h_0h_4h_8$\tikzmark{part2-figure_p_prime_j1-p43R} \\\hline
% r = 5
    $5$ &  & \tikzmark{part2-figure_p_prime_j1-p49L}$(27)\; q_0^2q_5\cdot h_2h_8$\tikzmark{part2-figure_p_prime_j1-p49R} &  \\\hline
% r = 4
    $4$ &  &  & \tikzmark{part2-figure_p_prime_j1-p50L}$q_0^2q_4\cdot h_2h_4h_8$\tikzmark{part2-figure_p_prime_j1-p50R} \\\hline
% r = 3
    $3$ &  & \tikzmark{part2-figure_p_prime_j1-p45L}$(25)\; q_1^3\cdot h_5h_8$\tikzmark{part2-figure_p_prime_j1-p45R} &  \\\hline
% r = 2
    $2$ &  & \tikzmark{part2-figure_p_prime_j1-p51L}$(28)\; q_0^2q_2\cdot h_5h_8$\tikzmark{part2-figure_p_prime_j1-p51R} & \tikzmark{part2-figure_p_prime_j1-p46L}$q_0q_1^2\cdot h_0h_5h_8$\tikzmark{part2-figure_p_prime_j1-p46R} \\\hline
% r = 1
    $1$ &  &  & \tikzmark{part2-figure_p_prime_j1-p52L}$q_0^2q_1\cdot h_1h_5h_8$\tikzmark{part2-figure_p_prime_j1-p52R} \\\hline
    \end{tabular}
\begin{tikzpicture}[overlay, remember picture, shorten >=.5pt, shorten <=.5pt, transform canvas={yshift=.25\baselineskip}]
    \draw [->] ({pic cs:part2-figure_p_prime_j1-p12R}) [red] to ({pic cs:part2-figure_p_prime_j1-p13L});
    \draw [->] ({pic cs:part2-figure_p_prime_j1-p14R}) [red,bend left=4] to ({pic cs:part2-figure_p_prime_j1-p15L});
    \draw [->] ({pic cs:part2-figure_p_prime_j1-p16R}) [red,bend left=4] to ({pic cs:part2-figure_p_prime_j1-p17L});
    \draw [->] ({pic cs:part2-figure_p_prime_j1-p29R}) [red,bend left=16] to ({pic cs:part2-figure_p_prime_j1-p30L});
    \draw [->] ({pic cs:part2-figure_p_prime_j1-p29R}) [red] to ({pic cs:part2-figure_p_prime_j1-p31L});
    \draw [->] ({pic cs:part2-figure_p_prime_j1-p29R}) [red,bend right=16] to ({pic cs:part2-figure_p_prime_j1-p32L});
    \draw [->] ({pic cs:part2-figure_p_prime_j1-p33R}) [red,bend left=4] to ({pic cs:part2-figure_p_prime_j1-p15L});
    \draw [->] ({pic cs:part2-figure_p_prime_j1-p33R}) [red,bend right=4] to ({pic cs:part2-figure_p_prime_j1-p34L});
    \draw [->] ({pic cs:part2-figure_p_prime_j1-p35R}) [red,bend right=4] to ({pic cs:part2-figure_p_prime_j1-p34L});
    \draw [->] ({pic cs:part2-figure_p_prime_j1-p36R}) [red] to ({pic cs:part2-figure_p_prime_j1-p37L});
    \draw [->] ({pic cs:part2-figure_p_prime_j1-p38R}) [red,bend left=16] to ({pic cs:part2-figure_p_prime_j1-p39L});
    \draw [->] ({pic cs:part2-figure_p_prime_j1-p38R}) [red] to ({pic cs:part2-figure_p_prime_j1-p40L});
    \draw [->] ({pic cs:part2-figure_p_prime_j1-p38R}) [red,bend right=16] to ({pic cs:part2-figure_p_prime_j1-p41L});
    \draw [->] ({pic cs:part2-figure_p_prime_j1-p42R}) [red,bend left=4] to ({pic cs:part2-figure_p_prime_j1-p17L});
    \draw [->] ({pic cs:part2-figure_p_prime_j1-p42R}) [red,bend right=4] to ({pic cs:part2-figure_p_prime_j1-p43L});
    \draw [->] ({pic cs:part2-figure_p_prime_j1-p44R}) [red,bend right=4] to ({pic cs:part2-figure_p_prime_j1-p43L});
    \draw [->] ({pic cs:part2-figure_p_prime_j1-p45R}) [red] to ({pic cs:part2-figure_p_prime_j1-p46L});
    \draw [->] ({pic cs:part2-figure_p_prime_j1-p47R}) [red] to ({pic cs:part2-figure_p_prime_j1-p48L});
    \draw [->] ({pic cs:part2-figure_p_prime_j1-p49R}) [red] to ({pic cs:part2-figure_p_prime_j1-p50L});
    \draw [->] ({pic cs:part2-figure_p_prime_j1-p51R}) [red] to ({pic cs:part2-figure_p_prime_j1-p52L});
    \draw [->] ({pic cs:part2-figure_p_prime_j1-p55R}) [red,bend left=16] to ({pic cs:part2-figure_p_prime_j1-p14L});
    \draw [->] ({pic cs:part2-figure_p_prime_j1-p55R}) [red] to ({pic cs:part2-figure_p_prime_j1-p33L});
    \draw [->] ({pic cs:part2-figure_p_prime_j1-p55R}) [red,bend right=16] to ({pic cs:part2-figure_p_prime_j1-p35L});
    \draw [->] ({pic cs:part2-figure_p_prime_j1-p56R}) [red,bend left=16] to ({pic cs:part2-figure_p_prime_j1-p16L});
    \draw [->] ({pic cs:part2-figure_p_prime_j1-p56R}) [red] to ({pic cs:part2-figure_p_prime_j1-p42L});
    \draw [->] ({pic cs:part2-figure_p_prime_j1-p56R}) [red,bend right=16] to ({pic cs:part2-figure_p_prime_j1-p44L});
\end{tikzpicture}
\par\endgroup
\caption{The BXSS in degree $(s,k,t) = (2,3,2^{3}\cdot 73+1)$ for \cref{CESS p prime}.
The vertical axis is the BX filtration and the horizontal axis is the $s$-degree. Red arrows indicate $d_1$-differentials, and blue arrows indicate $d_2$-differentials.}
\label{fig:figure_p_prime_j1}
\label{fig:figure_p_prime_j1-part2}
\end{figure}

    If $j=0$, the following classes and differentials changes:
    \[(2) d_1(q_{2}q_{6}^2\cdot h_{4}) = q_{1}q_{6}^2\cdot h_1h_{4} \]%d_1 boundary
    \[(5) d_1(q_1q_4q_6\cdot h_1h_6)=q_1q_3q_6\cdot h_1h_3h_6\]
    \[(6) d_1(q_1q_3q_7\cdot h_1h_3)=q_1q_3q_6\cdot h_1h_3h_6\]
    \[(7) q_0q_4q_7\cdot h_0h_1=0\]
    \[(8) d_1(q_2q_3^2\cdot h_7)=q_1q_3^2\cdot h_1h_7\]%d_1 boundary
    \[(11) \text{ is the same as } (2)\]
    \[(13) \text{ is the same as } (5)\]
    \[(14) \text{ is the same as } (6)\]
    \[(15) \text{ is the same as } (7)\]
    \[(16) d_1(q_1^2q_7\cdot h_3^2)=q_1^2q_6\cdot h_3^2h_6\]
    %\[(17) d_2(q_1^2q_6\cdot h_4h_6)=q_0^2q_6\cdot h_1h_4h_6\]
    \[(18) \text{ is the same as } (19)\]
    \[(20) \text{ is the same as } (8)\]
    \[(21) d_1(q_1^2q_5\cdot h_6^2)=q_1^2q_4\cdot h_4h_6^2\]
    %\[(22) d_2(q_1^2q_3\cdot h_3h_7)=q_0^2q_3\cdot h_1h_3h_7\]

    To account for the new linear relation $d_1(q_1q_4q_6\cdot h_1h_6+q_1q_3q_7\cdot h_1h_3)=0$, we have $d_1(q_1q_4q_7\cdot h_1)=q_1q_4q_6\cdot h_1h_6+q_1q_3q_7\cdot h_1h_3$

\begin{figure}[H]
\centering
\begingroup
\setlength{\tabcolsep}{6pt}
\renewcommand{\arraystretch}{1.1}
\begin{tabular}{|wc{1.2cm}||wc{2.5cm}|p{5cm}|wc{4cm}|} \hline
      $w$ & $s=1$ & \multicolumn{1}{p{5cm}|}{\centering $s=2$} & $s=3$ \\\hline\hline
% r = 16
    $16$ &  & \tikzmark{part1-figure_p_prime_j0-p1L}$(1)\; q_4q_6^2\cdot h_1^2$\tikzmark{part1-figure_p_prime_j0-p1R} &  \\\hline
% r = 15
    $15$ &  &  & \tikzmark{part1-figure_p_prime_j0-p2L}$q_3q_6^2\cdot h_1^2h_3$\tikzmark{part1-figure_p_prime_j0-p2R} \\\hline
% r = 14
    $14$ & \tikzmark{part1-figure_p_prime_j0-p3L}$q_2q_6^2\cdot h_4$\tikzmark{part1-figure_p_prime_j0-p3R} &  &  \\\hline
% r = 13
    $13$ &  & \tikzmark{part1-figure_p_prime_j0-p4L}$(2,11)\; q_1q_6^2\cdot h_1h_4$\tikzmark{part1-figure_p_prime_j0-p4R} &  \\
     &  & \tikzmark{part1-figure_p_prime_j0-p5L}$(3)\; q_3^2q_7\cdot h_1^2$\tikzmark{part1-figure_p_prime_j0-p5R} &  \\\hline
% r = 12
    $12$ & \tikzmark{part1-figure_p_prime_j0-p34L}$q_1q_4q_7\cdot h_1$\tikzmark{part1-figure_p_prime_j0-p34R} & \tikzmark{part1-figure_p_prime_j0-p7L}$(4,12)\; q_1q_4q_7\cdot h_0^2$\tikzmark{part1-figure_p_prime_j0-p7R} & \tikzmark{part1-figure_p_prime_j0-p6L}$q_3^2q_6\cdot h_1^2h_6$\tikzmark{part1-figure_p_prime_j0-p6R} \\\hline
% r = 11
    $11$ &  & \tikzmark{part1-figure_p_prime_j0-p11L}$(5,13)\; q_1q_4q_6\cdot h_1h_6$\tikzmark{part1-figure_p_prime_j0-p11R} & \tikzmark{part1-figure_p_prime_j0-p8L}$q_0q_4q_7\cdot h_0^3$\tikzmark{part1-figure_p_prime_j0-p8R} \\
     &  & \tikzmark{part1-figure_p_prime_j0-p13L}$(6,14)\; q_1q_3q_7\cdot h_1h_3$\tikzmark{part1-figure_p_prime_j0-p13R} & \tikzmark{part1-figure_p_prime_j0-p9L}$q_1q_3q_7\cdot h_0^2h_3$\tikzmark{part1-figure_p_prime_j0-p9R} \\
     &  &  & \tikzmark{part1-figure_p_prime_j0-p10L}$q_1q_4q_6\cdot h_0^2h_6$\tikzmark{part1-figure_p_prime_j0-p10R} \\\hline
% r = 10
    $10$ &  &  & \tikzmark{part1-figure_p_prime_j0-p12L}$q_1q_3q_6\cdot h_1h_3h_6$\tikzmark{part1-figure_p_prime_j0-p12R} \\\hline
% r = 9
    $9$ & \tikzmark{part1-figure_p_prime_j0-p35L}$q_1^2q_7\cdot h_4$\tikzmark{part1-figure_p_prime_j0-p35R} & \tikzmark{part1-figure_p_prime_j0-p20L}$(16)\; q_1^2q_7\cdot h_3^2$\tikzmark{part1-figure_p_prime_j0-p20R} &  \\\hline
% r = 8
    $8$ & \tikzmark{part1-figure_p_prime_j0-p14L}$q_2q_3^2\cdot h_7$\tikzmark{part1-figure_p_prime_j0-p14R} & \tikzmark{part1-figure_p_prime_j0-p22L}$(17)\; q_1^2q_6\cdot h_4h_6$\tikzmark{part1-figure_p_prime_j0-p22R} & \tikzmark{part1-figure_p_prime_j0-p21L}$q_1^2q_6\cdot h_3^2h_6$\tikzmark{part1-figure_p_prime_j0-p21R} \\
     &  & \tikzmark{part1-figure_p_prime_j0-p23L}$(18,19)\; q_0q_1q_7\cdot h_0h_4$\tikzmark{part1-figure_p_prime_j0-p23R} &  \\\hline
% r = 7
    $7$ &  & \tikzmark{part1-figure_p_prime_j0-p15L}$(8,20)\; q_1q_3^2\cdot h_1h_7$\tikzmark{part1-figure_p_prime_j0-p15R} & \tikzmark{part1-figure_p_prime_j0-p24L}$q_0q_1q_6\cdot h_0h_4h_6$\tikzmark{part1-figure_p_prime_j0-p24R} \\
     &  & \tikzmark{part1-figure_p_prime_j0-p16L}$(9,26)\; q_0^2q_7\cdot h_1h_4$\tikzmark{part1-figure_p_prime_j0-p16R} & \tikzmark{part1-figure_p_prime_j0-p25L}$q_0^2q_7\cdot h_0^2h_4$\tikzmark{part1-figure_p_prime_j0-p25R} \\
     &  & \tikzmark{part1-figure_p_prime_j0-p26L}$(21)\; q_1^2q_5\cdot h_6^2$\tikzmark{part1-figure_p_prime_j0-p26R} &  \\\hline
% r = 6
    $6$ & \tikzmark{part1-figure_p_prime_j0-p36L}$q_1^2q_4\cdot h_7$\tikzmark{part1-figure_p_prime_j0-p36R} &  & \tikzmark{part1-figure_p_prime_j0-p17L}$q_0^2q_6\cdot h_1h_4h_6$\tikzmark{part1-figure_p_prime_j0-p17R} \\
     &  &  & \tikzmark{part1-figure_p_prime_j0-p27L}$q_1^2q_4\cdot h_4h_6^2$\tikzmark{part1-figure_p_prime_j0-p27R} \\\hline
% r = 5
    $5$ &  & \tikzmark{part1-figure_p_prime_j0-p28L}$(22)\; q_1^2q_3\cdot h_3h_7$\tikzmark{part1-figure_p_prime_j0-p28R} &  \\
     &  & \tikzmark{part1-figure_p_prime_j0-p29L}$(23,24)\; q_0q_1q_4\cdot h_0h_7$\tikzmark{part1-figure_p_prime_j0-p29R} &  \\\hline
% r = 4
    $4$ &  & \tikzmark{part1-figure_p_prime_j0-p18L}$(10,27)\; q_0^2q_4\cdot h_1h_7$\tikzmark{part1-figure_p_prime_j0-p18R} & \tikzmark{part1-figure_p_prime_j0-p30L}$q_0q_1q_3\cdot h_0h_3h_7$\tikzmark{part1-figure_p_prime_j0-p30R} \\
     &  &  & \tikzmark{part1-figure_p_prime_j0-p31L}$q_0^2q_4\cdot h_0^2h_7$\tikzmark{part1-figure_p_prime_j0-p31R} \\\hline
% r = 3
    $3$ &  &  & \tikzmark{part1-figure_p_prime_j0-p19L}$q_0^2q_3\cdot h_1h_3h_7$\tikzmark{part1-figure_p_prime_j0-p19R} \\\hline
% r = 2
    $2$ &  &  &  \\\hline
% r = 1
    $1$ &  & \tikzmark{part1-figure_p_prime_j0-p32L}$(25,28)\; q_0^2q_1\cdot h_4h_7$\tikzmark{part1-figure_p_prime_j0-p32R} &  \\\hline
% r = 0
    $0$ &  &  & \tikzmark{part1-figure_p_prime_j0-p33L}$q_0^3\cdot h_0h_4h_7$\tikzmark{part1-figure_p_prime_j0-p33R} \\\hline
    \end{tabular}
\begin{tikzpicture}[overlay, remember picture, shorten >=.5pt, shorten <=.5pt, transform canvas={yshift=.25\baselineskip}]

    \draw [->] ({pic cs:part1-figure_p_prime_j0-p1R}) [red] to ({pic cs:part1-figure_p_prime_j0-p2L});
    \draw [->] ({pic cs:part1-figure_p_prime_j0-p3R}) [red] to ({pic cs:part1-figure_p_prime_j0-p4L});
    \draw [->] ({pic cs:part1-figure_p_prime_j0-p5R}) [red] to ({pic cs:part1-figure_p_prime_j0-p6L});
    \draw [->] ({pic cs:part1-figure_p_prime_j0-p7R}) [red,bend left=16] to ({pic cs:part1-figure_p_prime_j0-p8L});
    \draw [->] ({pic cs:part1-figure_p_prime_j0-p7R}) [red] to ({pic cs:part1-figure_p_prime_j0-p9L});
    \draw [->] ({pic cs:part1-figure_p_prime_j0-p7R}) [red,bend right=16] to ({pic cs:part1-figure_p_prime_j0-p10L});
    \draw [->] ({pic cs:part1-figure_p_prime_j0-p11R}) [red,bend left=4] to ({pic cs:part1-figure_p_prime_j0-p12L});
    \draw [->] ({pic cs:part1-figure_p_prime_j0-p13R}) [red,bend right=4] to ({pic cs:part1-figure_p_prime_j0-p12L});
    \draw [->] ({pic cs:part1-figure_p_prime_j0-p14R}) [red] to ({pic cs:part1-figure_p_prime_j0-p15L});
    \draw [->] ({pic cs:part1-figure_p_prime_j0-p16R}) [red,bend right=4] to ({pic cs:part1-figure_p_prime_j0-p17L});
    \draw [->] ({pic cs:part1-figure_p_prime_j0-p18R}) [red,bend right=4] to ({pic cs:part1-figure_p_prime_j0-p19L});
    \draw [->] ({pic cs:part1-figure_p_prime_j0-p20R}) [red] to ({pic cs:part1-figure_p_prime_j0-p21L});
    \draw [->] ({pic cs:part1-figure_p_prime_j0-p23R}) [red,bend left=8] to ({pic cs:part1-figure_p_prime_j0-p24L});
    \draw [->] ({pic cs:part1-figure_p_prime_j0-p23R}) [red,bend right=8] to ({pic cs:part1-figure_p_prime_j0-p25L});
    \draw [->] ({pic cs:part1-figure_p_prime_j0-p26R}) [red] to ({pic cs:part1-figure_p_prime_j0-p27L});
    \draw [->] ({pic cs:part1-figure_p_prime_j0-p29R}) [red,bend left=8] to ({pic cs:part1-figure_p_prime_j0-p30L});
    \draw [->] ({pic cs:part1-figure_p_prime_j0-p29R}) [red,bend right=8] to ({pic cs:part1-figure_p_prime_j0-p31L});
    \draw [->] ({pic cs:part1-figure_p_prime_j0-p32R}) [red] to ({pic cs:part1-figure_p_prime_j0-p33L});
    \draw [->] ({pic cs:part1-figure_p_prime_j0-p34R}) [red,bend left=8] to ({pic cs:part1-figure_p_prime_j0-p11L});
    \draw [->] ({pic cs:part1-figure_p_prime_j0-p34R}) [red,bend right=8] to ({pic cs:part1-figure_p_prime_j0-p13L});
    \draw [->] ({pic cs:part1-figure_p_prime_j0-p35R}) [red] to ({pic cs:part1-figure_p_prime_j0-p22L});
    \draw [->] ({pic cs:part1-figure_p_prime_j0-p36R}) [red] to ({pic cs:part1-figure_p_prime_j0-p28L});
\end{tikzpicture}
\par\endgroup
\caption{The BXSS in degree $(s,k,t) = (2,3,2^{2}\cdot 73+1)$ for \cref{CESS p prime}.
The vertical axis is the BX filtration and the horizontal axis is the $s$-degree. Red arrows indicate $d_1$-differentials, and blue arrows indicate $d_2$-differentials.}
\label{fig:figure_p_prime_j0}
\end{figure}

    Therefore, none of these classes survives to CESS, and there are no classes supporting differentials hitting $q_0\cdot T_j.$ 
    
    \end{proof}

\section{The $f_j$ Family}
We prove the following theorem:
\begin{thm}
There are Adams differentials
\[
d_3(f_{j+2}) = h_1h_{j+4}Q_3(j) \neq 0
\]
for $j\geq 0$.
\end{thm}

Recall that the form of the differential was shown to hold in \Cref{sec:differentials}, so all that remains is to show that the target is nonzero on the Adams $E_3$ page. The proof proceeds in two steps, following the program of Burklund-Xu. First we show, for $j\geq1$, that the element $h_0h_{j+3}Q_3(j-1)$ is a permanent cycle in the Burklund-Xu spectral sequence and detects a nonzero class in the abutment, an associated graded of the Cartan-Eilenberg $E_2$ page. We then show that $h_0h_{j+3}Q_3(j-1)$ is a permanent cycle in that spectral sequence and detects a nonzero class in the abutment, an associated graded of the Adams $E_2$ page. Since $h_0h_{j+3}Q_3(j-1)$ detects the classical Adams element $h_1h_{j+4}Q_3(j)$ under the doubling-up isomorphism
\[
\text{Ext}_\mathcal{A}\xrightarrow{\cong}\text{Ext}_\mathcal{P},
\]
we will have shown that $h_1h_{j+4}Q_3(j)$ is nonzero on the Adams $E_2$ page for $j\geq1$. We then analyze possible sources of $d_2$ differentials hitting these classes and show that no such differentials occur, so the $j\geq1$ case will follow. The $j=0$ case can be read off of \cite{charts}.

To see that $f_{j+1}$ survives to the Adams $E_3$ page, we use the differentials
\[
d_2(f_{j+1}) = h_0\cdot \text{Sq}^1(f_j) = 0
\]
for $j\geq1$ (see \Cref{sec:differentials}).

We now need to show that $h_1h_{j+4}Q_3(j)$ detects a nonzero class on the Adams $E_3$ page. First note that, since $h_1h_{j+4}Q_3(j)$ is the target of an Adams differential, it cannot support any Adams differentials. In particular, it must be a $d_2$ cycle so survives to the Adams $E_3$ page, even if it detects a trivial class there.

Note that $h_1h_4Q_3(0)\neq0$ on the Adams $E_3$ page, by \cite{charts}, so the $j=0$ case is handled. In order to show that $h_1h_{j+4}Q_3(j)$ is nontrivial on the Adams $E_2$ page for $j\geq1$, we will study the class
\[
h_0h_{j+4}Q_3(j)
\]
in the Burklund-Xu and Cartan-Eilenberg spectral sequences for $j\geq0$. We first show that this class survives the Burklund-Xu spectral sequence.

\begin{prop}
The class $h_0h_{j+4}Q_3(j)$ survives the Burklund-Xu spectral sequence to detect a nonzero class on the Cartan-Eilenberg $E_2$ page for all $j\geq0$.
\end{prop}
\begin{proof}
Since $h_0h_{j+4}Q_3(j)$ lies in Burklund-Xu degree 0, it cannot support any differentials in that spectral sequence. Furthermore, by the splitting described in \autoref{sec:recollections}, the only classes which could support differentials hitting $h_0h_{j+4}Q_3(j)$ similarly lie in Burklund-Xu degree 0 so cannot support any differentials, either. Then $h_0h_{j+4}Q_3(j)$ supports and receives no Burklund-Xu differentials so detects a nonzero class on the Cartan-Eilenberg $E_2$ page.
\end{proof}

We now determine the fate of this family in the Cartan-Eilenberg spectral sequence.

\begin{prop}
\label{CESS f}
The class $h_0h_{j+4}Q_3(j)$ survives the Cartan-Eilenberg spectral sequence to detect a nonzero class on the Adams $E_2$ page for all $j\geq1$. %j=0 case is given by Weinan's chart
\end{prop}
\begin{proof}
We first note that $h_0h_{j+4}Q_3(j)$ lies in Cartan-Eilenberg filtration 0 so cannot \textit{support} any Cartan-Eilenberg differentials, as every differential lowers filtration. Then it remains to show that this class does not \textit{receive} any differentials, either. This will follow from \Cref{prop:no-d3}, showing that it does not receive any Cartan-Eilenberg $d_3$ differential, and \Cref{prop:no-d5}, showing that it does not receive any Cartan-Eilenberg $d_5$ differential.

Indeed, since $h_0h_{j+4}Q_3(j)$ lies in Cartan-Eilenberg filtration 0, any differential source must lie in filtration 2, 4, or 6. We immediately rule out filtration 6 since the Ext filtration of such a class would be $s = 6-6 = 0$. By \cite{Benson}, the only such class on the Cartan-Eilenberg $E_2$ page is $q_0^6$, which supports no differentials because it is the only class which can detect $h_0^6$ on the Adams $E_2$ page.
\end{proof}

\begin{prop}
\label{prop:no-d3}
The class $h_0h_{j+4}Q_3(j)$ does not receive a Cartan-Eilenberg $d_3$ differential for any $j\geq1$.
\end{prop}
\begin{proof}
If $h_0h_{j+4}Q_3(j)$ receives a $d_3$ differential
\[
    d_3(x) = h_0h_{j+4}Q_3(j),
\]
then $x$ must lie in Cartan-Eilenberg filtration $k = 2$. The Ext component of $x$ must have filtration $s = 4$ since $h_0h_{j+4}Q_3(j)$ lies in total filtration $k+s = 6+1$. A general class in $(k,s) = (2,4)$ takes one of the following two forms:
\[
\begin{array}{c c c}
q_aq_b\cdot h_ch_dh_eh_f & q_aq_b\cdot h_dc_e & q_aq_b\cdot x_c,
\end{array}
\]
where $x_c$ can be $d_c,e_c,f_c,g_c,p_c,D_3(c),p'_c$.\\
In the first case, the total degree
\[
t(q_aq_b\cdot h_ch_dh_eh_f) = 2^{a+1}+2^{b+1}-2+2^{c+1}+2^{d+1}+2^{e+1}+2^{f+1},
\]
where we assume without loss of generality that $a\leq b$ and $c\leq d\leq e\le f$. Setting this equal to the total degree of $h_0h_{j+4}Q_3(j)$, we obtain
\[
2^a+2^b+2^c+2^d+2^e+2^f = 2^1+2^{j+3}+2^{j+4}+2^{j+6}
\]
after adding 2 and dividing by 2.\\
Up to reordering, $(a,b,c,d,e,f)$ can be one of the following:
\begin{itemize}
    \item $(0,0,j+2,j+2,j+4,j+6)$
    \item $(0,0,j+3,j+3,j+3,j+6)$
    \item $(0,0,j+3,j+4,j+5,j+5)$
    \item $(1,j+2,j+2,j+3,j+3,j+6)$
    \item $(1,j+2,j+2,j+4,j+5,j+5)$
    \item $(1,j+3,j+3,j+3,j+5,j+5)$
\end{itemize}
The possible classes are 
\begin{itemize}
    \item $q_{j+4}q_{j+6}\cdot h_0^2h_{j+2}^2$
    \item $q_{j+5}^2\cdot h_0^2h_{j+3}h_{j+4} = 0$
    \item $q_{j+5}^2\cdot h_1h_{j+2}^2h_{j+4}=q_{j+5}^2\cdot h_1h_{j+3}^3$ %class is 0 when j=1
    \item $q_{j+3}q_{j+6}\cdot h_0^2h_{j+3}^2$
    \item $q_{j+4}q_{j+5}\cdot h_0^2h_{j+3}h_{j+5}$
    \item $q_{j+2}q_{j+6}\cdot h_0^2h_{j+2}h_{j+4}$
    \item $q_{j+2}q_{j+6}\cdot h_0^2h_{j+2}h_{j+3}^2=0$
    \item $q_{j+2}q_{j+6}\cdot h_1h_{j+2}h_{j+3}^2=0$
    \item $q_{j+3}q_{j+5}\cdot h_1h_{j+3}^2h_{j+5}=q_{j+3}q_{j+5}\cdot h_1h_{j+4}^3$
    \item $q_{j+2}q_{j+6}\cdot h_1h_{j+1}^2h_{j+4}$
    \item $q_{j+3}q_{j+5}\cdot h_0^2h_{j+4}h_{j+5}=0$
    \item $q_{j+2}q_{j+5}\cdot h_1h_{j+2}h_{j+4}h_{j+5}=0$
    \item $q_{j+3}q_{j+4}\cdot h_0^2h_{j+5}^2$
    \item $q_{j+2}q_{j+4}\cdot h_0^2h_{j+2}h_{j+6}$
    \item $q_{j+2}q_{j+4}\cdot h_1h_{j+2}h_{j+5}^2$
    \item $q_{j+3}^2\cdot h_0^2h_{j+3}h_{j+6}$
    \item $q_{j+3}^2\cdot h_1h_{j+2}^2h_{j+6}$ %class is 0 when j=1
    \item $q_{j+3}^2\cdot h_1h_{j+3}h_{j+5}^2=0$
    \item $q_{j+2}q_{j+3}\cdot h_1h_{j+2}h_{j+3}h_{j+6}=0$
    \item $q_{j+2}^2\cdot h_0^2h_{j+4}h_{j+6}$
    \item $q_{j+2}^2\cdot h_1h_{j+3}^2h_{j+6}$
    \item $q_{j+2}^2\cdot h_1h_{j+4}h_{j+5}^2=0$
    \item $q_1q_{j+6}\cdot h_{j+2}^2h_{j+3}^2=0$
    \item $q_0q_{j+6}\cdot h_0h_{j+2}^2h_{j+4}=q_0q_{j+6}\cdot h_0h_{j+3}^3$ 
    \item $q_1q_{j+5}\cdot h_{j+2}^2h_{j+4}h_{j+5}=0$
    \item $q_1q_{j+5}\cdot h_{j+3}^3h_{j+5}=q_1q_{j+5}\cdot h_{j+3}h_{j+4}^3=0$
    \item $q_0q_{j+5}\cdot h_0h_{j+3}h_{j+4}h_{j+5}=0$
    \item $q_1q_{j+4}\cdot h_{j+2}^2h_{j+5}^2=0$
    \item $q_0q_{j+4}\cdot h_0h_{j+2}^2h_{j+6}$
    \item $q_0q_{j+4}\cdot h_0h_{j+3}h_{j+5}^2=0$
    \item $q_1q_{j+3}\cdot h_{j+3}^2h_{j+5}^2=0$
    \item $q_0q_{j+3}\cdot h_0h_{j+3}^2h_{j+6}$;
    \item $q_0q_{j+3}\cdot h_0h_{j+4}h_{j+5}^2=0$
    \item $q_1q_{j+2}\cdot h_{j+2}h_{j+3}^2h_{j+6}=0$
    \item $q_1q_{j+2}\cdot h_{j+2}h_{j+4}h_{j+5}^2=0$
    \item $q_0q_{j+2}\cdot h_0h_{j+2}h_{j+4}h_{j+6}$
    \item $q_0^2\cdot h_{j+3}^3h_{j+6}$
    \item $q_0^2\cdot h_{j+3}h_{j+4}h_{j+5}^2=0$
    \item $q_{j+4}q_{j+5}\cdot h_1h_{j+2}^2h_{j+5}$
    \item $q_{j+2}q_{j+4}\cdot h_1h_{j+1}^2h_{j+5}$
    \item $q_{j+1}q_{j+2}\cdot h_1h_{j+1}h_{j+4}h_{j+6}$
    \item $q_{j+1}^2\cdot h_1h_{j+2}h_{j+4}h_{j+6}$
    \item $q_1q_{j+2}\cdot h_{j+1}^2h_{j+4}h_{j+6}$ 
\end{itemize}
Most of these classes supports or killed by BX differentials:
\begin{enumerate}
    \item $d_1(q_{j+4}q_{j+6}\cdot h_0^2h_{j+2}^2)=q_{j+4}q_{j+5}\cdot h_0^2h_{j+2}^2h_{j+5}$ %=0 if j=0
    %\item $q_{j+5}^2\cdot h_0^2h_{j+3}h_{j+4} = 0$
    \item $q_{j+5}^2\cdot h_1h_{j+2}^2h_{j+4}=q_{j+5}^2\cdot h_1h_{j+3}^3$; $d_6(q_{j+5}^2\cdot h_1h_{j+3}^3)=q_{j+2}^2\cdot \langle h_{j+2}, h_{j+3}, h_{j+4}, h_1h_{j+3}^3 \rangle = q_{j+2}^2\cdot h_1h_{j+3}c_{j+2}$ %class is 0 when j=1
    %\item $q_{j+5}^2\cdot h_1h_{j+3}^3$
    \item $d_1(q_{j+3}q_{j+6}\cdot h_0^2h_{j+3}^2)=q_{j+3}q_{j+5}\cdot h_0^2h_{j+3}^2h_{j+5}=q_{j+3}q_{j+5}\cdot h_0^2h_{j+4}^3$
    \item $d_1(q_{j+4}q_{j+5}\cdot h_0^2h_{j+3}h_{j+5})=q_{j+3}q_{j+5}\cdot h_0^2h_{j+3}^2h_{j+5}=q_{j+3}q_{j+5}\cdot h_0^2h_{j+4}^3$
    \item $d_1(q_{j+3}q_{j+6}\cdot h_0^2h_{j+4})=q_{j+2}q_{j+6}\cdot h_0^2h_{j+2}h_{j+4}$
    %$d_2(q_{j+2}q_{j+6}\cdot h_0^2h_{j+2}h_{j+4})=q_{j+2}q_{j+4}\cdot h_0^2h_{j+2}h_{j+5}^2$ %d_1 target
    %\item $q_{j+2}q_{j+6}\cdot h_0^2h_{j+2}h_{j+3}^2=0$
    %\item $q_{j+2}q_{j+6}\cdot h_1h_{j+2}h_{j+3}^2=0$
    \item $q_{j+3}q_{j+5}\cdot h_1h_{j+3}^2h_{j+5}=q_{j+3}q_{j+5}\cdot h_1h_{j+4}^3$; $d_1(q_{j+3}q_{j+6}\cdot h_1h_{j+3}^2)=q_{j+3}q_{j+5}\cdot h_1h_{j+3}^2h_{j+5}$
    \item $d_1(q_{j+2}q_{j+6}\cdot h_1h_{j+1}^2h_{j+4})=q_{j+1}q_{j+6}\cdot h_1h_{j+1}^3h_{j+4}$
    %$d_1(q_{j+3}q_{j+5}\cdot h_1h_{j+4}^3)=q_{j+3}q_{j+4}\cdot h_1h_{j+4}^4$
    %\item $q_{j+3}q_{j+5}\cdot h_0^2h_{j+4}h_{j+5}=0$
    %\item $q_{j+2}q_{j+5}\cdot h_1h_{j+2}h_{j+4}h_{j+5}=0$
    \item $d_1(q_{j+3}q_{j+4}\cdot h_0^2h_{j+5}^2)=q_{j+2}q_{j+4}\cdot h_0^2h_{j+2}h_{j+5}^2$
    \item $d_2(q_{j+2}q_{j+4}\cdot h_0^2h_{j+2}h_{j+6})=q_{j+2}^2\cdot h_0^2h_{j+3}^2h_{j+6}$
    \item $d_1(q_{j+3}q_{j+4}\cdot h_1h_{j+5}^2) = q_{j+2}q_{j+4}\cdot h_1h_{j+2}h_{j+5}^2$%+q_{j+3}^2\cdot h_1h_{j+3}h_{j+5}^2$
    %\item $d_3(q_{j+2}q_{j+4}\cdot h_1h_{j+2}h_{j+5}^2)=q_{j+1}q_{j+2}\cdot \langle h_{j+1}, h_{j+2}, h_{j+3}, h_1h_{j+2}h_{j+5}^2\rangle$ %boundary
    \item $d_2(q_{j+3}^2\cdot h_0^2h_{j+3}h_{j+6})=q_{j+2}^2\cdot h_0^2h_{j+3}^2h_{j+6}$
    \item $d_6(q_{j+3}^2\cdot h_1h_{j+2}^2h_{j+6})=q_{j}^2\cdot \langle h_{j+1}, h_{j+2}, h_{j+3}, h_1h_{j+2}^2h_{j+6} \rangle = q_j^2\cdot h_1h_{j+6}c_{j+1}$ %class is 0 when j=1
    %\item $q_{j+3}^2\cdot h_1h_{j+3}h_{j+5}^2=0$
    %\item $q_{j+2}q_{j+3}\cdot h_1h_{j+2}h_{j+3}h_{j+6}=0$
    \item $d_2(q_{j+2}^2\cdot h_0^2h_{j+4}h_{j+6}) = q_{j+1}^2\cdot h_0^2h_{j+2}h_{j+4}h_{j+6}$
    \item $q_{j+2}^2\cdot h_1h_{j+3}^2h_{j+6}$; $d_2(q_{j+3}^2\cdot h_1h_{j+3}h_{j+6})=q_{j+2}^2\cdot h_1h_{j+3}^2h_{j+6}$
    %\item $q_{j+2}^2\cdot h_1h_{j+4}h_{j+5}^2=0$
    %\item $q_1q_{j+6}\cdot h_{j+2}^2h_{j+3}^2=0$
    \item $q_0q_{j+6}\cdot h_0h_{j+2}^2h_{j+4}=q_0q_{j+6}\cdot h_0h_{j+3}^3$; $d_1(q_1q_{j+6}\cdot h_{j+3}^3)=q_0q_{j+6}\cdot h_0h_{j+3}^3$ %q_0q_{j+4}\cdot \langle h_{j+4}, h_{j+5},  h_{j+3}^3\rangle
    %\item $q_0q_{j+6}\cdot h_0h_{j+3}^3$
    %\item $q_1q_{j+5}\cdot h_{j+2}^2h_{j+4}h_{j+5}=0$
    %\item $q_1q_{j+5}\cdot h_{j+3}^3h_{j+5}=q_1q_{j+5}\cdot h_{j+3}h_{j+4}^3=0$
    %\item $q_0q_{j+5}\cdot h_0h_{j+3}h_{j+4}h_{j+5}=0$
    %\item $q_1q_{j+4}\cdot h_{j+2}^2h_{j+5}^2=0$
    \item $d_1(q_1q_{j+4}\cdot h_{j+2}^2h_{j+6})=q_0q_{j+4}\cdot h_0h_{j+2}^2h_{j+6}$
    %\item $q_0q_{j+4}\cdot h_0h_{j+3}h_{j+5}^2=0$
    %\item $q_1q_{j+3}\cdot h_{j+3}^2h_{j+5}^2=0$
    \item $q_0q_{j+3}\cdot h_0h_{j+3}^2h_{j+6}$; $d_1(q_0q_{j+4}\cdot h_0h_{j+3}h_{j+6})=q_0q_{j+3}\cdot h_0h_{j+3}^2h_{j+6}$
    %\item $q_0q_{j+3}\cdot h_0h_{j+4}h_{j+5}^2=0$
    %\item $q_1q_{j+2}\cdot h_{j+2}h_{j+3}^2h_{j+6}=0$
    %\item $q_1q_{j+2}\cdot h_{j+2}h_{j+4}h_{j+5}^2=0$
    \item $q_0q_{j+2}\cdot h_0h_{j+2}h_{j+4}h_{j+6}$; $d_1(q_0q_{j+3}\cdot h_0h_{j+4}h_{j+6})=q_0q_{j+2}\cdot h_0h_{j+2}h_{j+4}h_{j+6}$
    %\item $q_0^2\cdot h_{j+3}^3h_{j+6}$
    %\item $q_0^2\cdot h_{j+3}h_{j+4}h_{j+5}^2=0$
    \item $d_1(q_{j+4}q_{j+6}\cdot h_1h_{j+2}^2)=q_{j+4}q_{j+5}\cdot h_1h_{j+2}^2h_{j+5}$
    \item $d_1(q_{j+2}q_{j+4}\cdot h_1h_{j+1}^2h_{j+5})=q_{j+1}q_{j+4}\cdot h_1h_{j+1}^3h_{j+5}+q_{j+2}q_{j+3}\cdot h_1h_{j+1}^2h_{j+3}h_{j+5}$
    \item $d_1(q_{j+1}q_{j+2}\cdot h_1h_{j+1}h_{j+4}h_{j+6})=q_{j+1}^2\cdot h_1h_{j+1}^2h_{j+4}h_{j+6}$
    \item $d_2(q_{j+2}^2\cdot h_1h_{j+4}h_{j+6})=q_{j+1}^2\cdot h_1h_{j+2}h_{j+4}h_{j+6}$
    \item $d_1(q_1q_{j+2}\cdot h_{j+1}^2h_{j+4}h_{j+6})=q_0q_{j+2}\cdot h_0h_{j+1}^2h_{j+4}h_{j+6}+q_1q_{n+1}\cdot h_{j+1}^3h_{j+4}h_{j+6}$
\end{enumerate}
To account for the linear relation $d_2(q_{j+2}q_{j+4}\cdot h_0^2h_{j+2}h_{j+6}+q_{j+3}^2\cdot h_0^2h_{j+3}h_{j+6})=0$, we need the differentials $d_1(q_{j+3}q_{j+4}\cdot h_0^2h_{j+6})=q_{j+2}q_{j+4}\cdot h_0^2h_{j+2}h_{j+6}+q_{j+3}^2\cdot h_0^2h_{j+3}h_{j+6}$\\
To account for the linear relation $d_1(q_{j+3}q_{j+6}\cdot h_0^2h_{j+3}^2+q_{j+4}q_{j+5}\cdot h_0^2h_{j+3}h_{j+5})=0$, we need the differential $d_1(q_{j+4}q_{j+6}\cdot h_0^2h_{j+3})=q_{j+3}q_{j+6}\cdot h_0^2h_{j+3}^2+q_{j+4}q_{j+5}\cdot h_0^2h_{j+3}h_{j+5}$.\\
For $j \ge 2$, all the classes above are distinct and nonzero. When $j=1$, $q_{j+5}^2\cdot h_1h_{j+3}^3$ and $q_{j+3}^2\cdot h_1h_{j+2}^2h_{j+6}$ are 0, so we don't need (2) or (11).

The class $q_0^2\cdot h_{j+2}^2h_{j+4}h_{j+6}=q_0^2\cdot h_{j+3}^3h_{j+6}$ survives to CESS, but it cannot support CESS differentials as products of permanent cycles $q_0$ and $h_j$ detecting $h_0$ and $h_{j+1}$ on Adams $E_2$ page.

In the second case, the total degree
\[
t(q_aq_b\cdot h_dc_e)=2^{a+1}+2^{b+1}-2+2^{d+1}+2^{e+1}+2^{e+2}+2^{e+4} = 2+2^{j+4}+2^{j+5}+2^{j+7}.
\]
%$2^{e+4}>4$, so we must have $e+4\le j+4,e+4=j+5,e+4=j+6$, or $e+4=j+7$.\\
If $e+4\leq j+4$, then $2^{e+1}+2^{e+2} \leq 4$, which is impossible.\\ 
%and $e=0$, so the only possibility is $e=j=0,a=3,b=5,d=0$.\\ %give an argument
%If $e+3=j+1$ and $j\ge 4$, then $2^{e+1}>4$, so we must have $e+1=j+3,e+1=j+4,e+1=j+6$, which are all impossible.\\
%If $e+3=j+3$ and $j\ge 2$, then $2^{e+1}>4$, so we must have $e+1=j+3,e+1=j+4,e+1=j+6$, which are all impossible.\\
If $e+4=j+5$, then $2^{j+2}+2^{j+3}+2^{j+5}+2^{a+1}+2^{b+1}+2^{d+1}=4+2^{j+4}+2^{j+5}+2^{j+7}$, and the possibilities are 
\begin{enumerate}
    \item $a=1,b=j+1,d=j+6$
    \item $a=1,b=j+6,d=j+1$
    \item $a=j+1,b=j+6,d=1$
    %\item $e=1,j=0,a=0,b=1,d=6$
    %\item $e=1,j=0,a=0,b=5,d=2$
    %\item $e=1,j=0,a=1,b=5,d=1$
    %\item $e=2,j=1,a=1,b=4,d=5$
    %\item $e=2,j=1,a=4,b=4,d=3$
\end{enumerate}
%$2^{e+1}>4$, so we must have $e+1=j+3,e+1=j+4,e+1=j+6$, which are all impossible.\\
If $e+4=j+6$, then $2^{j+3}+2^{j+4}+2^{j+6}+2^{a+1}+2^{b+1}+2^{d+1}=4+2^{j+4}+2^{j+5}+2^{j+7}$, which is impossible because the digit difference is greater than 2.\\
If $e+4=j+7$, then $2^{j+4}+2^{j+5}+2^{j+7}+2^{a+1}+2^{b+1}+2^{d+1}=4+2^{j+4}+2^{j+5}+2^{j+7}$, which is impossible.\\
In conclusion, the only possible classes killing $h_0h_{j+4}Q_3(j)$ for $j\geq 1$ in this case are
\begin{enumerate}
    \item $d_1(q_1q_{j+1}\cdot h_{j+6}c_{j+1})=q_0q_{j+1}\cdot h_0h_{j+5}c_{j+1}$
    \item $q_1q_{j+6}\cdot h_{j+1}c_{j+1}=0$
    \item $d_1(q_{j+1}q_{j+6}\cdot h_1c_{j+1})=q_{j+1}q_{j+5}\cdot h_1h_{j+5}c_{j+1}$

    %\item $q_3q_5\cdot h_0c_0$
    %\item $q_0q_1\cdot h_6c_1$
    %\item $q_0q_5\cdot h_2c_1$
    %\item $q_1q_5\cdot h_1c_1$
    %\item $q_1q_4\cdot h_5c_2=0$
    %\item $d_6(q_4^2\cdot h_3c_2)=q_1^2\cdot \langle h_2,h_3,h_4,h_3c_2 \rangle = q_1^2h_3d_2$ %j=1
\end{enumerate}
%Hence, the equation can only be satisfied when $j=0,1,2$, and the possibilities are $a=3,b=5,d=0,e=0$ and $a=0,b=1,d=6,e=1$, when $j=0$.
%exaust other possibilities
%$a=b=0$ and $d=j+3$, and $q_0^2\cdot c_{j+3}$ is a permanent cycle.

In the third case, the total degree
\[t(q_aq_b\cdot d_c)=2^{a+1}+2^{b+1}-2+2^{c+5}+2^{c+2}=2+2^{j+4}+2^{j+5}+2^{j+7}\]
After moving $-2$ to right hand side, we find both sides having 4 digits in the binary expansion, so the digits must match up exactly. As $2^{c+5}>4$, we must have $c+5=j+4,c+5=j+5$, or $c+5=j+7$.\\
If $c+5 = j+4$, then $2^{c+2}=4$. Hence, $c=0,j=1,a=5,b=7$.\\
If $c+5=j+5$, then $c=j=0$, which is not in the range of consideration.\\
If $c+5=j+7$, then the only possibilities are $a=1,b=j+4$.\\

\[t(q_aq_b\cdot e_c)=2^{a+1}+2^{b+1}-2+2^{c+5}+2^{c+3}+2^{c+1}=2+2^{j+4}+2^{j+5}+2^{j+7}\]
%$2^{c+5}>4$, so we must have $c+5=j+2,c+5=j+3,c+5=j+4,c+5=j+5$, or $c+5=j+6$.\\
if $c+5\leq j+4$, then $2^{c+3}<2^{j+3}$ we must have $2^{c+3}\leq 4$ or the digit difference is greater than 2, but $2^{c+3}\leq 4$ is impossible. 
%Then, $2^0+2^2+2^4+2^{a+1}+2^{b+1}=4+2^{j+3}+2^{j+4}+2^{j+6}$, which is impossible.\\
If $c+5=j+4$, then $2^{c+3}=2^{j+2}$ and $2^{c+1} = 4$, so the only possibility is $c=1, j=2,a=3,b=7$.\\
If $c+5=j+5$, then $2^{j+1}+2^{j+3}+2^{j+5}+2^{a+1}+2^{b+1}=4+2^{j+3}+2^{j+4}+2^{j+6}$ and the only possibility is $c=j=1,a=4,b=5.$\\
If $c+5=j+6$, then $2^{j+2}+2^{j+4}+2^{j+6}+2^{a+1}+2^{b+1}=4+2^{j+3}+2^{j+4}+2^{j+6}$, and the possibilities are $a=1,b=j+1,c=j+1$.\\

\[t(q_aq_b\cdot f_c)=2^{a+1}+2^{b+1}-2+2^{c+5}+2^{c+3}+2^{c+2}=2+2^{j+4}+2^{j+5}+2^{j+7}\]
%$2^{c+5}>4$, so we must have $c+5=j+2,c+5=j+3,c+5=j+4,c+5=j+5$, or $c+5=j+6$.\\
if $c+5\leq j+4$, then $2^{c+3}<2^{j+3}$ and $2^{c+3}=4$, which is impossible %Then, $2^1+2^2+2^4+2^{a+1}+2^{b+1}=4+2^{j+3}+2^{j+4}+2^{j+6}$, which is impossible.\\
If $c+5=j+5$, then $2^{c+2}<2^{j+3}$ and $2^{c+2}=4$, so the only possibility is $c=j=0$, which is not in the range of consideration.\\
If $c+5=j+6$, then $2^{j+3}+2^{j+4}+2^{j+6}+2^{a+1}+2^{b+1}=4+2^{j+4}+2^{j+5}+2^{j+7}$, which is impossible. %and the only possibility is $c=0,j=1,a=3,b=4.$\\
If $c+5=j+7$, then $2^{j+4}+2^{j+5}+2^{j+7}+2^{a+1}+2^{b+1}=4+2^{j+4}+2^{j+5}+2^{j+7}$, and the possibilities are $a=b=0,c=j+2$.

\[t(q_aq_b\cdot g_c)=2^{a+1}+2^{b+1}-2+2^{c+4}+2^{c+3}=2+2^{j+4}+2^{j+5}+2^{j+7}\]
After moving $-2$ to right hand side, we find both sides having 4 digits in the binary expansion, so the digits must match up exactly. If $c+4\leq j+4$, then $2^{c+3}=4$ which is impossible\\ %and the only possibility is $c=j=0,a=3,b=5$.\\
If $c+4=j+5$, then the equation gives $a=1,b=j+6,c=j+1$.\\
If $c+4=j+7$, then the equation $2^{a+1}+2^{b+1}+2^{j+6}+2^{j+7}=4+2^{j+4}+2^{j+5}+2^{j+7}$ cannot be satisfied.\\

\[t(q_aq_b\cdot 
p_c)=2^{a+1}+2^{b+1}-2+2^{c+6}+2^{c+3}+2^{c+2}=2+2^{j+4}+2^{j+5}+2^{j+7}\]
If $c+6 \leq j+5$, then $2^{c+3}<2^{j+3}$ and $2^{c+3}=4$, and the equation cannot be satisfied.\\
If $c+6=j+6$, then $2^{c+2}<2^{j+3}$ and $2^{c+2}=4$, so $c=j=0$, which is not in the range of consideration.
If $c+6=j+7$, then $2^{j+3}+2^{j+4}+2^{j+7}+2^{a+1}+2^{b+1}=4+2^{j+4}+2^{j+5}+2^{j+7}$, which cannot be satisfied.\\

\[t(q_aq_b\cdot D_3(c))=2^{a+1}+2^{b+1}-2+2^{c+7}+2^{j+1}=2+2^{j+4}+2^{j+5}+2^{j+7}\]
We must have $c+7\leq j+7$, then $2^{c+1}<2^{j+4}$ and $2^{c+1}=4$. The only possibilities are $c=1, j=3, a=6,b=9$, $c=1,j=4,a=8,b=10$, and $c=j=1,a=4,b=5$.\\

\[t(q_aq_b\cdot p'_c)=2^{a+1}+2^{b+1}-2+2^{c+7}+2^{c+4}+2^{c+1}=2+2^{j+4}+2^{j+5}+2^{j+7}\]
If $c+7 \leq j+6$, then $2^{c+4}<2^{j+4}$ and $2^{c+4}\le 4$, which is impossible.\\ 
We must have $c+7 = j+7$ and the only possibilities are $c=j=1, a=b=4$ (excluding $j=0$). %and $c=j=0, a=0,b=4$.\\

In conclusion, the only possible classes killing $h_0h_{j+4}Q_3(j)$ for $j\geq 1$ in this case are
\begin{enumerate}
    \item $d_1(q_3q_6\cdot d_0)=q_2q_6\cdot h_2d_0$ %j=1;h_2d_0=h_0e_0
    %\item $d_1(q_5q_7\cdot d_0)=q_5q_6\cdot h_6d_0$ %j=2;
    \item $d_1(q_1q_{j+4}\cdot d_{j+2})=q_0q_{j+4}\cdot h_0d_{j+2}+q_1q_{j+3}\cdot h_{j+3}d_{j+2}$ %j>=1;h_0d_2=h_3D_3(0)
    \item $d_1(q_3q_7\cdot e_1)=q_2q_7\cdot h_2e_1+q_3q_6\cdot h_6e_1$ %j=2; h_2e_1=h_1f_1=h_0^2c_2
    \item $d_1(q_4q_5\cdot e_1)=q_3q_5\cdot h_3e_1+q_3q_4\cdot h_4e_1$ %j=1
    \item $d_1(q_1q_{j+1}\cdot e_{j+1})=q_0q_{j+1}\cdot h_0e_{j+1}$ %j>=1
    %\item $d_1(q_2q_6\cdot f_0)=q_2q_5\cdot h_5f_0$ %j=1
    \item $d_1(q_1q_{j+6}\cdot g_{j+1})=q_0q_{j+6}\cdot h_0g_{j+1}+q_1q_{j+5}\cdot h_{j+5}g_{j+1}$ %j>=1; h_0g_1=h_2e_0;
    \item $d_1(q_6b_9\cdot D_3(1))=q_5q_9\cdot h_5D_3(1)$ %j=4
    \item $d_1(q_8q_{10}\cdot D_3(1))=q_8q_9\cdot h_9D_3(1)$ %j=5
    \item $d_1(q_4q_5\cdot D_3(1))=q_3q_5\cdot h_3D_3(1)+q_4^2\cdot h_4D_3(1)$ %j=2
    \item $d_4(q_4^2\cdot p_1')=q_2^2\cdot \langle h_3,h_4,p_1' \rangle = q_2^2 \cdot h_1e_3$ (Bruner's resolver (\cite{BrunerExt}) gives $\langle h_3,h_4,p_1' \rangle=h_1e_3$) %d_2=0, look for d_4 or higher
    %$d_8(q_4^2\cdot p_1')=q_0^2\cdot \langle h_1,h_2,h_3,h_4,p'_2 \rangle = q_0^2\cdot V_1$ %j=2
    %\item $q_0q_4\cdot p_0'$ %d_1=0, look for d_2 or higher
    %$d_4(q_0q_4\cdot p_0')=q_0^2\cdot \langle h_0,h_1,h_2,h_3,p'_1 \rangle = q_0^2\cdot V_0$ %j=1
    %$d_2(q_0q_4\cdot p_1')=q_0q_2\cdot \langle h_2,h_3,p'_1 \rangle$
\end{enumerate}

The class $q_0^2\cdot f_{j+2}$ survives to CESS, but it cannot support CESS differentials as products of permanent cycles $q_0$ and $f_{j+2}$ detecting $h_0$ and $f_{j+3}$ on Adams $E_2$ page.

\begin{figure}[H]
\centering
\begingroup
\setlength{\tabcolsep}{4pt}
\renewcommand{\arraystretch}{1.15}
\begin{tabular}{|wc{1.2cm}||wc{3.1cm}|p{4.6cm}|wc{4.6cm}|} \hline
      $w$ & $s=3$ & \multicolumn{1}{p{4.6cm}|}{\centering $s=4$} & $s=5$ \\\hline\hline
% w = 2j+10
    $2j+10$ & \tikzmark{part1-figure_f_d3_1-n49L}$\; q_{j+4}q_{j+6}\cdot h_1h_{j+2}^2$\tikzmark{part1-figure_f_d3_1-n49R} & \tikzmark{part1-figure_f_d3_1-n1L}$\; q_{j+4}q_{j+6}\cdot h_0^2h_{j+2}^2$\tikzmark{part1-figure_f_d3_1-n1R} &  \\
     & \tikzmark{part1-figure_f_d3_1-n51L}$\; q_{j+4}q_{j+6}\cdot h_0^2h_{j+3}$\tikzmark{part1-figure_f_d3_1-n51R} & \tikzmark{part1-figure_f_d3_1-n3L}$\; q_{j+5}^2\cdot h_1h_{j+3}^3$\tikzmark{part1-figure_f_d3_1-n3R} &  \\\hline
% w = 2j+9
    $2j+9$ & \tikzmark{part1-figure_f_d3_1-n42L}$\; q_{j+3}q_{j+6}\cdot h_0^2h_{j+4}$\tikzmark{part1-figure_f_d3_1-n42R} & \tikzmark{part1-figure_f_d3_1-n50L}$q_{j+4}q_{j+5}\cdot h_1h_{j+2}^2h_{j+5}$\tikzmark{part1-figure_f_d3_1-n50R} & \tikzmark{part1-figure_f_d3_1-n2L}$q_{j+4}q_{j+5}\cdot h_0^2h_{j+2}^2h_{j+5}$\tikzmark{part1-figure_f_d3_1-n2R} \\
     & \tikzmark{part1-figure_f_d3_1-n10L}$\; q_{j+3}q_{j+6}\cdot h_1h_{j+3}^2$\tikzmark{part1-figure_f_d3_1-n10R} & \tikzmark{part1-figure_f_d3_1-n5L}$\; q_{j+3}q_{j+6}\cdot h_0^2h_{j+3}^2$\tikzmark{part1-figure_f_d3_1-n5R} &  \\
     &  & \tikzmark{part1-figure_f_d3_1-n7L}$\; q_{j+4}q_{j+5}\cdot h_0^2h_{j+3}h_{j+5}$\tikzmark{part1-figure_f_d3_1-n7R} &  \\\hline
% w = 2j+8
    $2j+8$ &  & \tikzmark{part1-figure_f_d3_1-n8L}$q_{j+2}q_{j+6}\cdot h_0^2h_{j+2}h_{j+4}$\tikzmark{part1-figure_f_d3_1-n8R} & \tikzmark{part1-figure_f_d3_1-n6L}$q_{j+3}q_{j+5}\cdot h_0^2h_{j+4}^3$\tikzmark{part1-figure_f_d3_1-n6R} \\
     &  & \tikzmark{part1-figure_f_d3_1-n11L}$q_{j+3}q_{j+5}\cdot h_1h_{j+3}^2h_{j+5}$\tikzmark{part1-figure_f_d3_1-n11R} &  \\
     &  & \tikzmark{part1-figure_f_d3_1-n47L}$\; q_{j+2}q_{j+6}\cdot h_1h_{j+1}^2h_{j+4}$\tikzmark{part1-figure_f_d3_1-n47R} &  \\\hline
% w = 2j+7
    $2j+7$ & \tikzmark{part1-figure_f_d3_1-n32L}$\; q_{j+3}q_{j+4}\cdot h_0^2h_{j+6}$\tikzmark{part1-figure_f_d3_1-n32R} & \tikzmark{part1-figure_f_d3_1-n12L}$\; q_{j+3}q_{j+4}\cdot h_0^2h_{j+5}^2$\tikzmark{part1-figure_f_d3_1-n12R} &  \\
     & \tikzmark{part1-figure_f_d3_1-n15L}$\; q_{j+3}q_{j+4}\cdot h_1h_{j+5}^2$\tikzmark{part1-figure_f_d3_1-n15R} & \tikzmark{part1-figure_f_d3_1-n45L}$\; q_{j+1}q_{j+6}\cdot h_1c_{j+1}$\tikzmark{part1-figure_f_d3_1-n45R} & \tikzmark{part1-figure_f_d3_1-n48L}$q_{j+1}q_{j+6}\cdot h_1h_{j+1}^3h_{j+4}$\tikzmark{part1-figure_f_d3_1-n48R} \\\hline
% w = 2j+6
    $2j+6$ & \tikzmark{part1-figure_f_d3_1-n22L}$q_{j+3}^2\cdot h_1h_{j+3}h_{j+6}$\tikzmark{part1-figure_f_d3_1-n22R} & \tikzmark{part1-figure_f_d3_1-n13L}$\; q_{j+2}q_{j+4}\cdot h_0^2h_{j+2}h_{j+6}$\tikzmark{part1-figure_f_d3_1-n13R} & \tikzmark{part1-figure_f_d3_1-n9L}$q_{j+2}q_{j+4}\cdot h_0^2h_{j+2}h_{j+5}^2$\tikzmark{part1-figure_f_d3_1-n9R} \\
     &  & \tikzmark{part1-figure_f_d3_1-n17L}$\; q_{j+3}^2\cdot h_0^2h_{j+3}h_{j+6}$\tikzmark{part1-figure_f_d3_1-n17R} & \tikzmark{part1-figure_f_d3_1-n46L}$q_{j+1}q_{j+5}\cdot h_1h_{j+5}c_{j+1}$\tikzmark{part1-figure_f_d3_1-n46R} \\
     &  & \tikzmark{part1-figure_f_d3_1-n18L}$\; q_{j+3}^2\cdot h_1h_{j+2}^2h_{j+6}$\tikzmark{part1-figure_f_d3_1-n18R} &  \\
     &  & \tikzmark{part1-figure_f_d3_1-n16L}$q_{j+2}q_{j+4}\cdot h_1h_{j+2}h_{j+5}^2$\tikzmark{part1-figure_f_d3_1-n16R} &  \\
     &  & \tikzmark{part1-figure_f_d3_1-n52L}$\; q_{j+2}q_{j+4}\cdot h_1h_{j+1}^2h_{j+5}$\tikzmark{part1-figure_f_d3_1-n52R} &  \\\hline
% w = 2j+5
    $2j+5$ &  &  & \tikzmark{part1-figure_f_d3_1-n53L}$q_{j+1}q_{j+4}\cdot h_1h_{j+1}^3h_{j+5}$\tikzmark{part1-figure_f_d3_1-n53R} \\
     &  &  & \tikzmark{part1-figure_f_d3_1-n54L}$q_{j+2}q_{j+3}\cdot h_1h_{j+1}^2h_{j+3}h_{j+5}$\tikzmark{part1-figure_f_d3_1-n54R} \\\hline
% w = 2j+4
    $2j+4$ & \tikzmark{part1-figure_f_d3_1-n57L}$\; q_{j+2}^2\cdot h_1h_{j+4}h_{j+6}$\tikzmark{part1-figure_f_d3_1-n57R} & \tikzmark{part1-figure_f_d3_1-n20L}$\; q_{j+2}^2\cdot h_0^2h_{j+4}h_{j+6}$\tikzmark{part1-figure_f_d3_1-n20R} & \tikzmark{part1-figure_f_d3_1-n4L}$q_{j+2}^2\cdot h_1h_{j+3}c_{j+2}$\tikzmark{part1-figure_f_d3_1-n4R} \\
     &  & \tikzmark{part1-figure_f_d3_1-n23L}$\; q_{j+2}^2\cdot h_1h_{j+3}^2h_{j+6}$\tikzmark{part1-figure_f_d3_1-n23R} & \tikzmark{part1-figure_f_d3_1-n14L}$q_{j+2}^2\cdot h_0^2h_{j+3}^2h_{j+6}$\tikzmark{part1-figure_f_d3_1-n14R} \\\hline
% w = 2j+3
    $2j+3$ &  & \tikzmark{part1-figure_f_d3_1-n55L}$\; q_{j+1}q_{j+2}\cdot h_1h_{j+1}h_{j+4}h_{j+6}$\tikzmark{part1-figure_f_d3_1-n55R} &  \\\hline
% w = 2j+2
    $2j+2$ &  & \tikzmark{part1-figure_f_d3_1-n58L}$q_{j+1}^2\cdot h_1h_{j+2}h_{j+4}h_{j+6}$\tikzmark{part1-figure_f_d3_1-n58R} & \tikzmark{part1-figure_f_d3_1-n21L}$q_{j+1}^2\cdot h_0^2h_{j+2}h_{j+4}h_{j+6}$\tikzmark{part1-figure_f_d3_1-n21R} \\
     &  &  & \tikzmark{part1-figure_f_d3_1-n56L}$q_{j+1}^2\cdot h_1h_{j+1}^2h_{j+4}h_{j+6}$\tikzmark{part1-figure_f_d3_1-n56R} \\\hline
% w = 2j+1
    $2j+1$ &  &  &  \\\hline
% w = 2j
    $2j$ &  &  & \tikzmark{part1-figure_f_d3_1-n19L}$q_j^2\cdot h_1h_{j+6}c_{j+1}$\tikzmark{part1-figure_f_d3_1-n19R} \\\hline
    \multicolumn{4}{|c|}{} \\\hline
\end{tabular}
\begin{tikzpicture}[overlay, remember picture, shorten >=.5pt, shorten <=.5pt, transform canvas={yshift=.25\baselineskip}]
    \draw [->] ({pic cs:part1-figure_f_d3_1-n1R}) [red] to ({pic cs:part1-figure_f_d3_1-n2L});
    \draw [->] ({pic cs:part1-figure_f_d3_1-n3R}) [green!60!black,out=15,in=180] to ({pic cs:part1-figure_f_d3_1-n4L});
    \draw [->] ({pic cs:part1-figure_f_d3_1-n5R}) [red,bend left=5] to ({pic cs:part1-figure_f_d3_1-n6L});
    \draw [->] ({pic cs:part1-figure_f_d3_1-n7R}) [red,bend right=5] to ({pic cs:part1-figure_f_d3_1-n6L});
    \draw [->] ({pic cs:part1-figure_f_d3_1-n42R}) [red] to ({pic cs:part1-figure_f_d3_1-n8L});
    \draw [->] ({pic cs:part1-figure_f_d3_1-n10R}) [red] to ({pic cs:part1-figure_f_d3_1-n11L});
    \draw [->] ({pic cs:part1-figure_f_d3_1-n12R}) [red,bend right=5] to ({pic cs:part1-figure_f_d3_1-n9L});
    \draw [->] ({pic cs:part1-figure_f_d3_1-n13R}) [blue,out=0,in=150] to ({pic cs:part1-figure_f_d3_1-n14L});
    \draw [->] ({pic cs:part1-figure_f_d3_1-n15R}) [red] to ({pic cs:part1-figure_f_d3_1-n16L});
    \draw [->] ({pic cs:part1-figure_f_d3_1-n17R}) [blue,bend right=5] to ({pic cs:part1-figure_f_d3_1-n14L});
    \draw [->] ({pic cs:part1-figure_f_d3_1-n18R}) [green!60!black,out=0,in=160] to ({pic cs:part1-figure_f_d3_1-n19L});
    \draw [->] ({pic cs:part1-figure_f_d3_1-n20R}) [blue] to ({pic cs:part1-figure_f_d3_1-n21L});
    \draw [->] ({pic cs:part1-figure_f_d3_1-n22R}) [blue] to ({pic cs:part1-figure_f_d3_1-n23L});
    \draw [->] ({pic cs:part1-figure_f_d3_1-n32R}) [red,bend left=12] to ({pic cs:part1-figure_f_d3_1-n13L});
    \draw [->] ({pic cs:part1-figure_f_d3_1-n32R}) [red,bend right=12] to ({pic cs:part1-figure_f_d3_1-n17L});
    \draw [->] ({pic cs:part1-figure_f_d3_1-n45R}) [red,out=0,in=180] to ({pic cs:part1-figure_f_d3_1-n46L});
    \draw [->] ({pic cs:part1-figure_f_d3_1-n47R}) [red] to ({pic cs:part1-figure_f_d3_1-n48L});
    \draw [->] ({pic cs:part1-figure_f_d3_1-n49R}) [red] to ({pic cs:part1-figure_f_d3_1-n50L});
    \draw [->] ({pic cs:part1-figure_f_d3_1-n51R}) [red,bend left=12] to ({pic cs:part1-figure_f_d3_1-n5L});
    \draw [->] ({pic cs:part1-figure_f_d3_1-n51R}) [red,bend right=12] to ({pic cs:part1-figure_f_d3_1-n7L});
    \draw [->] ({pic cs:part1-figure_f_d3_1-n52R}) [red,bend left=8] to ({pic cs:part1-figure_f_d3_1-n53L});
    \draw [->] ({pic cs:part1-figure_f_d3_1-n52R}) [red,bend right=8] to ({pic cs:part1-figure_f_d3_1-n54L});
    \draw [->] ({pic cs:part1-figure_f_d3_1-n55R}) [red] to ({pic cs:part1-figure_f_d3_1-n56L});
    \draw [->] ({pic cs:part1-figure_f_d3_1-n57R}) [blue] to ({pic cs:part1-figure_f_d3_1-n58L});
\end{tikzpicture}
\par\endgroup
\end{figure}

\begin{figure}[H]
\centering
\begingroup
\setlength{\tabcolsep}{4pt}
\renewcommand{\arraystretch}{1.15}
\begin{tabular}{|wc{1.2cm}||wc{3.1cm}|p{4.6cm}|wc{4.6cm}|} \hline
      $w$ & $s=3$ & \multicolumn{1}{p{4.6cm}|}{\centering $s=4$} & $s=5$ \\\hline\hline
% w = j+7
    $j+7$ & \tikzmark{part2-figure_f_d3_1-n59L}$\; q_1q_{j+6}\cdot h_{j+3}^3$\tikzmark{part2-figure_f_d3_1-n59R} & \tikzmark{part2-figure_f_d3_1-n39L}$\; q_1q_{j+6}\cdot g_{j+1}$\tikzmark{part2-figure_f_d3_1-n39R} &  \\\hline
% w = j+6
    $j+6$ &  & \tikzmark{part2-figure_f_d3_1-n24L}$q_0q_{j+6}\cdot h_0h_{j+3}^3$\tikzmark{part2-figure_f_d3_1-n24R} & \tikzmark{part2-figure_f_d3_1-n40L}$q_0q_{j+6}\cdot h_0g_{j+1}$\tikzmark{part2-figure_f_d3_1-n40R} \\
     &  &  & \tikzmark{part2-figure_f_d3_1-n41L}$q_1q_{j+5}\cdot h_{j+5}g_{j+1}$\tikzmark{part2-figure_f_d3_1-n41R} \\\hline
% w = j+5
    $j+5$ & \tikzmark{part2-figure_f_d3_1-n60L}$\; q_1q_{j+4}\cdot h_{j+2}^2h_{j+6}$\tikzmark{part2-figure_f_d3_1-n60R} & \tikzmark{part2-figure_f_d3_1-n33L}$\; q_1q_{j+4}\cdot d_{j+2}$\tikzmark{part2-figure_f_d3_1-n33R} &  \\\hline
% w = j+4
    $j+4$ & \tikzmark{part2-figure_f_d3_1-n28L}$q_0q_{j+4}\cdot h_0h_{j+3}h_{j+6}$\tikzmark{part2-figure_f_d3_1-n28R} & \tikzmark{part2-figure_f_d3_1-n26L}$q_0q_{j+4}\cdot h_0h_{j+2}^2h_{j+6}$\tikzmark{part2-figure_f_d3_1-n26R} &  \\
     &  &  & \tikzmark{part2-figure_f_d3_1-n34L}$q_0q_{j+4}\cdot h_0d_{j+2}$\tikzmark{part2-figure_f_d3_1-n34R} \\
     &  &  & \tikzmark{part2-figure_f_d3_1-n35L}$q_1q_{j+3}\cdot h_{j+3}d_{j+2}$\tikzmark{part2-figure_f_d3_1-n35R} \\\hline
% w = j+3
    $j+3$ & \tikzmark{part2-figure_f_d3_1-n30L}$q_0q_{j+3}\cdot h_0h_{j+4}h_{j+6}$\tikzmark{part2-figure_f_d3_1-n30R} & \tikzmark{part2-figure_f_d3_1-n29L}$\; q_0q_{j+3}\cdot h_0h_{j+3}^2h_{j+6}$\tikzmark{part2-figure_f_d3_1-n29R} &  \\
     &  & \tikzmark{part2-figure_f_d3_1-n36L}$\; q_1q_{j+2}\cdot e_{j+2}$\tikzmark{part2-figure_f_d3_1-n36R} &  \\
     &  & \tikzmark{part2-figure_f_d3_1-n61L}$\; q_1q_{j+2}\cdot h_{j+1}^2h_{j+4}h_{j+6}$\tikzmark{part2-figure_f_d3_1-n61R} &  \\\hline
% w = j+2
    $j+2$ &  & \tikzmark{part2-figure_f_d3_1-n31L}$\; q_0q_{j+2}\cdot h_0h_{j+2}h_{j+4}h_{j+6}$\tikzmark{part2-figure_f_d3_1-n31R} & \tikzmark{part2-figure_f_d3_1-n37L}$q_0q_{j+2}\cdot h_0e_{j+2}$\tikzmark{part2-figure_f_d3_1-n37R} \\
     &  & \tikzmark{part2-figure_f_d3_1-n43L}$\; q_1q_{j+1}\cdot h_{j+6}c_{j+1}$\tikzmark{part2-figure_f_d3_1-n43R} & \tikzmark{part2-figure_f_d3_1-n62L}$q_0q_{j+2}\cdot h_0h_{j+1}^2h_{j+4}h_{j+6}$\tikzmark{part2-figure_f_d3_1-n62R} \\
     &  &  & \tikzmark{part2-figure_f_d3_1-n63L}$q_1q_{j+1}\cdot h_{j+1}^3h_{j+4}h_{j+6}$\tikzmark{part2-figure_f_d3_1-n63R} \\\hline
% w = j+1
    $j+1$ &  &  & \tikzmark{part2-figure_f_d3_1-n44L}$q_0q_{j+1}\cdot h_0h_{j+6}c_{j+1}$\tikzmark{part2-figure_f_d3_1-n44R} \\\hline
    \end{tabular}
\begin{tikzpicture}[overlay, remember picture, shorten >=.5pt, shorten <=.5pt, transform canvas={yshift=.25\baselineskip}]
    \draw [->] ({pic cs:part2-figure_f_d3_1-n59R}) [red] to ({pic cs:part2-figure_f_d3_1-n24L});
    \draw [->] ({pic cs:part2-figure_f_d3_1-n60R}) [red] to ({pic cs:part2-figure_f_d3_1-n26L});
    \draw [->] ({pic cs:part2-figure_f_d3_1-n28R}) [red] to ({pic cs:part2-figure_f_d3_1-n29L});
    \draw [->] ({pic cs:part2-figure_f_d3_1-n30R}) [red] to ({pic cs:part2-figure_f_d3_1-n31L});
    \draw [->] ({pic cs:part2-figure_f_d3_1-n33R}) [red,bend left=12] to ({pic cs:part2-figure_f_d3_1-n34L});
    \draw [->] ({pic cs:part2-figure_f_d3_1-n33R}) [red,out=0,in=180] to ({pic cs:part2-figure_f_d3_1-n35L});
    \draw [->] ({pic cs:part2-figure_f_d3_1-n36R}) [red,bend left=8] to ({pic cs:part2-figure_f_d3_1-n37L});
    \draw [->] ({pic cs:part2-figure_f_d3_1-n39R}) [red,bend left=8] to ({pic cs:part2-figure_f_d3_1-n40L});
    \draw [->] ({pic cs:part2-figure_f_d3_1-n39R}) [red,out=0,in=180] to ({pic cs:part2-figure_f_d3_1-n41L});
    \draw [->] ({pic cs:part2-figure_f_d3_1-n43R}) [red] to ({pic cs:part2-figure_f_d3_1-n44L});
    \draw [->] ({pic cs:part2-figure_f_d3_1-n61R}) [red,bend left=8] to ({pic cs:part2-figure_f_d3_1-n62L});
    \draw [->] ({pic cs:part2-figure_f_d3_1-n61R}) [red,out=0,in=180] to ({pic cs:part2-figure_f_d3_1-n63L});
\end{tikzpicture}
\par\endgroup
\caption{The BXSS in degree $(s,k,t) = (4,2,2^{j+4}\cdot 11+2)$ for \cref{prop:no-d3}, where $j\geq 2$.
The vertical axis is the BX filtration $w$, and the horizontal axis is the Adams filtration $s$.
Red, blue, and green arrows indicate $d_1$, $d_2$, and $d_6$, respectively. Rows are grouped into symbolic $w$-grading bands.}
\label{fig:f_d3_1}
\label{fig:f_d3_1-part2}
\end{figure}

\end{proof}

\begin{prop}\label{prop:no-d5}
The class $h_0h_{j+4}Q_3(j)$ does not receive a Cartan-Eilenberg $d_5$ differential for any $j\geq1$.
\end{prop}
\begin{proof}
If $h_0h_{j+4}Q_3(j)$ receives a Cartan-Eilenberg $d_5$ differential
\[
d_5(x) = h_0h_{j+4}Q_3(j),
\]
then $x$ must lie in Cartan-Eilenberg filtration $k = 4$. The Ext component of $x$ must have filtration $s=2$ since $h_0h_{j+4}Q_3(j)$ lies in total filtration $k+s = 6+1$. A general class in $(k,s) = (4,2)$ takes the following form:
	\[
		q_aq_bq_cq_d\cdot h_eh_f
	\]

Since Cartan-Eilenberg differentials preserve internal degree, we have
	\[
		\begin{split}
			t(h_0h_{j+4}Q_3(j)) &= 2^1+2^{j+4}+2^{j+5}+2^{j+7} \\
			&= 2^{a+1}+2^{b+1}+2^{c+1}+2^{d+1}-4+2^{e+1}+2^{f+1} \\
			&= t(q_aq_bq_cq_d\cdot h_eh_f).
		\end{split}
	\]
	Adding 4 to both sides of the middle equality, and then dividing by 2 yields
	\[
		2^0+2^1+2^{j+3}+2^{j+4}+2^{j+6} = 2^a+2^b+2^c+2^d+2^e+2^f.
	\]
	Note that $2^0+2^1$ and $2^0+2^0+2^0$ are the only ways of writing $2^0+2^1$ as a sum of powers of 2. Then there are four possibilities for the indices $a,b,c,d,e,f$, up to reordering of $a,b,c,d$ and $e,f$:
	\begin{enumerate}
		\item $[a,b,c,d,e,f] = [0,0,0,j+3,j+4,j+6]$
		\item $[a,b,c,d,e,f] = [0,1,j+2,j+2,j+4,j+6]$
		\item $[a,b,c,d,e,f] = [0,1,j+3,j+3,j+3,j+6]$
		\item $[a,b,c,d,e,f] = [0,1,j+3,j+4,j+5,j+5]$
	\end{enumerate}

	In case (i), we obtain the following possibilities by brute-force calculation:
	\begin{enumerate}
		\item $q_0^3q_{j+3}\cdot h_{j+4}h_{j+6}$
		\item $q_0^3q_{j+4}\cdot h_{j+3}h_{j+6}$
		\item $q_0^3q_{j+6}\cdot h_{j+3}h_{j+4} = 0$
		\item $q_0^2q_{j+3}q_{j+4}\cdot h_0h_{j+6}$
		\item $q_0^2q_{j+3}q_{j+6}\cdot h_0h_{j+4}$
		\item $q_0^2q_{j+4}q_{j+6}\cdot h_0h_{j+3}$
		\item $q_0q_{j+3}q_{j+4}q_{j+6}\cdot h_0^2$
	\end{enumerate}
	
    Removing (3) since it is trivial, the remaining possibilities all support Burklund-Xu $d_1$ differentials, so they do not survive to the Cartan-Eilenberg $E_2$ page.\newline

	In case (ii), we obtain the following possibilities:
	\begin{enumerate}
		\item $q_0q_1q_{j+2}^2\cdot h_{j+4}h_{j+6}$
		\item $q_0q_1q_{j+2}q_{j+4}\cdot h_{j+2}h_{j+6}$
		\item $q_0q_1q_{j+2}q_{j+6}\cdot h_{j+2}h_{j+4}$
		\item $q_0q_1q_{j+4}q_{j+6}\cdot h_{j+2}^2$
		\item $q_0q_{j+2}^2q_{j+4}\cdot h_1h_{j+6}$
		\item $q_0q_{j+2}^2q_{j+6}\cdot h_1h_{j+4}$
		\item $q_0q_{j+2}q_{j+4}q_{j+6}\cdot h_1h_{j+2}$
		\item $q_1q_{j+2}^2q_{j+4}\cdot h_0h_{j+6}$
		\item $q_1q_{j+2}^2q_{j+6}\cdot h_0h_{j+4}$
		\item $q_1q_{j+2}q_{j+4}q_{j+6}\cdot h_0h_{j+2}$
		\item $q_{j+2}^2q_{j+4}q_{j+6}\cdot h_0h_1 = 0$
	\end{enumerate}
	As above, we remove (11) since it is trivial. All of the remaining possibilities support Burklund-Xu $d_1$ differentials except (6). However, this class does support a $d_2$ differential
    \[
        \begin{split}
            d_2(q_0q_{j+2}^2q_{j+6}\cdot h_1h_{j+4}) &= q_0q_{j+6}d_2(q_{j+2}^2\cdot h_1h_{j+4}) + q_0q_{j+2}^2d_2(q_{j+6}\cdot h_1h_{j+4}) \\
            &= q_0q_{j+6}q_{j+1}^2\cdot h_1h_{j+2}h_{j+4} + q_0q_{j+2}^2q_{j+4}\cdot\langle h_{j+4},h_{j+5},h_1h_{j+4}\rangle \\
            &= q_0q_{j+6}q_{j+1}^2\cdot h_1h_{j+2}h_{j+4} + q_0q_{j+2}^2q_{j+4}\cdot h_1h_{j+5}^2,
        \end{split}
    \]
    by \Cref{prop:leibniz}. These terms do not cancel on the Burklund-Xu $E_2$ page since they are not component-wise equal and the neither term is a summand in a $d_1$ differential. The latter claim follows from the fact that every term in a $d_1$ differential has at least one index, such as $j+1$, common between its $q$ and Ext components. Since neither term has common indices for any value of $j$, the claim follows. In particular, the target of this $d_2$ differential must be nontrivial for all $j\geq0$, so possibility (6) does not survive to the Cartan-Eilenberg $E_2$ page.

	In case (iii), we obtain the following possibilities:
	\begin{enumerate}
		\item $q_0q_1q_{j+3}^2\cdot h_{j+3}h_{j+6}$
		\item $q_0q_1q_{j+3}q_{j+6}\cdot h_{j+3}^2$
		\item $q_0q_{j+3}^3\cdot h_1h_{j+6}$
		\item $q_0q_{j+3}^2q_{j+6}\cdot h_1h_{j+3}$
		\item $q_1q_{j+3}^3\cdot h_0h_{j+6}$
		\item $q_1q_{j+3}^2q_{j+6}\cdot h_0h_{j+3}$
		\item $q_{j+3}^3q_{j+6}\cdot h_0h_1 = 0$
	\end{enumerate}

    Removing (7) since it is trivial, all of the remaining possibilities support Burklund-Xu $d_1$ differentials for all $j\geq0$ except (3), which does so only for $j\geq1$. If $j=0$, then (3) becomes
    \[
        q_0q_3^3\cdot h_1h_6,
    \]
    which does not support a nonzero $d_1$ differential. However, we can use the Leibniz rule to compute
    \[
        d_2(q_0q_3^3\cdot h_1h_6) = q_0q_1q_3^2\cdot\langle h_1,h_2,h_1h_6\rangle + q_0q_2^2q_3\cdot h_1h_3h_6,
    \]
    which we claim is nontrivial. Indeed, while the second term is the target of the differential
    \[
        d_1(q_0q_2^2q_4\cdot h_1h_6),
    \]
    the first term can be evaluated using knowledge of the Massey product
    \[
        \langle h_1,h_2,h_1h_6\rangle = h_2^2h_6,
    \]
    whose indeterminacy vanishes. We thus obtain
    \[
        d_2(q_0q_3^3\cdot h_1h_6) = q_0q_1q_3^2\cdot h_2^2h_6,
    \]
    which is a $d_1$ cycle since
    \[
        d_1(q_0q_1q_3^2\cdot h_2^2h_6) = q_0^2q_3^2\cdot h_0h_2^2h_6 = 0,
    \]
    by the relation $h_jh_{j+2}^2 = 0$ in Theorem 1.2 of \cite{Chen}. Furthermore, the target
    \[
        q_0q_1q_3^2\cdot h_2^2h_6
    \]
    cannot receive a $d_1$ differential since there is no matching index between the $q$ and Ext components. In other words, the intersection $\{0,1,3\}\cap\{2,6\}$ is empty. Then the $j=0$ element of (3) supports a nonzero Burklund-Xu $d_2$ differential.

	In case (iv), we obtain the following possibilities:
	\begin{enumerate}
		\item $q_0q_1q_{j+3}q_{j+4}\cdot h_{j+5}^2$
		\item $q_0q_1q_{j+3}q_{j+5}\cdot h_{j+4}h_{j+5} = 0$
		\item $q_0q_1q_{j+4}q_{j+5}\cdot h_{j+3}h_{j+5}$
		\item $q_0q_1q_{j+5}^2\cdot h_{j+3}h_{j+4} = 0$
		\item $q_0q_{j+3}q_{j+4}q_{j+5}\cdot h_1h_{j+5}$
		\item $q_0q_{j+3}q_{j+5}^2\cdot h_1h_{j+4}$
		\item $q_0q_{j+4}q_{j+5}^2\cdot h_1h_{j+3}$
		\item $q_1q_{j+3}q_{j+4}q_{j+5}\cdot h_0h_{j+5}$
		\item $q_1q_{j+3}q_{j+5}^2\cdot h_0h_{j+4}$
		\item $q_1q_{j+4}q_{j+5}^2\cdot h_0h_{j+3}$
		\item $q_{j+3}q_{j+4}q_{j+5}^2\cdot h_0h_1 = 0$
	\end{enumerate}

    Removing the trivial classes, all of the remaining possibilities support Burklund-Xu $d_1$ differentials for all $j\geq0$ except (6), which does so only for $j\geq1$. If $j = 0$, then (6) becomes
    \[
        q_0q_3q_5^2\cdot h_1h_4,
    \]
    which does not support a nonzero $d_1$ differential. However, we can use the Leibniz rule to compute
    \[
    \begin{split}
        d_2(q_0q_3q_5^2\cdot h_1h_4) &= q_0q_1q_5^2\cdot\langle h_1,h_2,h_1h_4\rangle + q_0q_3q_4^2\cdot h_1h_4h_5 \\
        &= q_0q_1q_5^2\cdot h_2^2h_4,
    \end{split}
    \]
    which we claim is nontrivial. Indeed, this class is a $d_1$ cycle since
    \[
        d_1(q_0q_1q_5^2\cdot h_2^2h_4) = q_0^2q_5^2\cdot h_0h_2^2h_4 = 0,
    \]
    by the relation $h_jh_{j+2}^2 = 0$ in Theorem 1.2 of \cite{Chen}. Furthermore, the target
    \[
        q_0q_1q_5^2\cdot h_2^2h_4 = q_0q_1q_5^2\cdot h_3^3
    \]
    cannot receive a $d_1$ differential since there is no matching index between the $q$ and Ext components. In other words, the intersection $\{0, 1, 5\}\cap\{2, 3, 4\}$ is empty. Then the $j=0$ element of (6) supports a nonzero Burklund-Xu $d_2$ differential.
\end{proof}

Then we have shown that $h_1h_{j+4}Q_3(j)$ is nonzero on the Adams $E_2$ page for $j\geq1$. In order to show that it detects a nontrivial class on the Adams $E_3$ page, we need to show that it is not the target of an Adams $d_2$ differential. To that end, we first determine all possible sources of such a differential:

\begin{prop}\label{prop:possible-d2}
The only Adams $E_2$ classes $x$ which can possibly satisfy
\[
d_2(x) = h_1h_{j+4}Q_3(j)
\]
are linear combinations of the classes
\[
\begin{array}{c c c}
h_0f_{j+2} &&h_0h_{j+3}^3h_{j+6}
\end{array}
\]
for $j\geq0$.
\end{prop}
\begin{proof}
In order to enumerate all such classes $x$, we first note that Adams $d_2$ differentials raise internal degree by 1 so that
\[
t(x) = t(h_1h_{j+4}Q_3(j))-1 = 2^1+2^{j+3}+2^{j+4}+2^{j+6}-1 = 2^0+2^{j+3}+2^{j+4}+2^{j+6} =: t_\text{source},
\]
where we are using the internal degree as elements in $\text{Ext}_\mathcal{A}$, as opposed to $\text{Ext}_\mathcal{P}$, as in the rest of the paper. The difference is immaterial, but it should be noted if the reader is carefully following degrees. We choose to use $\text{Ext}_\mathcal{A}$ since we are studying the Adams $E_2$ page directly, rather than using Burklund-Xu representatives.

We first consider classes $x$ which are indecomposable. Using the presentation in \cite{Chen}, we rule out all but $x = D_1(i)$ since the remaining indecomposables have internal degrees with fewer than four terms in their 2-adic expansions. In order for $x = D_1(i)$ to have internal degree $t_\text{source}$, we must have $i=0$. However, this would give
\[
2^0+2^3+2^4+2^5 = 2^0+2^{j+3}+2^{j+4}+2^{j+6},
\]
which is not possible since matching the highest powers would force $j = -1$. Then $x$ cannot be an indecomposable element.

We now consider $x = h_a\cdot y$, where $y$ is an indecomposable in Adams filtration 4. Using the presentation in \cite{Chen}, we rule out all but $y = e_b$, $y = f_b$, $y = p_b$, and $y = p_b'$ since the remaining indecomposables have fewer than three terms in their 2-adic expansions.

If $x = h_a\cdot e_b$, then
\[
t(h_a\cdot e_b) = 2^a+2^b+2^{b+2}+2^{b+4} = 2^0+2^{j+3}+2^{j+4}+2^{j+6} = t_\text{source},
\]
meaning that $a=0$ or $b=0$. If $a=0$, then $b=j+3$, giving
\[
2^0+2^{j+3}+2^{j+5}+2^{j+7} = 2^0+2^{j+3}+2^{j+4}+2^{j+6},
\]
which is not possible for any value of $j$. If $b=0$, then
\[
2^a+2^0+2^2+2^4 = 2^0+2^{j+3}+2^{j+4}+2^{j+6},
\]
which is not possible for any values of $a$ and $j\geq0$. Then $x \neq h_a\cdot e_b$.

If $x = h_a\cdot f_b$, then
\[
t(h_a\cdot f_b) = 2^a+2^{b+1}+2^{b+2}+2^{b+4} = 2^0+2^{j+3}+2^{j+4}+2^{j+6} = t_\text{source},
\]
which is possible only if $a=0$ and $b=j+2$ so that $x = h_0f_{j+2}$. Classes of this form have the proper internal degree, so we must consider them.

If $x = h_a\cdot p_b$, then
\[
t(h_a\cdot p_b) = 2^a+2^b+2^{b+2}+2^{b+5} = 2^0+2^{j+3}+2^{j+4}+2^{j+6} = t_\text{source},
\]
which is not possible because $2^b$ and $2^{b+2}$ are two factors of 2 apart on the left, whereas no two powers on the right can be two factors of 2 apart without violating $j\geq0$. Then $x \neq h_a\cdot p_b$.

If $x = h_a\cdot p_b'$, then
\[
t(h_a\cdot p_b') = 2^a+2^b+2^{b+3}+2^{b+6} = 2^0+2^{j+3}+2^{j+4}+2^{j+6} = t_\text{source},
\]
meaning that $a=0$ or $b=0$. If $a=0$, then
\[
2^0+2^b+2^{b+3}+2^{b+6} = 2^0+2^{j+3}+2^{j+4}+2^{j+6},
\]
which is not possible because, after subtracting 1 from both sides, the smallest powers of 2 on the left are three factors of 2 apart while the smallest powers of 2 on the right are only one factor apart. Then $b=0$ so that
\[
2^a+2^0+2^3+2^6 = 2^0+2^{j+3}+2^{j+4}+2^{j+6},
\]
which is possible only if $j=0$ and $a=j+4=4$. In this case, we have
\[
x = h_4\cdot p_0' = h_0f_2
\]
by the presentation in \cite{Chen}. This element belongs to the family already uncovered above.

We now consider elements of the form
\[
x = h_\ell h_mc_n,
\]
where we take $\ell\leq m$ without loss of generality, by commutativity. In this case, we have
\[
t(x) = 2^\ell+2^m+2^n+2^{n+1}+2^{n+3} = 2^0+2^{j+3}+2^{j+4}+2^{j+6},
\]
meaning that $\ell = 0$ or $n = 0$. If $\ell = 0$, then
\[
2^0+2^m+2^n+2^{n+1}+2^{n+3} = 2^0+2^{j+3}+2^{j+4}+2^{j+6},
\]
which is possible only if $m\in\{n,n+1,n+3\}$. If $m = n$ or $m = n+1$, then the smallest powers of 2 on the left would differ by more than one factor of 2, %This should be more detailed
preventing the left-hand side from matching the right-hand side. If $m = n+3$, then
\[
2^0+2^n+2^{n+1}+2^{n+4} = 2^0+2^{j+3}+2^{j+4}+2^{j+6}
\]
forces $n = j+2$ so that
\[
x = h_0h_{j+5}c_{j+2},
\]
but this element is identically 0, by the presentation in \cite{Chen}. Then we conclude that $\ell\neq0$.

If $n = 0$, then
\[
2^\ell+2^m+2^0+2^1+2^3 = 2^0+2^{j+3}+2^{j+4}+2^{j+6},
\]
where, by our conclusion above that $\ell\neq0$, we may take $1\leq\ell\leq m$. Then we may subtract 1 from both sides and then divide by 2 to obtain
\[
2^{\ell-1}+2^{m-1}+2^0+2^2 = 2^{j+2}+2^{j+3}+2^{j+5}.
\]
Since the lowest possible value of $j+2$ is $0+2=2$, the ordering $1\leq\ell\leq m$ forces $\ell-1 = 0$ and $m-1 = 1$ so that the powers of 2 on the left-hand side combine to form $2^{j+2}$. In this case,
\[
 2^2+2^2 = 2^0+2^1+2^0+2^2 = 2^{\ell-1}+2^{m-1}+2^0+2^2 = 2^{0+2}+2^{0+3}+2^{0+5},
\]
a clear contradiction. Then $x \neq h_\ell h_mc_n$.

The final possibility to consider is
\[
x = h_ah_bh_ch_dh_e,
\]
where we take $a\leq b\leq c\leq d\leq e$ without loss of generality, by commutativity. In this case, we have
\[
t(x) = 2^a+2^b+2^c+2^d+2^e = 2^0+2^{j+3}+2^{j+4}+2^{j+6},
\]
meaning that $a=0$. Then
\[
2^0+2^b+2^c+2^d+2^e = 2^{j+3}+2^{j+4}+2^{j+6},
\]
meaning that two of the powers of 2 on the left-hand side combine to produce a single power of 2. The ordering $b\leq c\leq d\leq e$ can be used to reduce the number of possibilities for $(b,c,d,e)$:
\[
\begin{array}{c c c c c}
(j+2,j+2,j+4,j+6) &&(j+3,j+3,j+3,j+6) &&(j+3,j+4,j+5,j+5)
\end{array}
\]
Note that the first two tuples above produce the same element
\[
x = h_0h_{j+2}^2h_{j+4}h_{j+6} = h_0h_{j+3}^3h_{j+6},
\]
by the relation $h_{j+2}^2h_{j+4} = h_{j+3}^3$, while the third element produces a trivial element
\[
x = h_0h_{j+3}h_{j+4}h_{j+5}^2 = 0,
\]
by the relation $h_{j+3}h_{j+4} = 0$. As a result, the only elements of the form
\[
x = h_ah_bh_ch_dh_e
\]
with the proper internal degree are those of the form
\[
x = h_0h_{j+3}^3h_{j+6}.
\]

In summary, the only classes $x$ as in the statement of the proposition are linear combinations of the classes
\[
\begin{array}{c c c}
h_0f_{j+2} &&h_0h_{j+3}^3h_{j+6}
\end{array}
\]
for $j\geq0$.
\end{proof}

\begin{cor}
The class $h_1h_{j+4}Q_3(j)$ is defined and nonzero on the Adams $E_3$ page for all $j\geq0$.
\end{cor}
\begin{proof}
By \Cref{prop:identifysq0y0}, we may write
\[
    h_{j+4}Q_3(j) = (\text{Sq}^0)^j(y_0),
\]
so Bruner's formula gives
\[
    d_2(h_{j+4}Q_3(j)) = h_0(\text{Sq}^0)^j\text{Sq}^1(y_0)
\]
for all $j\geq0$. Since $h_1$ is a permanent cycle, the Leibniz rule gives
\[
    d_2(h_1h_{j+4}Q_3(j)) = h_1h_0(\text{Sq}^0)^j\text{Sq}^1(y_0) = 0
\]
for all $j\geq0$, by the relation $h_0h_1 = 0$. Then $h_1h_{j+4}Q_3(j)$ is a $d_2$ cycle. To see that it is not a $d_2$ boundary, we use the charts in \cite{charts} to settle the $j=0$ case and \Cref{prop:possible-d2} to consider all possible classes supporting such a $d_2$ differential for $j\geq1$. For first possible differential in \Cref{prop:possible-d2}, note that
\[
d_2(h_0f_{j+2}) = h_0d_2(f_{j+2}) = 0,
\]
by the Leibniz rule and the fact that $f_{j+2}$ is a $d_2$ cycle for all $j\geq1$. For the second possible differential in \Cref{prop:possible-d2}, we have that
\[
\begin{split}
d_2(h_0h_{j+3}^3h_{j+6}) &= h_0d_2(h_{j+3})h_{j+3}^2h_{j+6} + h_0h_{j+3}^3d_2(h_{j+6}) \\
&= h_0^2h_{j+2}^2h_{j+3}^2h_{j+6} + h_0^2h_{j+3}^3h_{j+5}^2 \\
&= 0,
\end{split}
\]
by the Leibniz rule, the Hopf invariant one differential $d_2(h_{j+1}) = h_0h_j^2$, and the relations $h_jh_{j+1} = 0$ and $h_jh_{j+2}^2 = 0$. Then the only classes which can possibly support $d_2$ differentials with target $h_1h_{j+4}Q_3(j)$ are themselves $d_2$ cycles. The claim follows.
\end{proof}

The arguments above show that
\[
d_3(f_{j+2}) = h_1h_{j+4}Q_3(j)
\]
for all $j\geq 1$ and that the target $h_1h_{j+4}Q_3(j)$ is defined and nonzero on the Adams $E_3$ page. In other words, the only remaining case is $j=0$. This case is handled in Lemma 5.25 of \cite{IWX} by an argument relying on nearby differentials. Then we have shown that
\[
d_3(f_{j+2}) = h_1h_{j+4}Q_3(j) \neq 0
\]
for all $j\geq0$, as claimed.

\section{Appendix}\label{sec:appendix}

%Here we provide code snippets for the various automated tasks needed throughout the paper.

%\subsection{Enumerating Possible Classes}

%In this subsection, we provide a very simple Python script that can be used to enumerate classes in a known finite range of degrees, perhaps with modifications to match different use cases:

%\begin{lstlisting}
%def deg(a,b,c,d,e):
%    output = 0
%    for x in [a,b,c,d,e]:
%        output += 2**(x+1)
%    return output

%output = []
%for a in range(lo,hi+1):
%    for b in range(a,hi+1):
%        for c in range(b,hi+1):
%            for d in range(c,hi+1):
%                for e in range(lo,hi+1):
%                    if deg(a,b,c,d,e) == targetDeg and sum([a,b,c,d]) > -4:
%                        output.append([a,b,c,d,e])
%print(len(output),"\n")
%for x in output:
%    print(x,sum(x[:4]),"\n")
%\end{lstlisting}

%To use this script, set \lstinline{lo} to be the index of the smallest possible power of 2 considered, \lstinline{hi} to be the index of the highest possible power, and \lstinline{targetDeg} to be the degree under consideration. As an example, suppose that one is searching for all
%\[
%q_iq_jq_kq_l\cdot h_m
%\]
%which could participate in a Burklund-Xu differential with target $q_{j-1}^4\cdot h$. Since \textbf{(is this a good example?)}

%\[
%t(q_iq_j\cdot h_k) = 2^{i+1}+2^{j+1} - 2
%\]

%\subsection{Bruner and BFCC Resolvers}\label{sec:resolvers}

\subsection{Ext Resolver Software}

Here we use version 1.9.5 of Bob Bruner's software (\cite{BrunerExt}), available at \url{http://www.rrb.wayne.edu/papers/}, and the SSeq software developed and maintained by Beauvais-Feisthauer, Chatham, and Chua (\cite{BFCC24}), available at \url{https://github.com/SpectralSequences/sseq}. The reason for using both programs is that, while both compute resolutions needed for the Adams $E_2$ page for the sphere, the class being studied in this paper is known only as $6_{16}$ in Bruner's notation. At the same time, only the BFCC resolver currently has support for computing the algebraic Steenrod operation $\text{Sq}^0$.

After building Bruner's program, there will be a directory labeled \lstinline{A/S0}. In it, there is a chart of the first 104 stems through filtration 20, and in (stem,filtration)-bidegree (38,6) of that chart, one can find a class labeled $6_{16}$.

After building the BFCC software, one can navigate to the directory labeled \lstinline{ext}. If one runs the command \lstinline{cargo run --example bruner}, the user will be prompted to specify a maximal stem $n$ through which to perform a translation from the BFCC class names to the Bruner class names. The user need not specify a save directory since the output will be printed to the screen. If one uses 38 as the maximal stem $n$, then the following two lines will appear in the output:

\begin{lstlisting}
x_(38,6,0) = [1,1]
x_(38,6,1) = [0,1]
\end{lstlisting}

Since Bruner's program produces two classes $6_{16}$ and $6_{17}$ in the same bidegree, these two output lines indicate that the BFCC class $x_{38,6,0}$ corresponds to the Bruner class $6_{16}+6_{17}$ while the BFCC class $x_{38,6,1}$ corresponds to the Bruner class $6_{17}$. As a result, the BFCC class $x_{38,6,0}+x_{38,6,1}$ corresponds to the Bruner class $6_{16}$.

Now that the Bruner class $6_{16}$ has been identified with the BFCC class $x_{38,6,0}+x_{38,6,1}$, we may compute $\text{Sq}^0$ on this class. To do so, one may stay in the \lstinline{ext} directory and run \lstinline{cargo run --example sq0}. The user will be prompted to specify a module, with the default module being that for the sphere. Next, the user will be prompted to specify a save directory if desired, but again the output will be printed to the screen all the same. For the maximal $n$, one may input 80 since this will compute a sufficiently large range for our purposes, and the filtration $s$ can be left at the default value of 7. The following two lines will appear in the output:

\begin{lstlisting}
Sq^0 x_(38,6,0) = [0,1]
Sq^0 x_(38,6,1) = [1,0]
\end{lstlisting}

Since $\text{Sq}^0$ is linear, these two lines indicate that
\[
\text{Sq}^0(x_{38,6,0}+x_{38,6,1}) = x_{82,6,1}+x_{82,6,0}
\]
since $\text{Sq}^0$ doubles internal degree and preserves filtration.

Now it remains to translate this linear combination back into Bruner's notation so that we may use another feature of Bruner's software to determine the identity of this class. Again running \lstinline{cargo run --example bruner}, this time through maximal stem $n = 82$, the following two lines appear in the output:

\begin{lstlisting}
x_(82,6,0) = [0,1]
x_(82,6,1) = [1,1]
\end{lstlisting}

Just as above, these two lines indicate that the sum of these two classes corresponds to Bruner's $6_{63}$. In other words,
\[
\text{Sq}^0(6_{16}) = 6_{63}.
\]

In order to see some of the $\text{Ext}_\mathcal{A}(\mathbb{F}_2,\mathbb{F}_2)$-module structure, one can, again in the directory \lstinline{S0}, open the file named \lstinline{himults}. If the class $5_{39}$ were to support an $h_3$-multiple, then this file would contain a line reading
\begin{lstlisting}
s g 5 39 3
\end{lstlisting}
for some \lstinline{s} and \lstinline{g}. Since there is no such line in \lstinline{himults}, we conclude that the class $5_{39}$ does not support multiplication by $h_3$ and therefore must be $Q_3(0)$. In order to identify $h_4Q_3(0)$, we look for this multiplication in \lstinline{himults}. The only line of the form
\begin{lstlisting}
s g 5 39 4
\end{lstlisting}
has $s = 6$ and $g = 63$, indicating that $6_{63}$ is the only nonzero term in $h_4Q_3(0)$. Therefore, the class $6_{63}$ must be $h_4Q_3(0)$.

\subsection{Massey Products}

Here we present a proof of those Massey products needed beyond what can be deduced from the calculations of \cite{IWX}. Namely, we show that $c_i\in\langle h_i,h_{i+1},h_{i+2},h_{i+1}^2\rangle$ and $V_i\in\langle h_i,h_{i+1},h_{i+2},h_{i+3},p_{i+1}'\rangle$ for all $i\geq0$.

First note that Table 3 of \cite{IWX} gives a fourfold Massey product representation
\[
    c_0 = \langle h_1^2,h_0,h_1,h_2\rangle
\]
of $c_0$. Although cyclic permutation of factors does not generally preserve Massey products, it does so in this case so that
\[
    c_0 \in \langle h_0,h_1,h_2,h_1^2\rangle.
\]
To see this, recall that $\text{Ext}_\mathcal{A}$ can be identified with the cohomology of the Lambda algebra $\Lambda$ (see section 3.3 of \cite{Ravenel-green}). The factors in the Massey product above have the following cohomology representatives:
\[
\begin{array}{c c c}
    h_0 = [\lambda_0] &h_1 = [\lambda_1] &h_2 = [\lambda_3]
\end{array}
\]
The first two consecutive products in $\langle h_0,h_1,h_2,h_1^2\rangle$ are trivial in $\Lambda$:
\[
    \lambda_0\lambda_1 = \lambda_1\lambda_3 = 0
\]
Then we define $a_{01} = a_{12} = 0$ in the first layer of our defining system. For the final consecutive product
\[
    \lambda_3\lambda_1^2
\]
in the Massey product, we define $a_{23} = \lambda_5\lambda_1$ so that
\[
    d(a_{23}) = d(\lambda_5\lambda_1) = \lambda_3\lambda_1^2.
\]
To define the next layer of our defining system for $\langle h_0,h_1,h_2,h_1^2\rangle$, we note that
\[
    \lambda_0 a_{12} + a_{01}\lambda_3 = 0,
\]
so we may choose $a_{02} = 0$. For the other threefold subbracket, we note that
\[
    \lambda_1 a_{23} + a_{12}\lambda_1^2 = \lambda_1\lambda_5\lambda_1 + 0 = \lambda_1\lambda_5\lambda_1.
\]
Defining $a_{13} = \lambda_3\lambda_5$, we have
\[
    d(a_{13}) = d(\lambda_3\lambda_5) = \lambda_3^2\lambda_1 = \lambda_1\lambda_5\lambda_1,
\]
using the relations in $\Lambda$. We now have all the elements of our defining system, so we compute
\[
    \begin{split}
        \langle h_0,h_1,h_2,h_1^2\rangle &\ni \lambda_0 a_{13} + a_{01}a_{23} + a_{02}\lambda_1^2 \\
        &= \lambda_0\lambda_3\lambda_5 + 0 + 0 \\
        &= \lambda_0\lambda_3\lambda_5 \\
        &= \lambda_2\lambda_3^2.
    \end{split}
\]
The $\Lambda$ representative for $c_0$ listed in (1.1) of \cite{Chen}, before Theorem 1.2, is precisely this element, giving
\[
    c_0 \in \langle h_0,h_1,h_2,h_1^2\rangle,
\]
as claimed. As a result, we have
\[
    c_i \in \langle h_i,h_{i+1},h_{i+2},h_{i+1}^2\rangle
\]
for all $i\geq0$. We may then shuffle to obtain
\[
    h_1h_{j+3}c_{j+2} \in \langle h_{j+2},h_{j+3},h_{j+4},h_1h_{j+3}^3\rangle
\]
and
\[
    h_kh_{j+6}c_{j+1} \in \langle h_{j+1},h_{j+2},h_{j+3},h_kh_{j+2}^2h_{j+6}\rangle
\]
when $j\geq2$.

We now show that
\[
    \langle h_i,h_{i+1},h_{i+2},h_{i+3},p_{i+1}'\rangle \ni V_i
\]
by first focusing on the $i=0$ case. Note, as above, that $h_i$ has $\Lambda$ algebra representative $\lambda_{2^i-1}$. Furthermore, these classes satisfy
\[
    \lambda_0\lambda_1 = \lambda_1\lambda_3 = \lambda_3\lambda_7 = 0
\]
before passing to $H^\ast\Lambda$, so we may choose
\[
    a_{01} = a_{12} = a_{23} = 0
\]
in a defining system for the Massey product in question. It takes slightly more work to find $a_{34}$, but one can check that
\[
    a_{34} = \lambda_7\lambda_{17}\lambda_{95}\lambda_{31} + \lambda_7\lambda_{17}\lambda_{39}\lambda_{87} + \lambda_1\lambda_{15}\lambda_{39}\lambda_{95}
\]
satisfies
\[
    da_{34} = \lambda_7\lambda_1(\lambda_{79}\lambda_{31}^2+\lambda_{31}^2\lambda_{79}),
\]
so we use this definition of $a_{34}$ in the defining system.

For the next layer of the defining system, we may choose
\[
    a_{02} = a_{13} = 0
\]
since $a_{01}$, $a_{12}$, and $a_{23}$ were trivial. One may check that
\[
    a_{24} = \lambda_{33}\lambda_{35}\lambda_{39}\lambda_{47}
\]
satisfies
\[
    da_{24} = \lambda_3a_{34},
\]
so we use this definition of $a_{24}$ in the defining system. We now need to choose nullhomotopies $a_{14}$ and $a_{03}$ of the fourfold subproduct representatives.

For $a_{14}$, this means finding a nullhomotopy of
\[
    a_1a_{24}+a_{12}a_{34}+a_{13}a_4 = \lambda_1\lambda_{33}\lambda_{35}\lambda_{39}\lambda_{47}+0+0 = \lambda_1\lambda_{33}\lambda_{35}\lambda_{39}\lambda_{47}.
\]
One can check that
\[
    a_{14} = \lambda_{31}\lambda_{47}^2\lambda_{31}
\]
is such a nullhomotopy.

Similarly, for $a_{03}$, we seek a nullhomotopy of
\[
    a_0a_{13}+a_{01}a_{23}+a_{02}a_3 = 0+0+0 = 0.
\]
Then we may once again choose the trivial nullhomotopy
\[
    a_{03} = 0
\]
in the defining system.

We now have a complete defining system with
\[
\begin{split}
    \langle h_0,h_1,h_2,h_3,p_1'\rangle &\ni a_0a_{14}+a_{01}a_{24}+a_{02}a_{34}+a_{03}a_4 \\
    &= \lambda_0\lambda_{31}\lambda_{47}^2\lambda_{31}+0+0+0 \\
    &= \lambda_0\lambda_{31}\lambda_{47}^2\lambda_{31}.
\end{split}
\]
Examining the presentation of the Adams 5-line in \cite{Chen}, we see that this is a Lambda algebra representative for $V_0$. Indeed, we may apply the relation
\[
    \lambda_{15}\lambda_{63} = \lambda_{47}\lambda_{31}
\]
to rewrite Chen's
\[
    V_0 = [\lambda_0\lambda_{31}\lambda_{47}\lambda_{15}\lambda_{63}]
\]
as
\[
    V_0 = [\lambda_0\lambda_{31}\lambda_{47}^2\lambda_{31}].
\]

Applying $\text{Sq}^0$ then gives
\[
    \langle h_i,h_{i+1},h_{i+2},h_{i+3},p_{i+1}'\rangle \ni V_i.
\]

\printbibliography

\end{document}